\documentstyle[12pt]{article}
\begin{document}
\title{The fiber of the Griffiths map for
the non-hyperelliptic Fano threefolds of genus 6}
\author{Atanas Iliev}
\date{ }
\maketitle
\begin{abstract}
Among smooth non-rational Fano $3$-folds, the non-hyperelliptic
Fano $3$-fold $X = X_{10}$ of degree $10$ (= of genus $6$)
has the unique property to admit
a non-trivial orbit of birationally equivalent threefolds,
inside its moduli space ${\cal X}_{10}$.
Here we prove that these orbits are,
in fact, the same as the fibers of the Griffiths intermediate
jacobian map on ${\cal X}_{10}$.  This coincidence leads up to
the main result of the paper:  The general fiber of the
Griffiths map on ${\cal X}_{10}$ is a union of two irreducible
$2$-dimensional components $F_1$ and $F_2$ of birationally
equivalent threefolds.  Moreover, each $F_i$
is isomorphic to the family of conics on any threefold
$X \in F_i$,  and any threefold of $F_i$ can be obtained
from a fixed $X(0) \in F_i$ by a ``quadruple cubic projection
from a conic on $X(0)$.  As an application, we give a negative
answer to the Tjurin's conjecture: {\it The general $X_{10}$
is birational to a quartic double solid}.
\end{abstract}

\bigskip

{\sc Contents.}

\smallskip

I. \ \  Introduction $....................................
.........................................................$   1

A. \ The prehomogeneous Fano fourfold $W$, and the
   moduli space ${\cal X}_{10}$ $..........$   8

L. \ Line-transformations $...............................
................................................$  16

C. \ Conic-transformations $..............................
...............................................$  26

N. \ Nodal $X_{10}$ $.....................................
..........................................................$  32

F. \ The fiber of the Griffiths intermediate jacobian map
   $j:{\cal X}_{10} \rightarrow {\cal A}_{10}$ $.......$  43

References $.............................................
........................................................$  50


\bigskip

\centerline{{\bf I.  Introduction.}}

{\bf (I.0)}
{\it The unhandy Fano threefold $X_{10}$.}

\ According to the classification
[Isk1] - [Isk2]
of Fano-Iskovskikh,
there are $18$ types of Fano threefolds of the
$1^{-st}$ species ( = $rk.{\bf Pic} = 1$).
For $8$ of these types:
${\bf P}^3$, the quadric $Q_2$,
the complete intersections of two quadrics $Y_4$,
the del Pezzo threefold $Y_5 = {\bf P}(T_{{\bf P}^2}) \cong$
any smooth complete intersection of the Grassmannian
$G = G(2,5) \subset {\bf P}^9$ and a subspace ${\bf P}^6$,
and the threefolds
$X_{2g-2} \subset {\bf P}^{g+1},
g = 7,9,10,12$  (as well the exceptional
$X'_{22}$ - see [MU], [F]),
any threefold is rational.

The general element of the next $10$ families
is non-rational (see e.g. [B1]);
and there are two basic methods to their study:

{\bf (CG).} \ The Clemens-Griffiths approach -- studying the
relation between the Fano threefold $X$ and its
principally polarized intermediate jacobian
$(J(X),{\Theta})$, and

{\bf (FI).} \ The approach of Fano-Iskovskikh -- studying the
birational isomorphisms $X \rightarrow Y$,
in the class of smooth Fano $3$-folds.

Three -- among the $10$ ``non-rational'' types --
the cubic $Y_3$, the quartic double solid $Y_2$,
and the complete intersection of two quadrics $X_8$,
have relatively simple jacobians, and
the approach (CH) has been successfully
applied for them, especially -- to prove the Torelli
theorem for each of them, in the form:
The singular theta divisor ${\Theta}$ determines unique the
threefold (see e.g. [CG], [Tju],
[B2], [Vo], [C2], [De1], [D2]).
Moreover, the threefold $X_{14}$ is birational to
a cubic $Y_3$, and its properties can be reduced to
properties of $Y_3$ (see e.g. [Isk3, Ch.3],
[Tr], [Pu]).

As far as the author is informed,  the approach (FI),
which works -- in general --
for threefolds which are ``far from ${\bf P}^3$'',
was successful to describe the birational
isomorphisms for $5$ of the rest $6$ types:
the quartic $X_4$,
the sextic double solid $X_2$,
the double Veronese cone $Y_1$,
the complete intersection of a quadric and a cubic $X_6$,
and
the double quadric $X'_4$.

As a summary, any birational map from one of these
threefolds $X$ to a Fano $3$-fold of the $1^{-st}$ specie
must be either a biregular automorphism of the
threefold $X$, or -- at most -- a composition of
birational involutions of $X$ and a biregular
automorphism
(see [IM], [I3], [I5], [H1]).

What was left is the type of Fano threefolds of
degree $10$, or -- equivalently -- {\it the Fano
threefolds of genus $6$} = the genus
of the canonical curve on the threefold
(see [Isk1], [Isk2]).
There are two kinds of such threefolds:

{\it The non-hyperelliptic kind}
$X_{10} = G \cap {\bf P}^7 \cap Q$ =
the complete intersection of the
Grassmannian $G = G(2,5) \subset {\bf P}^9$
with a subspace ${\bf P}^7$ and a quadric $Q$, and

{\it The hyperelliptic kind} -- the Gushel threefold
$X'_{10} \rightarrow Y_5$ =  the double covering
of the Fano threefold $Y_5 = G \cap {\bf P}^6$,
branched along a quadratic section of $Y_5$.
Moreover, any smooth $X'_{10}$ is a smooth
projective deformation of $X_{10}$,
i.e. $X'_{10}$ is ``less general'' than $X_{10}$
(see [Gus]).

{\bf (I.1)}
{\it The existence of a non-trivial birational orbit
-- a unique property of $X_{10}$.}

\ In this paper we study the ``more general''
non-hyperelliptic Fano threefold $X_{10}$.

First, we shall try to answer the question why
it is hard to apply for $X_{10}$ either of the
approaches (CG) and (FI).

The approach (CG), being applicable
for Fano's with ``simple'' jacobians should work
for Fano's $X$ of non-rational types and with small
$h^{1,2} = dim \ J(X)$.  Indeed,
$h^{1,2}(Y_3) = 5$, and $h^{1,2}(Y_2) = 10$,
$h^{1,2}(X_8) = 14$ -- which is conspicuously less than
e.g. $h^{1,2}(X_4) = 30$.  Moreover, $J(Y_3)$ and
$J(X_8)$ are Prym varieties -- which, in particular,
simplifies their descriptions,  while studying the
$10$-dimensional $J(Y_2)$ had taken lot of efforts
([W], [Ti1], [Ti2], [Ti3], [C1],
[Vo], [D2], [C2]).

Concerning $X_2$, $X_4$, $X'_4$, $Y_1$, and $X_6$ --
the field of the approach (FI) -- one has:
$h^{1,2}(X_2) = 52$,
$h^{1,2}(X_4) = H^{1,2}(X'_4) = 30$,
$h^{1,2}(Y_1) = 21$, and $h^{1,2}(X_6) = 20$;
and one should be prevented from the large
$h^{1,2}$ to find an easy way to apply (CG).

As regards $X_{10}$ and $X'_{10}$, one has
$h^{1,2}(X_{10}) = h^{1,2}(X'_{10}) = 10$,
and the Fano threefolds of degree $10$
should be nearer to such threefolds like
the quartic double solid $Y_2$
(see (I.4)).

However, as it follows from the forthcoming,
the Torelli theorem does not hold for $X_{10}$
-- i.e. one cannot recognize $X_{10}$ from
its p.p. intermediate jacobian
$J(X_{10}, {\Theta})$.
In other words,  if
${\cal X}_{10}$ is the moduli space
of $X_{10}$,  then the
Griffiths intermediate
jacobian map
$\ j:{\cal X}_{10} \rightarrow {\cal A}_{10} \ $
(see e.g. [CG], [Tju])
has a non-trivial fiber
$j^{-1}({\xi})$, over the general
${\xi} =(J(X_{10}),{\Theta}) \in {\cal J}_{10} = Image(j)$.

In order to answer the question:  Why $X_{10}$ is far
from threefolds like the quartic and the sextic double
solid, we introduce the notion  {\it birational orbit}:

{\bf (*).}  Let $X$ be a Fano threefold
of the $1^{-st}$ specie, and of any
of the $10$ ``non-rational'' types (say $X_4$,...).
Define
{\it the birational orbit of $X$} :=:
$Orb_{bir}(X)$ to be the set of
all the Fano threefolds $Y$ -- of the same type
(i.e. $X_4$,...) -- which are birational to $X$:
mod.(biregular equivalence).

Turning to the ``(FI)-threefolds'':
$X_4,X'_4,X_6,X_2$, and $Y_1$,
one immediately concludes that the birational
orbit of the general element of any of these
types is trivial --
any such birational isomorphism must transform $X$
into itself.

Even for the ``(CG)-threefolds'':
$Y_3$, $Y_2$, and $X_8$ --
which, in particular, have large groups of birational
automorphisms -- one can immediately see that the
Torelli theorem, together with the generic non-rationality,
prohibit the non-triviality of the birational orbit
of their general elements $X$.  Indeed, let $X$
be a general element of any of these $3$ types,
and let $X' \in Orb_{bir}(X)$; in particular,
$X'$ is also non-rational.
As it is well-known,
the jacobian of the general $3$-fold $X$ of
any of these $3$
types is an indecomposable abelian variety,
which is  not an jacobian of a curve.
Therefore the birationality of $X$ and $X'$ implies
that $J(X')$ is a product of $J(X)$ and an jacobian
$A$ of a (possibly -- non-connected) curve.  Now,
the identity
$dim \ J(X) = dim \ J(X')$ implies that $A = 0$.
That is $J(X) = J(X') \ \Rightarrow X = X'$ --
by the Torelli theorem.

All this brings up the question:

{\it Is $Orb_{bir}$ a trivial notion?} -- No, it isn't:
$Orb_{bir}(X_{10})$ is not trivial; and this
non-triviality is the property, by which $X_{10}$
differs from the rest non-rational types of
Fano threefolds.  Here, among all, we shall find
the birational orbit of the general $X_{10}$.

{\bf (I.2)}
{\it Structure of the paper.}

Let $X = X_{10} \in {\cal X}_{10}$ be a general
non-hyperelliptic Fano threefold of degree $10$.
Let
${\xi} = j(X) = (J(X), {\Theta})$
be its intermediate jacobian, and let
$j^{-1}({\xi}) \subset {\cal X}_{10}$
be the fiber of $j$ through $X$.

As one may expect, the inclusion
$Orb_{bir}(X) \subset j^{-1}(j(X))$
takes place from trivial reasons (see above).
The main result
of this paper -- Theorem (F.6) --
tells, in particular, that these two sets coincide.
That is -- if the general Fano threefolds $X,Y$
have the same jacobian, then they are birational
to each other.

The proof of the coincidence
$\{ \mbox{ the fiber of } j \mbox{ through }X \}$ =
$\{ \mbox{ the birational orbit of } X \}$,
is carried out in two steps:

{\bf Step (i).}  Studying the
fiber of the intermediate jacobian map $j$,
defined on the moduli space of $X_{10}$;
and this yields
the topological structure of the fiber -- it is
a union of two irreducible surfaces
$F \cup \overline{F}$;

{\bf Step (ii).}
Discovering of two $2$-dimensional irreducible
surfaces, inside the
birational orbit of a fixed general $X = X_{10}$;
and this, together with the inclusion
$Orb_{bir}(X) \subset$
$\{ \mbox{ the fiber } F \cup \overline{F}
\mbox{ through } X \}$,
yields
the algebraic structure of the components
$F$ and $\overline{F}$ -- each of them is
birational to the Fano surface of conics on any
of its elements.

In order to carry out step {\bf (i)},
in section A we describe first an appropriate model
of ${\cal X}_{10}$, as follows:
The general $X_{10}$ is a quadratic section
of a smooth $4$-fold of degree $5$
$W = W({\bf P}^7) = G \cap {\bf P}^7 \subset {\bf P}^7$.
It is known that such $W$ is unique, and -- as it
turns out -- $W$ is a prehomogeneous space for the
action of its $8$-dimensional automorphism group
$Aut(W)$ = an extension of ${\bf PGL}(2)$ by
a semidirect product
${\sf G}_m \triangleright {\sf G}_a^4$
(see (A.4.3)-(A.4.4)).  The action of $Aut(W)$ cuts the
space $\mid {\cal O}_W (2) \mid \cong {\bf P}^{30}$
into $8$-dimensional orbits of projectively equivalent
Fano threefolds of degree $10$;
and we see that the obtained $22$-dimensional
orbifold if birational to ${\cal X}_{22}$
-- see (A.6.2).

In order to find $dim \ j({\cal X}_{10})$,
and the topological structure of the fiber of $j$,
we study, in section N, the nodal
$X_{10}$, and its projection $Z_8 \subset {\bf P}^6$
through the node -- which is a degenerate intersection
of $3$ quadrics which contains a quadratic surface.
The observation is that the
jacobian of the general nodal $X_{10}$ is
a Prym variety of a double covering of a general plane
sextic.  It follows that the codimension one subspace
${\cal J}_{10}^{nodal} \subset {\cal J}_{10} = j({\cal X}_{10})$,
swept out by the jacobians of the
nodal $X_{10}$, is the same as the $19$-dimensional
space  ${\cal P}_6$  of Prymians of the
double coverings of plane sextics
-- see (N.3).

Now, we degenerate the result
of Friedman-Smith [FS] -- declaring
an (1:1) correspondence between general
complete
intersections $Z_8$ of $3$ quadrics,
and the corresponding double coverings of their
determinantal septics --
to general
$Z_8$ which contain quadratic surface,
and the corresponding double coverings
of their discriminantal sextics
(see section N).
As a result, we obtain that
the component $F$ of the fiber of $j$, through the nodal $X$,
is two-dimensional.
Moreover, the result of A.Verra
[Ve], which shows that the degree of the Prym map
for plane sextics is ${\it two}$,  implies that the fiber
of $j$ through the nodal $X_{10}$ has two components
-- both of dimension $2$.

The same is true also for the general fiber of $j$
-- it consists of two irreducible surfaces.
In particular,
$dim \ {\cal J}_{10} = dim \ {\cal X}_{10} - 2 = 20$,
where ${\cal J}_{10} \subset {\cal A}_{10}$
is the image of $j$ -- see above.

The next is to carry out step {\bf (ii).},
i.e., to consider the fiber of $j$ as
a birational orbit of threefolds,
and to describe this orbit as an
algebraic variety.  This is realized in sections
L, C, and F, as follows:

Let
$j^{-1}({\xi})$, ${\xi} = (J,{\Theta}) \in {\cal J}_{10}$,
be the general fiber of $j$, and let
$j^{-1}({\xi}) = F \cup \overline{F}$,
where $F$ and $\overline{F}$ are its two components.

Now, the question is to find an appropriate
$2$-dimensional family of $X_{10}$, which are birational
to each other. The base of this family will describe
one of these two components, say $F$,
of the fiber of the Griffiths map $j$.

As a summary, the main result
of section C, and of section F,
regarding this last question is:  Such a base is the
Fano surface ${\cal F}(X)$ of conics on any threefold
$X = X_{10}$ on the fiber of $j$.
Moreover, if
$X \in \mbox{ the component } F$
then ${\cal F}(X)$ and $F$ are birational to each other
-- see (F.5.2)-(F.5.3).

In order to prove this, we show in sections C and F
that the general conic $q$, on the general
$X \in {\cal X}_{10}$
defines a birational map to another element of
${\cal X}_{10}$:
${\alpha}_q: X \rightarrow X_q$.
Moreover, if $p$ and $q$ are two such non-coincident
conics on $X$, then $X_p \neq X_q$
-- see Proposition (F.2).

This implies that the Fano surface of any element
$X \in j^{-1}({\xi})$ is birational to one of the
components, say $F$, of $j^{-1}({\xi})$, and we show that
$X \in F$ \ -- see Proposition (F.5).
The residue component $\overline{F}$ of
$j^{-1}({\xi}) = F + \overline{F}$, is isomorphic to the
Fano surface of any line transform $X_l$ of $X$
-- by (F.5) and (F.5.2).
In particular, all the elements
of $j^{-1}({\xi})$ are birational to each other
(see sections L and C).

{\bf (I.3)}
{\it The line and the conic transformations
of $X_{10}$ -- analogs of the classical
double projections.}

Let $X = X_{10}$ be general, let
${\cal F}(X)$ be the Fano surface of conics
$q \subset X$, let
${\Gamma}(X)$ be the curve of lines $l \subset X$,
and let $F$ be the component of the fiber of
$j$ through $X$.

It is right to mention that the conic
transformations${\alpha}_q : X \rightarrow X_q$
-- which describe the component $F$,
as well the line transformations
${\beta}_l$
-- which interchange $F$ and $\overline{F}$,
relate the birational properties of the
threefolds $X_{10}$ and the birational
properties of the rational types of
Fano threefolds.

More precisely, the map
${\beta}_l:X \rightarrow X_l$
is defined by the system
$\mid {\cal O}_X (2) -3.l \mid$,
while the map
${\alpha}_q:X \rightarrow X_q$
is defined by the system
$\mid {\cal O}_X (3) - 4.q \mid$.
The point is that these maps
are {\it birational isomorphisms}
{\it between non-biregular}
Fano threefolds.

As it follows from (I.0) and (I.1),
an existence of such birationalities,
between non-biregular threefolds,
never happens, in general,  for non-rational
types of Fano threefolds -- except, of course,
these which come from the well known
birationalities between the cubics
and the threefolds $X_{14}$
(see also (I.5) -- Remark (a.i)).
Turning to rational classes of Fano
threefolds, one can find the three classical
birational double projections from a line:

$X_{16} \rightarrow {\bf P}^3$,
$X_{18} \rightarrow {Q}_2$,
and
$X_{22} \rightarrow Y_5$.

The double projections, discovered
in its own time by G.Fano, and used
-- in particular -- to establish
the rationality of
$X_{16},X_{18}$, and $X_{22}$,
have been reexamined recently
as generating examples of $D$-flops
over the projections from a line
(see e.g. [Isk7]).

The same is the situation with the
``triple quadratic projection''
from a line $l \subset X_{10}$,
and with the
``quadruple cubic projection''
from a conic $q \subset X_{10}$:
There exists a special $D$-flop
$X' \rightarrow X^+$,
which resolves, mod. codimension $2$,
the projection of $X_{10}$
from the line $l$ (or -- from
the conic $q$) -- see sections L and C.

{\bf (I.4)}
{\it A negative answer to a Tjurin's conjecture.}

In addition,  it is worth to mention that
the results of the paper imply, in particular,
an answer of the Tjurin's question
(see [Tju, Ch.3:Sect.2]):

{\bf (**).}
{\sl
Is the variety $X_{10}$ birational to
a quartic double solid $Y_2$?
}

The answer is {\it no} -- in general:
The image ${\cal J}_{10} = j({\cal X}_{10})$
is $20$-dimensional -- see e.g. Theorem (F.6),
while the well-known image
${\cal J}_2 = j(\{ \mbox{Quartic double solids} \})$
is $19$-dimensional -- see e.g [C1], [D2], [Ve].

\newpage

{\bf (I.5)}
{\it Remarks.}

{\bf (a).}
As it follows from the preceding, the only two other
non-rational types of Fano's of the $1^{-st}$ specie,
for which the birational orbit of the general element
might be non-trivial, are ${\cal X}_{14}$, and
${\cal X}'_{10}$.
Here we discuss separately either of these two cases:

{\bf (i).}
{\it The birational orbit of $X_{14}$.}

As it is well-known -- see e.g. [Gus], [Isk6] --
the Fano $3$-fold $X_{14}$ is an intersection of
$G(2,6) \subset {\bf P}^{14}$ and a $codim. \ 5$ subspace
${\bf P}^9$.
That is, the moduli space
${\cal X}_{14}$
should be birational to the $15$-dimensional quotient
space $G(5,15)/Aut \ G(2,6)$.

Here we present one conjecturable decription of the
fiber of the Griffiths map on ${\cal X}_{10}$:

Let ${\cal Y}_3$ be the $10$-dimensional moduli space
of cubic threefolds.  As one can see,
there exists a generically $1:1$-correspondence $\Sigma$
between the $20$-dimensional spaces
${\cal C}^1_5({\cal X}_{14})$ =
$\{ (X_{14},C) : C \subset X_{14}, deg(C) = 5, p_a(C) =1 \}$,
and
${\cal C}^1_5({\cal Y}_3)$ =
$\{ (Y_3,C) : C \subset Y_3, deg(C) = 5, p_a(C) =1 \}$;
and this correspondence $\Sigma$ conforms to birational
isomorphisms
(see e.g. [Isk3, Ch.3], [Pu]).

Let
$p':{\Sigma} \rightarrow {\cal C}^1_5({\cal X}_{14})$,
$p'':{\Sigma} \rightarrow {\cal C}^1_5({\cal Y}_3)$,
$p_{14}:{\cal C}^1_5({\cal X}_{14}) \rightarrow {\cal X}_{14}$,
and
$p_3:{\cal C}^1_5({\cal Y}_3) \rightarrow {\cal Y}_3$
be the natural projections.
Then, by our conjecture, the $5$-dimensional fiber of
$j:{\cal X}_{14} \rightarrow j({\cal Y}_3)$,
above the general $(J,{\Theta}) = j(Y_3)$,
has a component of maximal dimension,
which is birational to
$p_{14} \circ (p'')^{-1} \circ (p_3)^{-1}(Y_3)$;
and this fiber coincides with the birational orbit,
in ${\cal X}_{14}$, of any of its elements.

Anyway, the description of the birational orbit of
$X_{14}$, is a problem of parameterization of
the birational isomorphisms from a fixed cubic
$Y_3$, to Fano's of type ${\cal X}_{14}$.
Note that, besides the birational isomorphisms
$Y_3 \rightarrow X_{14}(C)$ -- defined by elliptic
quintics $C = C^1_5 \subset Y_3$
(ibid.),
are known also birational isomorphisms
$Y_3 \rightarrow X_{14}(E)$ -- defined by elliptic
quartics $E \subset Y_3$
(see [Tr]).

{\bf (ii).}
{\it The birational orbit of $X'_{10}$.}

One may assume that, being a smooth deformation
of $X_{10}$, the Gushel threefold should admit
-- by similarity -- a $2$-component fiber,
consisting of two Fano surfaces.
However, such a deformation does not live
in the model ${\cal X}_{10}$, of the moduli
space of non-hyperelliptic Fano's of degree $10$,
considered in this paper (see section A):
The smooth deformation from non-hyperelliptic
$X_{10}$-s to Gushel $3$-folds lives in the
moduli space ${\cal X}_{10}^c$  -- of Fano's
of degree $10$, regarded as embedded in
a cone $C(W)$ over the fourfold
$W = W_5 \subset {\bf P}^7$
(see [Gus], and section A).
Anyway, ${\cal X}_{10}^c$ is naturally
birational to our model ${\cal X}_{10}$.

In particular, all the generic properties of
the elements of
${\cal X}_{10}^c$ must be the same as the generic
properties of the elements of  ${\cal X}_{10}$.
This suggests to expect a $2$-dimensional
fiber of the Griffiths map $j^c$ on
${\cal X}_{10}^c$ --  also on
some special subvarieties of ${\cal X}_{10}^c$.

As it follows from the definition of
${\cal X}'_{10} \subset {\cal X}_{10}^c$,
this space is birational the moduli space
of their branched loci, i.e.:
${\cal X}'_{10}$ is birational to
the moduli space ${\cal F}_6$
of these $K3^{-s}$ \ $S_{10} \subset {\bf P}^6$,
which are quadratic sections
of the del Pezzo $3$-fold $Y_5$.
As is shown by Mukai
(see [Mu]),
the general $K3$-surface
$S_{10} \subset {\bf P}^6$
can be embedded as a quadratic section of $Y_5$,
and -- therefore -- the space
${\cal F}_6$ is birational to
the $19$-dimensional quotient
$Symm^{12} \ {\bf C}^2 \oplus Symm^8 \ {\bf C}^2 /
{\bf PGL}(2)$.
In short, the projective moduli of
$S_{10} \subset Y_5$
are the same as the quotient of
$\mid {\cal O}_{Y_5} (2) \mid \cong {\bf P}^{22}$,
by the group $Aut(Y_5) \cong {\bf PGL}(2)$
(see e.g. [MU] --  or [F], [PS1], [PS2]).
In particular,
${\cal X}'_{10}$ is of codimension $3$
in ${\cal X}_{10}^c$.

Let $X'_{10}$ be a general Gushel $3$-fold.
On the one hand, as one may expect,
the general conic transformations,
as well the general line transformations,
should bring $X'_{10}$ outside
${\cal X}'_{10}$, but they must leave
the birational images of $X'_{10}$
inside ${\cal X}_{10}^c$.
It is natural to assume that these
orbits are $2$-dimensional, i.e. --
of the same dimension as the general
birational orbits.
On the other hand, as far as the author can judge,
the subvariety  ${\cal J}'_{10}$ of intermediate
jacobians of Gushel $3$-folds should be a divisor
in the $20$-fold ${\cal J}_{10}^c$ =
$\{ j(X):X \in {\cal X}_{10}^c \}$.
One way to see this is to degenerate
the general nodal $X_{10}$ to a general nodal
Gushel $3$-fold.
Then -- see section N -- the Pfaff line
becomes totally tangent to the determinantal
sextic of the threefold, and one has to see
that the $18$-fold, of Prym varieties over
general plane sextics
with totally tangent line, is the same as
the set of abelian parts of
jacobians of general nodal Gushel $3$-folds.
Therefore, the ``birational extension''
$Orb_{bir}({\cal X}'_{10})$
(of all birational images of Gushel's,
in ${\cal X}_{10}^c$)
should be also a divisor in
${\cal X}_{10}^c$, and
the general birational orbit
$Orb_{bir}(X'_{10}) \subset {\cal X}_{10}^c$
should intersect ${\cal X}'_{10}$
in a finite number, say $N \ge 1$, of points.

If all this is true, then the integer $N$
will be the same as the degree of the Griffiths
map on ${\cal X}'_{10}$.
\ The question is:
$N = ?$.

{\bf (b).}
{\it The example of Donagi.}

Turning to the examples of non-smooth Fano's, one
can find the excellent example, due to Donagi,
of a non-trivial fiber of the Griffiths map
(see e.g. [C3]):

{\sl The general fiber
$j^{-1}(J,{\Theta})$
of the Griffiths map
$j:{\cal DS'} \rightarrow {\cal A}_4$
for quartic double solids with $6$ nodes is
a $32$-sheeted covering of a cubic threefold
$Y_3$ defined uniquely by $(J,{\Theta})$}.

\bigskip

\centerline{* \ * \ *}

\bigskip

The author would like to thank
to V.Iskovskikh and V.Shokurov --
for introducing into the subject,
to A.Verra and R.Donagi --
for the numerous helpful conversations,
and to S.Mukai -- for the help in
finding the group $Aut(W_5)$.
He is much obliged to A.Tikhomirov,
D.Logachev, S.Hashin, and Ju.Prokhorov
-- for their results on Fano's, which
inspired him to write this paper.


\newpage

\centerline{
{\bf A.
The prehomogeneous Fano fourfold $W$,
}}

\centerline{
{\bf
and the moduli space ${\cal X}_{10}$.
}}

{\bf (A.0)}
Let
$0 \le k \le n$ be inetegers, and let
$G(k+1,n+1) = G(k+1,{\bf C}^{n+1})$ =
$G(k:n) = G(k:{\bf P}^n)$
be the Grassmannian of $k+1$-dimensional
subspaces
${\bf C}^{k+1} \subset {\bf C}^{n+1}$ =
the Grassmannian of projective $k$-spaces
${\bf P}^k = {\bf P}({\bf C}^{k+1}) \subset
{\bf P}({\bf C}^{n+1}) = {\bf P}^n$.

Denote by
${\wedge}^{k+1} \ {\bf P}^n$
the projective space
${\bf P}({\wedge}^{k+1} \ {\bf C}^{n+1})$,
$0 \le k \le n$.
Let
${{\bf C}^{n+1}}^*$
be the dual space of
${\bf C}^{n+1}$,
with respect to a pairing
$(.,.):{\bf C}^{n+1} \times {\bf C}^{n+1}
\rightarrow {\bf C}$;
and denote by
${{\bf P}^n}^*$
the space
${\bf P}({{\bf C}^{n+1}}^*)$.

Let
${\bf P}^0(0) \subset {\bf P}^1(0) \subset ...
\subset {\bf P}^n$
be a fixed flag in ${\bf P}^n = {\bf P}^n(0)$.
It is well-known
(see e.g. [GH, $\S 1.5$])
that the complex Schubert cycles
${\sigma}_{a_1,...,a_{k+1}}$ =
${\sigma}_{a_1,...,a_{k+1}}
({\bf P}^0(0),{\bf P}^1(0),...,{\bf P}^n)$ =
$\{ {\bf P}^k \subset {\bf P}^n:
dim \ {\bf P}^k \cap {\bf P}^{n-k+i-a_i}(0) \ge i-1,
i = 0,...,k+1 \}$
$\subset G(k:{\bf P}^n)$,
$n-k \ge a_1 \ge a_2 \ge ... \ge a_{k+1}$,
generate the integer homology
of $G(k:{\bf P}^n)$.
Moreover,
${\sigma}_{a_1,...,a_{k+1}} \subset
G(k:{\bf P}^n)$
is of real codimension
$2(a_1 + ... + a_{k+1})$;
and, in the homology ring
$H^*(G(k:{\bf P}^n),{\bf Z})$,
the cycles
${\sigma}_{a_1,..,a_{k+1}}$
intersect each other by the
formulae of Pieri and Giambelli
([GH, $\S 1.5$]).

We shall use these formulae,
as well the shorten notation
of a Schubert cycle -- in which
we shall omit the entry ${\bf P}^k(0)$,
if it imposes an empty condition
on ${\bf P}^k$, or if the cycle
is considered to be in general position.

Fix
${\bf P}^4$
and the Grassmannian
$G = G(1:{\bf P}^4)$.
The wedge-map
${\bf P}^1 \subset {\bf P}^1
\mapsto
{\wedge}^2 \ {\bf P}^1 \subset {\wedge}^2 \ {\bf P}^4$
defines the Pl\"ucker embedding
$pl:G \rightarrow {\bf P}^9 :=
G(0: {\wedge}^2 \ {\bf P}^4) =
{\wedge}^2 \ {\bf P}^4$.
The isomorphic image
$G \subset {\bf P}^9$
is a smooth $6$-fold with a hyperplane section
defined by the linear system
$\mid {\sigma}_{1,0} \mid$, and of degree
$5 = deg \ {\sigma}_{1,0}^6$.

For a given $7$-space
${\bf P}^7 \subset {\bf P}^9$,
define
$W({\bf P}^7) := G \cap {\bf P}^7$.
The variety $W({\bf P}^7)$ is a $4$-fold,
since $G$ does not contain subspaces
of dimension $\ge 4$;
moreover $deg \ W({\bf P}^7) = deg \ G = 5$.

Let
${\cal D} =
\{ {\bf P}^7 : W({\bf P}^7) \mbox{ is smooth } \}$.
It is well-known
(see e.g. [Pr])
that all the fourfolds
$W({\bf P}^7)$
are projectively equivalent to
each other; in fact, the same
assertion follows from the
uniqueness of the orbit representation
of any smooth $W({\bf P}^7)$ -- see below.

Therefore, the projective equivalency class
of these smooth fourfolds defines the ''unique''
smooth fourfold
$W = W_5 \subset {\bf P}^7$ --
the Fano $4$-fold of degree $5$
in ${\bf P}^7$.

{\bf (A.1)}
{\it The family of planes on the fourfold $W_5$.}

The family of planes in
$G(1:{\bf P}^4 = G \subset {\bf P}^9$
has two components:

(1). the component of ${\sigma}$-planes
     ${\sigma}_{3,1}(x,{\bf P}^3)$ =
     $\{ {\bf P}^1 :
     x \in {\bf P}^1 \subset {\bf P}^3 \}$
     ,
     where
     $(x \in {\bf P}^3) \in
     \mbox{ the flag variety } F(0:3:{\bf P}^4)$,
and

(2). the component of ${\rho}$-planes
     ${\sigma}_{2,2}({\bf P}^2)$ =
     $\{ {\bf P}^1 :
     {\bf P}^1 \subset {\bf P}^2 \}$
     ,
     where
     ${\bf P}^2 \in G(2:{\bf P}^4)$.

{\bf Lemma.}
{\sl
Let
$W = W_5 = G \cap {\bf P}^7 \subset {\bf P}^7$
be the Fano fourfold of degree $5$
in ${\bf P}^7$.
Then:

{\bf (i).}
There exists a unique ${\rho}$-plane
${\bf P}^2_o = {\sigma}_{2,2}({\bf P}^2(0)) =
{\wedge}^2 \ {\bf P}^2(0)$,
which lie on $W$;

{\bf (ii).}
Let
${\bf P}^2(0) \subset {\bf P}^4$
be the ''plane'' of the unique
${\rho}$-plane
${\bf P}^2_o \subset W$.
The family of ${\sigma}$-planes on $W$
form a smooth rational $1$-dimensional
family
$\{ {\bf P}^2_t =
{\sigma}_{3,1}(x(t),{\bf P}^3(t)) \}$.
The family of their centers
$\{ x(t) \in {\bf P}^4 \}$,
sweeps out a smooth conic
$q(0) \subset {\bf P}^2(0) \subset {\bf P}^4$,
and the family of their ''spaces''
$\{ {\bf P}^3(t) \subset {\bf P}^4 \}$
coincides with the pencil of
${\bf P}^3-s$, in ${\bf P}^4$,
through ${\bf P}^2(0)$.

{\bf (iii).}
Let
$q(0) \subset {\bf P}^2(0)$
be the conic of centers from (ii),
and let
$q_o \subset {\bf P}^2_o$
be the conic of tangent lines
to $q(0)$ -- the dual conic of $q(0)$.

Let
$x(t) \in q(0)$ be the center of the
${\sigma}$-plane
${\bf P}^2_t$,
let $l(t)$ be the tangent line to $q(0)$
in ${\bf P}^0(t)$,
let
$x_t = {\wedge}^2 \l(t) \in q_o$
be the corresponding point in
$q_o \subset {\bf P}^2_o \subset W$,
and let $l_t$ be the tangent line to
$q_o$ in $x_t$.
Then the
${\sigma}$-plane ${\bf P}^2_t =
{\sigma}_{3,1}(x(t),{\bf P}^3(t))$
intersects the unique ${\rho}$-plane
${\bf P}^2_o = {\wedge}^2 \ {\bf P}^2(0) \subset W$
along the line $l_t$;

{\bf (iv).}
Let
$Y_o = \cup \{ {\bf P}^2_t \} \subset W$
be the union of all the ${\sigma}$-planes on $W$.
Then $Y_o$ is a hyperplane section of $W$ --
the unique effective divisor of the non-complete
linear system
$\mid {\cal O}_{W}(2) - 2.{\bf P}^2_o \mid$.
}

{\bf Remark.}
We will not write down here a proof
(see e.g. [Pu])
of Lemma (A.1), since it will follow from
the forthcoming coordinate study of $W$ --
where one can point at to the
${\rho}$-plane, to the
${\sigma}$-planes on $W$,
and also -- to the special
hyperplane section $Y_o \subset W$
(see e.g. (A.2.4) and thereafter).

{\bf (A.2)}
{\it The coordinate representation of $W$.}

{\bf (A.2.1)}
Denote by
$<M> \ \subset {\bf P}^n$
the projective-linear span of the subset
$M \subset {\bf P}^n$.
Let
$\{ e_0,...,e_4 \}$
and
$\{ x_0,...,x_4 \}$
be a dual pair of bases of
${\bf C}^5$ and ${{\bf C}^5}^*$,
with respect to the pairing $(.,.)$;
i.e., in the notation (A.0),
$<e_0,...,e_4> = {\bf P}^4$,
and
$<x_0,...,x_4> = {{\bf P}^4}^*$.
We call $\{ e_i \}$ and $\{ x_i \}$
also bases of ${\bf P}^4$ and of
${{\bf P}^4}*$.

Then
$e_{ij} = e_i \wedge e_j = -e_j \wedge e_i = -e_{ji}$,
and
$\{ e_{ij} : 0 \le i < j \le 4 \}$
is a base
of ${\wedge}^2{\bf P}^4 = {\bf P}^9$;
the same is true for
$x_{ij} = x_i \wedge x_j$,
and for the base
$\{ x_{ij}:0 \le j < j \le 4 \}$,
of
${\wedge}^2{{\bf P}^4}^*$,
which is naturally dual to
$\{ e_{ij} : 0 \le i < j \le 4 \}$.

{\bf (A.2.2)}
{\bf Lemma.}
{\sl
Let $W = {\bf G} \cap {\bf P}^7$
be the Fano $4$-fold of degree $5$.
Then one can choose the base vectors
$(e_i,x_j)$ in such way that
${\bf P}^7 \subset {\bf P}^9$
can be represented as the common set of
zeros (the base space) of the pencil
of hyperplane sections
$\{ H(t) = t_0.H_0 + t_1.H_1 :
(t) = (t_0:t_1) \in {\bf P}^1 \}$,
where

$H_0 = x_{03} - x_{14}, H_{1} = x_{04} - x_{23}$.
}

{\bf Proof.}
Since all the smooth $W({\bf P}^7)$ are projectively
equivalent to each other, it is sufficient to prove
that the intersection
$G \cap (H_{0} = H_{1} = 0)$
is smooth.
{\bf q.e.d.}

{\bf (A.2.3)}
{\bf Remark.}
Let
$\ \hat{}:H^0({\wedge}^2{\bf P}^4, {\cal O}(1))
\rightarrow Skew({\bf C}^5,{{\bf C}^5}^*)$,
$H \mapsto \hat{H}$
be the natural linear isomorphism
between the space
$H^0({\wedge}^2{\bf P}^4, {\cal O}(1))$ =
$({\wedge}^2 \ {\bf C}^5)^*$,
of hyperplane sections of ${\wedge}^2{\bf P}^4$,
and the space
of skew-linear maps
$Skew({\bf C}^5,{{\bf C}^5}^*)$.

Then the base space of the pencil
$<H_0,H_1> \subset {\bf P}({\wedge}^2 \ {\bf C}^5)^*
= {\wedge}^2{{\bf P}^4}^*$
will intersect smoothly the Grassmannian $G$
{\it iff}
all the elements of the pencil
$<\hat{H_0},\hat{H_1}>$
are skew-symmetric maps of maximal rank (=4)
(see e.g. [Put]).
In coordinates $(e_i,x_j)$,
all the elements $h_{ij}(t)$,
of the skew-symmetric matrix
$\hat{H(t)} = (h_{ij}(t))$,
are zero, except
$h_{03} = -h_{30} = -h_{14} = h_{40} = t_0$,
$h_{04} = - h_{40} = h_{23} = - h_{32} = t_{1}$,
i.e.
$rank \ \hat{H(t)} = 4,
\forall (t) = (t_0:t_1) \in {\bf P}^1$.
Therefore $W = G \cap (H_0 = H_1 = 0)$
is smooth.

{\bf (A.2.4)}
{\it The ${\rho}$-plane ${\bf P}^2_o$,
and the conic of centers $q(0)$}.

In the notation of (A.2.3),
the projective kernels
$x(t) = {\bf P}(Ker \ \hat{H(t)}) =
(-t_0t_1:-t_1^2:t_0^2:0:0)$,
sweep-out a smooth conic $q(0)$, in the plane
${\bf P}^2(0) := <e_0,e_1,e_2> \ \subset {\bf P}^4$ =
${\bf P}^4(x_0:...:x_4)$.

By the special choice of $H_0$ and $H_1$, the space
${\bf P}^7 = (H_0 = H_1 = 0)$ =
$<e_{01}, e_{01}, e_{12},
 e_{03} + e_{14}, e_{04} + e_{23},
 e_{13}, e_{24}, e_{34}>$.
It is clear that the plane
${\bf P}^2_o = <e_{01},e_{02},e_{12}>$ =
${\wedge}^2 \ {\bf P}^2(0)$ =
${\wedge}^2 \ <e_0,e_1,e_2>$
is the unique ${\rho}$-plane on $W$
(see Lemma (A.1)(i)).

It is not hard to see that
the conic of kernels $q(0)$
is the same as the conic of centers
from Lemma (A.1)(ii).
Now, (A.1)(iii) can be verified directly;
and the hyperplane section $Y_o$, from (A.5)(iv),
is $Y_o = W \cap (x_{34} = 0)$ --
see also the proof of (A.5.1).


{\bf (A.3)}
{\it The group $Aut \ G(1:{\bf P}^4).$}

The group
$Aut \ G(k:{\bf P}^n)$
is well-known.
It is naturally isomorphic
either to
$Aut \ {\bf P}^n = {\bf PGL}(n+1)$ --
if $n+1 \neq 2(k+1)$,
or
$Aut \ {\bf P}^n$ is a normal subgroup of index $2$
in $Aut \ G(k:{\bf P}^n)$ -- if $n+1 = 2(k+1)$
(see e.g.
[M.J.Cowen, Proc.A.M.S.,106:1(1989),99-106],
[W.Kaup, Math.Z.,144(1975),75-96]
).

In fact,
$G(k:{\bf P}^n) = G(k+1,n+1) \cong$
the homogeneous space
$GL(n+1)/GL(k+1) \times GL(n-k)$,
and
either the automorphism group of $G(k+1,n+1)$
is isomorphic to
${\bf P}(GL(n+1)) = {\bf PGL}(n+1)$,
or $k+1 = n-1$ and the isomorphism -- between
the factors of the group at the quotient --
generate an additional involution in
$G(k+1,n+1)$
(see e.g. [W.L.Chow,Ann.of Math.,50(1949),32-67]).

Let, in particular,
$G = G(1:{\bf P}^4)$;
and let, as usual,
${\bf P}^4$ be the projectivization
of the vector space ${\bf C}^5$.

The natural isomorphism
${\wedge}^2:Aut \ {\bf P}^4 \rightarrow Aut \ G$
is the projectivization of the linear map
$A \mapsto {\wedge}^2 \ A$,
$A \in GL(5) = Aut({\bf C}^5)$,
defined on the decomposable bivectors
$u \wedge v$, in ${\bf C}^5$, by:
${\wedge}^2 \ A \ (u \wedge v) =
A(u) \wedge A(v)$.
Since $G$ coincides with the projectivization
of the variety of decomposable bivectors,
${\wedge}^2 \ A$ is an automorphism of $G$.
The statement is that these
${\wedge}^2 \ A$  are the automorphisms of $G$.

{\bf (A.4)}
{\it The group of projective automorphisms $Aut_o(W)$.}

We shall find the group $Aut_o(W)$ of these
biregular automorphisms of $W$ which are induced
by automorphisms of $G \supset W$.

{\bf (A.4.1)}
By definition,

$Aut_o(W)$ =
$\{ {\alpha} \in Aut \ {\bf P}^7:
\exists {\alpha}' \in Aut \ G = G(1:{\bf P}^4)
\mbox{ s.t. } {\alpha} = {\alpha}'_{\mid W} \}$.

Fix the bases $(e_i,x_j)$
as in lemma (A.2.2).
In particular,

{\bf (1).}
${\bf P}^7 = (x_{03} - x_{14} = x_{04} = x_{23} = 0$ =
$<v_1 := e_{01}, v_2 := e_{01}, v_3 := e_{12},
  v_4 := e_{03} + e_{14}, v_5 := e_{04} + e_{23},
  v_6 := e_{13}, v_7 := e_{24}, v_8 := e_{34}>$.

As it follows from (A.3),
$Aut_o(W)$ is naturally embedded in
$Aut \ {\bf P}^4 = {\bf PGL}(5)$,
and

{\bf (2).}
$Aut_o(W)$ =
$\{ A \in Aut \ {\bf P}^4 :
{\wedge}^2 \ A \ ({\bf P}^7) = {\bf P}^7 \}$ =

$\{ A = (a_{ij}) \in {\bf PGL}(5):
{\wedge}^2 A(v_i)_{03} - {\wedge}^2 A(v_i)_{14} =
{\wedge}^2 A(v_i)_{04} - {\wedge}^2 A(v_i)_{23} = 0,
i = 1,...,8 \}.$

{\bf (3).}
In the coordinates
$(e_i,x_j)$,

${\wedge}^2(A) = {\wedge}^{2}(a_{kl}):
e_{ij} \mapsto
{\Sigma}_{0 \le k < l \le 4}({a_{ik}a_{jl} - a_{il}a_{jk}})e_{kl}$.

The substitution of (3) in (2)
gives a system of $16$ quadratic equations
for the elements $a_{kl}$ of $A$ --
two equations for any base vector $v_i$
of ${\bf P}^7$.

In order to solve this system, we use the following

{\bf (A.4.2)}
{\bf Lemma.}
{\sl
The group
$Aut_o(W) \subset {\bf PGL}(5)$
is a subgroup of the parabolic group
$St \ <e_0,e_1,e_2> \subset {\bf PGL}(5)$.
}

{\bf Proof.}
According to (A.2.4), the plane
${\bf P}^2_o = <e_{01},e_{02},e_{12}>$ =
${\wedge}^2 \ <e_0,e_1,e_2>$ is the
unique ${\rho}$-plane in $W$.
Since the elements of $Aut_o(W)$
are projective-linear maps of
${\bf P}^7 \supset W$,  they leave invariant
the set ${\cal P}$ of planes on $W$.
According to
lemma (A.1),
${\cal P} = C \cup P$ is a disjoint union of a smooth rational
curve $C$ -- representing the family of
${\sigma}$-planes on $W$,
and a point $P$ -- representing the unique
${\rho}$-plane on $W$.
Therefore,  any element ${\alpha}$ of $Aut_o(W)$ defines
an automorphism of
${\cal P} = C \cup P$,
i.e., $\alpha$ must leave the point $P$
and the curve $C$ invariant. Since
$\alpha = {\wedge}^2 \ A$,
for some $A \in {\bf PGL}(5)$,
the canonical isomorphism
$<e_0,e_1,e_2> \cong {\wedge}^2<e_0,e_1,e_2>$
gives that $A$ must leave the plane
$<e_0,e_1,e_2>$
invariant whenever ${\wedge}^2 \ A$ leaves the plane
$<e_{01},e_{02},e_{12}>$
invariant
{\bf q.e.d.}

{\bf (A.4.3)}
{\it Coordinate description of $Aut_o(W)$.}

Now, the system of quadratic equations (A.4.1)(2)-(3),
for the elements
$a_{ij}$ of $A \in Aut_o(W) \subset {\bf PGL}(5)$,
can be easily solved:

In the coordinates $(e_i,x_j)$,
the $(5 \times 5)$-matrix $A$ takes the form

\[
A \ = \ \left(\begin{array}{cc}
X & U  \\
O &  G
\end{array} \right)
 mod.{\bf C}^*, \mbox{ where :}
\]


\[
G =
\left( \begin{array}{cc}
a_{33} & a_{43} \\
a_{34} & a_{44}
\end{array} \right)  =  \left( \begin{array}{cc}
a & b \\
c & d
\end{array} \right)
\in \mbox{ } {\bf PGL}(2) = {\bf PSL}(2) = SL(2)/{\bf Z_2},
\]

and we can let $det(G) = ad - bc = 1$;


\[
X = {\lambda}.Symm^2(G), \lambda \in {\bf C}^*,
\mbox{ where }
Symm^2(G) =
\left( \begin{array}{ccc}
ad + bc & ac & bd \\
2ab & a^2 & b^2 \\
2cd &c^2 & d^2
\end{array} \right) ;
\]


and the $2 \times 3$-matrix

\[
U =
\left( \begin{array}{cc}
a_{30} & a_{40} \\
a_{31} & a_{41} \\
a_{32} & a_{42}
\end{array} \right)
\]

lies in the $4$-space
${\bf C}^4(G)$, defined by the condition:
${\wedge}^2 A (e_{34}) \in {\bf P}^7$,
for $v_8 = e_{34}$, i.e.:

{\bf (*).}

${\wedge}^2 A (e_{34})_{03}     -
 {\wedge}^2 A (e_{34})_{14}$    =
$(a_{30}a_{43} - a_{33}a_{40})  -
 (a_{31}a_{44} - a_{34}a_{41})$ =
$ba_{30} - aa_{40} - da_{31} + ca_{41} = 0$,

${\wedge}^2 A (e_{34})_{04}     -
 {\wedge}^2 A (e_{34})_{23}$    =
$(a_{30}a_{44} - a_{34}a_{40})  -
 (a_{32}a_{43} - a_{33}a_{42})$ =
$da_{30} - ca_{40} - ba_{32} + aa_{42} = 0$.


\[
{\bf (A.4.4)}
\mbox{ Denote the matrix }
A =
\left( \begin{array}{cc}
{\lambda}.Symm^2(G) & U \\
0 & G
\end{array} \right)
\mbox{ by }
A = ({\lambda},U,G).
\]


\bigskip

{\bf Proposition.}
{\sl
The group $Aut_o(W)$ is an extension of
${\bf PGL}(2)$
by the semidirect product
${\sf G}_m \ \triangleright \ {\sf G}_a^4$:

{\bf (*).}
\ $1 \rightarrow {\sf G}_m \ \triangleright \ {\sf G}_a^4
\rightarrow Aut_o(W) \rightarrow 0$,

where
${\sf G}_m$ is the multiplicative group of ${\bf C}$,
${\sf G}_a^4$ is the additive group of the complex $4$-space,
and:

{\bf (a).}
The induced, by the matrix multiplication in $GL(5)$,
action of the multiplicative subgroup
${\sf G}_m = \{ \left[{\lambda}\right]
:= ({\lambda},0;E_2): \lambda \in {\bf C}^* \}
\subset Aut_o(W)$,
on the
additive subgroup

\[
{\sf G}_a^4 =
\{ \left[ u \mid v \mid x \mid y \right] := (1,E_{2};U), \ U =
\left( \begin{array}{cc}
-u & -v \\
 v &  x \\
 y &  u
\end{array} \right)
: (u,v,x,y) \in {\bf C}^4 \}
\]
is: \
$\left[{\lambda}\right]: \left[ u \mid v \mid x \mid y \right] )$
$\mapsto \left[ \lambda .u \mid \lambda .v \mid \lambda .x \mid \lambda .y
\right]$;

{\bf (b).}
In the notation of (a):

(b.1). The multiplicative product, in
${\sf G}_m \subset Aut_o(W)$,
induced by the matrix multiplication in $GL(5)$, is:
$\left[ {\lambda}_1 \right]. \left[ {\lambda}_2 \right] =
\left[ {\lambda}_1 . {\lambda}_2 \right]$,

(b.2). The  additive product, in
${\sf G}_a^4 \subset Aut_o(W)$,
induced by the matrix multiplication in $GL(5)$, is:
$\left[ u_1 \mid v_1 \mid x_1 \mid y_1 \right] .
 \left[ u_2 \mid v_2 \mid x_2 \mid y_2 \right]$ =
$\left[ u_1.u_2 \mid v_1.v_2 \mid x_1.x_2 \mid y_1.y_2 \right]$.
}

{\bf Proof.}
Straightforward.

For example, equations (A.4.3)(*) imply
that
$A = (a_{ij}) = (1,E_2;U) \in  Aut_o(W)$,
i.e., $A \in {\sf G}_a^4 \subset Aut_o(W)$,
{\it iff}
$a_{30} + a_{42} = a_{40} + a_{31} = 0$,
i.e. {\it iff} $U$ has the form prescribed in (a).

{\bf (A.5)}
{\it The orbits of $Aut_o(W)$ on $W$.}

{\bf (A.5.0)}
As it follows from the proof of
lemma (A.4.2),
the action of $Aut_o(W)$ on $W$ must leave
invariant the unique ${\rho}$-plane
${\bf P}^2_o = <e_{01},e_{02},e_{12}>
\subset W$,
as well -- the ${\bf P}^1$-family
$\{ {\bf P}^2_t \}$ of ${\sigma}$-planes
on $W$.  Since the ${\sigma}$-planes
on $W$ sweep out the special hyperplane section
$Y_o$ of $W$
(see (A.1)),
$Aut_o(W)$ leaves the threefold
$Y_o \subset W$ invariant.
It can be also seen that
the special conic $q_o$ is also
an $Aut_o(W)$-invariant subvariety of $W$
(see below).
In fact, we shall prove more:
the subsets $W - Y_o$, $Y_o - {\bf P}^2_o$,
${\bf P}^2_o$, and $q_o$ are orbits of the
action of $Aut_o(W)$ on $W$.

Let $v \in W$.
Denote by $Orb \ v \subset W$ the {\it orbit}
of $v$, under the action of the group
$Aut_o(W)$, i.e.:

$Orb \ v = \{ {\alpha}(v):{\alpha} \in Aut_o(W) \}$.

{\bf (A.5.1)}
{\bf Proposition.}
{\sl
Let $W = G \cap {\bf P}^7$,
$G = G(1:{\bf P}^4) \subset
{\wedge}^2 \ {\bf P}^4 = {\bf P}^9$,
${\bf P}^7 = (x_{03} - x_{14} = x_{04} - x_{23} = 0)$
be the coordinate representation of the Fano fourfold
$W = W_5 \subset {\bf P}^7$.
Let ${\bf P}^2_o$ and $\{ {\bf P}^2_t \}$
be the unique ${\rho}$-plane and the pencil
of ${\sigma}$-planes on $W$, let
$Y_o = \cup \{ {\bf P}^2_t \} \supset {\bf P}^2_o$
be the special hyperplane section of $W$, and let
$q_o \subset {\bf P}^2_o$ be the special conic on $W$
(see (A.1)).
Then, in the coordinates $(e_i,x_j)$ from
(A.2.2):

{\bf (1).}
${\bf P}^2_o =
<e_{01}, e_{02}, e_{12}>$, and:

-- the conic
$q_o \subset {\bf P}^2(x_{01}:x_{02}:x_{12})$
has the equation:
$q_o = (x_{12}^2 - 4.x_{01}x_{02} = 0)$;

-- the hyperplane section
$Y_o \subset W$ is defined by:
$Y_o = (x_{34} = 0)$.

{\bf (2).}
The disjoint union

$W = (W - Y_o) \cup (Y_o - {\bf P}^2_o) \cup
{\bf P}^2_o) \cup q_o$

coincides with the orbit stratification of $W$,
under the action of the group $Aut_o(W)$.

In particular,
$W - Y_o = Orb \ e_{34}$,
$Y_o - {\bf P}^2_o = Orb \ e_{13}$,
${\bf P}^2_o - q_o = Orb \ e_{12}$, and
$q_o = Orb \ e_{02}.$
}

{\bf Proof.}
The proposition will be proved, if we show
that the orbits of
$e_{34}, e_{13}, e_{12}, e_{02}$,
and the corresponding subsets in (2) coincide.

To begin with, we shall find $Orb \ e_{34}$.

{\it Lemma.}
The stabilizer $St = St(e_{34}) \subset Aut_o(W)$
coincides with the extension

{\bf (**).}
\ $1 \rightarrow {\sf G}_m \rightarrow St \rightarrow
{\bf PGL}(2) \rightarrow 0$,

induced by restriction from $(A.4.4)(*)$,
on the subgroup
${\sf G}_m \cong {\sf G}_m \ \triangleright \ \{ (1,0;E_2) \} \subset
{\sf G}_m \ \triangleright \ {\sf G}_a^4$,
i.e., in the notation from (A.4.4),i.e.:

$St(e_{34}) = \{ ({\lambda},0;G):
{\lambda} \in {\bf C}^* , G \in {\bf PGL}(2) \}$
$=: {\sf G}_m \odot {\bf PGL}(2)$.

{\it Proof.}
Straightforward.

As it follows from the lemma,

$Orb \ e_{34} = {C}^4 (e_{34})$ =
$\left[ u \mid v \mid x \mid y \right] (e_{34}):
((u,v,x,y) \in {\bf C}^4$,

where the action of
$\left[ u \mid v \mid x \mid y \right]$
on $e_{34}$ is induced by the action of $Aut_o(W)$ on
$W$ -- i.e. -- by the $2$-nd wedge product
of the matrix multiplication in $GL(5)$.

Since

\[
\left[ u \mid v \mid x \mid y \right] (e_{34}) =
{\wedge}^2 (1,E_2;U)\ (e_{34}), \mbox{ where } U =
\left( \begin{array}{cc}
-u & -v \\
 v &  x \\
 y &  u
\end{array} \right), \mbox{ then:}
\]

$\left[ u \mid v \mid x \mid y \right] (e_{34})$ =

$(v^2 - ux).e_{01} + (vy - u^2).e_{02} + (uv - xy).e_{12}
+ v.(e_{03} + e_{14}) - u.(e_{04} + e_{23})
- x.e_{13} + y.e_{24} + e_{34}$,
where (u,v,x,y)$\in {\bf C}^4$.

Therefore
$Orb \ e_{34} = W - (x_{34} = 0)$ --
since the coefficient at $e_{34}$ is
always non-zero.

Clearly,
the action of the $4$-dimensional additive
group ${C}^4$, on the $4$-dimensional orbit
$Orb \ e_{34}$, is transitive.
Therefore the complement
$W - Orb \ e_{34}$ =
$W \cap (x_{34} = 0) \subset W$
can be the only $Aut_o(W)$ invariant threefold
in $W$.  However, there is already one
$Aut_o(W)$-invariant threefold in $W$ --
the special hyperplane section
$Y_o$ = (the union of ${\sigma}$-planes on $W$).
Therefore $Y_o = W \cap (x_{34} = 0)$.

The orbits of $e_{13} \in Y_o - {\bf P}^2_o$,
$e_{12} \in {\bf P}^2_o$, and $e_{02} \in q_o$
can be found in a similar way.

For example, the action of $Aut_o(W)$,
on the orbit of any $v \in {\bf P}^2_o$, is reduced
to the action of $Aut_o(W)$ on the
invariant plane
${\bf P}^2_o = <e_{01},e_{02},e_{12}>$.
In particular,

$Orb \ e_{12} = Aut_o(W) (e_{01})$ =
${\wedge}^2 \ Symm^2 \ ({\bf PGL}(2)) \ (e_{01})$ =

$\{ {\wedge}^2 \ Symm^2 \ G  \ (e_{12}) : G \in {\bf PGL}(2) \}$
=
$\{ - ab.e_{01} + cd.e_{02} + (ad + bc).e_{12} : ad - bc = 1 \}
/{\bf C}^*$,
where
$a = a_{33}, b = a_{43}, c = a_{34}, d = a_{44}$
are the coefficients of the matrix $G$;
and we have let $ad-bc = 1$,
since ${\bf PGL}(2) = {\bf PSL}(2)$.

Therefore,
$Orb \ e_{12}$ =
$(-ab:cd:ad+bc): ad-bc \neq 0$ =
${\bf P}^2_o - (x_{12}^2 - 4.x_{01}x_{02} = 0)$,

since
$(ad+bc)^2 - 4.(-ab).cd = (ad -bc)^2$, which
is always $\neq 0$ for $G \in {\bf PGL}(2)$.

In fact, the group $Symm^2 \ {\bf PGL}(2)$,
which carries out the action of
$Aut_o(W)$ on ${\bf P}^2_o$,
is the same as the projective orthogonal group
${\bf PSO}(3)$ attached to the quadratic form
$Q_o = x_{12}^2 - 4.x_{01}x_{02}$.
Therefore, the reduced action of $Aut_o(W)$ on
${\bf P}^2_o$ has two orbits:
${\bf P}^2_{o} - (Q_o = 0)$, and
$(Q_o = 0)$.

Now, it is not hard to see that the conic
$(Q_{0} = 0)$ coincides with the special conic
$q_o \subset {\bf P}^2_o$.
It is enough to see that
$q_o \subset{\bf P}^2_o$ is $Aut_o(W)$-invariant,
since our ${\bf PSO}(3)$
has $(Q_o = 0)$ as its unique invariant conic.

According to
the proof of (A.4.2),
the action of $Aut_o(W)$, on $W$,
induces an action on the rational pencil
of ${\sigma}$-planes on $W$,
as well --  an action on the conic
$q(0) \subset <e_0,e_1,e_2> \subset {\bf P}^4$,
which $q(0)$ is swept-out by their vertices
(see (A.1)).
Since
$q_o \subset {\wedge}^2<e_0,e_1,e_2>$
is, in fact, the same as the dual
conic to $q(0)$, $q_o$ is also
invariant.

In addition,
the action of $Aut_o(W)$ on $q_o$
is transitive, since the action of its restriction
${\bf PSO}(3) = {\bf PSO}(3,Q_o)$,
on $q_o = (Q_o = 0)$ is transitive.
{\bf q.e.d.}

It follows from (A.5.1) that
$W$ contains a dense $Aut_o(W)$-orbit,
e.g. -- the orbit of $e_{34}$.
That is, $W$ is a {\it prehomogeneous
space} under the action of $Aut_o(W)$
(see e.g. [M.Sato, T.Kimura, Nagoya Math. J.,
65 (1977), 1-155], or [Gua]).
In
particular,
$W$ is a closure of the quotient space
$Orb(e_{34}) = Aut_o(W)/({\sf G}_m \odot {\bf PGL}(2))$,
and any biregular automorphism of $W$ is induced
by a multiplication by a fixed element
of $Aut_o(W)$.

This implies the following:

{\bf (A.5.2)}
{\bf Corollary.}
{\sl
Let $Aut(W)$ be the group of biregular automorphisms
of the Fano $4$-fold $W = W_5 \subset {\bf P}^7$.
Then $Aut(W) = Aut_o(W)$.
}

{\bf (A.6)}
{\it The moduli space ${\cal X}_{10}$ of $X_{10}$.}

{\bf (A.6.0)}
{\it Definition of ${\cal X}_{10}$.}

As is known, (see (I.0)), any non-hyperelliptic Fano
threefold $X_{10}$ of degree $10$ is a complete
intersections of $G = G(1:{\bf P}^4) \subset {\bf P}^9$
by a supspace ${\bf p}^7$ and a quadric.
Moreover, if the intersection $G \cap {\bf P}^7$ is
smooth, then $G \cap {\bf P}^7$ is projectively
equivalent to the unique smooth Fano $4$-fold
$W = W_5$; i.e., the general $X_{10}$ is an effective
divisor of the linear system $\mid {\cal O}_W (2) \mid$.
This justifies the choice to let:

${\cal X}_{10} :=
\{ X_{10} \subset W \}/ \cong$

be the moduli space of Fano $3$-folds of
degree $10$, where $\cong$ denotes the
congruence relation induced by the biregular
equivalence of threefolds.

{\bf (A.6.1)}
{\bf Lemma.}
{\sl The general $X_{10}$ has no non-trivial biregular
automorphisms, i.e. $Aut \ X_{10} = \{ id \}$.
}

{\bf Proof.}
Let ${\Gamma}$ be the automorphism group of the general
$X_{10}$.  Then ${\Gamma}$ will be a subgroup
of $Aut \ X_{10}$ for any $X_{10}$.
In particular, let $X_{10} = X_{10}(0)$ be a general
intersection of $W$ and a reducible quadric
${\cal Q}(0) = H' + H''$, i.e. -- $H'$ and $H''$ are
two hyperplanes in ${\bf P}^7$ -- in general position.
Therefore $X_{10}(0) = Y' \cup Y''$,
where $Y' = W \cap H'$, $Y'' = W \cap H''$.
In particular,
$Aut \ (Y' \cup Y'') \supset {\Gamma}$.

Since $H'$ and $H'$ are general, $Y'$ and $Y''$
are projectively equivalent to the unique smooth
Fano threefold $Y = Y_5 \subset {\bf P}^6$
(see e.g. [Isk], [MU], [PS1], [FN]).
As is well-known that the action of $Aut \ Y_5$
stratifies $Y_5$ into a disjoint union of $3$ orbits:

$Y_5 = (Y_5 - S_o) \cup (S_o - C_o) \cup C_o$,

where $C_o$ is a special rational normal sextic
on $Y_5$, and $S_o$ is a special quadratic
section of $Y_5$;  $S_o$ coincides with the union
of the tangent lines of $C_o$.
Moreover, the restriction
$Aut \ Y_5 \rightarrow Aut \ C_o = {\bf PGL}(2)$
is an isomorphism (ibid.)

Let $S_5 = Y' \cap Y''$ be the intersection of
the components of $X_{10}(0)$.
Since $H'$ and $H''$ are general, $S_5$
is a smooth del Pezzo surface of degree $5$,
embedded anticanonically in
${\bf P}^5_o = H' \cap H''$.
Let the ``flags''
$F' := C'_o \subset S''_o \subset Y'$,
and
$F'':= C'_o \subset S''_o \subset Y''$
be the copies of
$F := C_o \subset S_o \subset Y_5$,
on $Y'$ and $Y''$,
and let ${\alpha} \in Aut \ (Y' \cup Y'')$.
Then: either ${\alpha}$ interchanges
$F'$ and $F''$,
or ${\alpha}$ acts trivially on the
pair $(F',F'')$.
Moreover, ${\alpha}$ must leave the intersection
$S_5 = Y' \cap Y''$ invariant.

In particular, let
$\{ x'_i \} = \{ x'_1,...,x'_6 \} = C' \cap S_5$,
and
$\{ x''_i \} = \{ x''_1,...,x''_6 \} = C'' \cap S_5$.
Then:
either ${\alpha}$ interchanges
the $6$-tuples
$\{ x'_i \}$ and $\{ x'_i \}$,
or ${\alpha}$ leaves each of them invariant --
depending on the action of $\alpha$
on the pair $(F',F'')$.

In either of the two possible cases,
$\alpha$ defines, by restriction,
an automorphism ${\alpha}_o$ of the curve
$C_o \cong {\bf P}^1$.
Let
$\{ x'_{io} \} \subset C_o$
and
$\{ x''_{io} \} \subset C_o$
be the copies of the $6$-tuples
$\{ x'_i \} \subset C'_o$
and
$\{ x''_i \} \subset C''_o$,
on $C_o$.
Therefore: either ${\alpha}_o$
interchanges
$\{ x'_{io} \}$ and $\{ x''_{io} \}$,
or ${\alpha}_o$ leaves either of them
invariant.
However, the general choice of
$H'$ and $H'$ implies that
$\{ x'_{i} \}$ and $\{ x''_{i} \}$,
are general hyperplane sections of
the rational normal sextics
$C'$ and $C''$.  Therefore
$\{ x'_{io} \}$ and $\{ x''_{io} \}$,
are general hyperplane sections
(i.e. -- general $6$-tuples of points)
of  $C_o \cong {\bf P}^1$.
Therefore,
${\alpha}_o \in {\bf PGL}(2)$
can be only the identity
of ${\bf PGL}(2)$.  This implies
that ${\alpha}$ does not interchange
$F'$ and $F''$, i.e.
${\alpha}_{\mid Y'} \subset Aut \ Y'$, and
${\alpha}_{\mid Y''} \subset Aut \ Y'$.
Moreover,
${\alpha}_{\mid Y'} = id_{Y'}$,
and
${\alpha}_{\mid Y''} = id_{Y''}$ --
since
${\alpha}_{\mid C'} = id$,
and
${\alpha}_{\mid C''} = id$ -- see above.
{\bf q.e.d.}

{\bf (A.6.2)}
{\bf Corollary.}
{\sl
Let $Aut(W) = Aut_o(W)$ be the $8$-dimensional
automorphism group of the Fano $4$-fold
$W = W_5 \subset {\bf P}^7$ (see (A.5.2)),  and let
${\cal X}_{10} = \mid {\cal O}_W (2) \mid / {\cong}$
be the moduli space of (non-hyperelliptic)
Fano $3$-folds of degree $10$.
Then
the moduli space ${\cal X}_{10}$ is birational
to the $22$-dimensional quotient space
$\tilde{\cal X}_{10} = \mid {\cal O}_W (2) \mid /Aut_o(W)$.
}


\newpage


\centerline{{\bf L.  Line-transformations.}}

{\bf (L.0)}
Let $X = X_{10} \in {\cal X}_{10}$ be general.
In this section we prove that the general line
$l \subset X$ defines a birational isomorphism
${\beta}_{l}:X \rightarrow X_{l}$, where
$X_{l} \in {\cal X}_{10}$. It will follow from the
forthcoming that:

{\bf (i).} $X_{l}$ is not biregular to $X$ (see section F);

{\bf (ii).} the map ${\beta}_{l}$ is ``involutive'',
      i.e., the pair $(X,l)$ defines uniquely
      a line $\underline{l} \subset X_{l}$
      such that
      $X = {\beta}_{\underline{l}}(X_{l})$;

{\bf (iii).} $X$ and $X_{l}$ have the same intermediate
      jacobian.

In particular, (i),(ii), and (iii) imply that the
map ${\beta}:{\Gamma}(X) \rightarrow {\cal X}_{10}$,
${\beta} : l \mapsto X_{l}$
defines generically an embedding of the curve of lines
${\Gamma}(X)$ into the fiber $j^{-1}(j(X))$ of the
Griffiths intermediate jacobian map
$j:{\cal X}_{10} \rightarrow {\cal J}_{10}$.

The description of the map ${\beta}_{l}$ is based
on the existence of a special $D$-flop $\rho$,
over the image of the projection of $X$ from $l$.
The discourse repeats in many ways the modern
translation of the double projection from a line
(see [Isk6],[Isk7]); however, the difference
between line transformations and double projections,
reflects
in the special choice of the divisor $D$ defining
$\rho$.  Note that, in our particular situation,
the double projection from $l$ contracts $X$ to a point.

{\bf (L.1)}
{\it The $1$-dimensional family of lines ${\Gamma}(X)$.}

{\bf (L.1.1)}
{\bf Lemma.}
{\sl
Let $X = X_{10} \subset {\bf P}^{7}$ be general,
and let $\Gamma = {\Gamma}(X)$ be the family of
lines on $X$. Then:

{\bf (i).}
$\Gamma$ is a smooth irreducible curve
(of genus $71$);

{\bf (ii).}
Let $l \in \Gamma$ be general. Then the normal bundle
$N_{l \mid X} \cong
{\cal O}_{{\bf P}^{1}}{(-1)} \oplus {\cal O}_{{\bf P}^{1}}$;

{\bf (iii).}
Let $l \in \Gamma$ be general. Then there are exactly
$11$ lines $l_{1} ,..., l_{11}$ on $X$ which intersect
$l$, $\ l_{i} \neq l$, $1 \le i \le 11$.
}

{\bf (L.1.2)}
{\bf Remark.}
The statement from (iii) is well-known
(see e.g.
[Pu], or

[D.Markushevich, Math. USSR - Sbornik,
44(2), (1981), 239-260]).

Part (ii) will be a direct corollary from general
results about lines on Fano $3$-folds of index $1$ --
provided (i) is true, esp. -- if one proves that $\Gamma$
is a smooth curve (see e.g. [Isk6, Prop.7.5]).

Part (i) will be true, if one knows the same for some
particular smooth deformations of $X_{10}$, in an
appropriate compactification of ${\cal X}_{10}$.
Such kind of smooth deformation from the general
$X_{10}$ is the Gushel $3$-fold $X_{10}'$ --
the double covering of the smooth Fano $3$-fold
$Y_{5} = G(2,5) \cap {\bf P}^{6}$, branched along an
intersection of $Y_{5}$ and a quadric
(see e.g. [Gus]).
Now, the fact that $\Gamma$ is a smooth irreducible curve
for the general Gushel $3$-fold, has been proved in
[Ili2]. This implies (i).

{\bf (L.2)}
{\it The projection from a general line $l \subset X$.}

{\bf (L.2.1)}
{\it Notation.}

Let $l \subset X = X_{10}$ be a general line
(in particular $N_{l \mid X} =
{\cal O}_{{\bf P}^{1}}{(-1)} \oplus {\cal O}_{{\bf P}^{1}}$).
Let $H_{X}$ be the hyperplane section of $X$, and let
${\pi}_{l}: X \rightarrow X''$
be the rational map defined by the
non-complete linear system
$\mid H_{X} - l \mid$, i.e.
${\pi}_{l}$ is the projection from $l$.
Let
${\sigma}_{l}: X' \rightarrow X$
be the blow-up of $l$, and let
$L' = L'_{X} = {\sigma}^{-1}(l) \subset X'$
be the exceptional divisor of
${\sigma}_{l}$;
in particular, the ruled surface $L'_{X}$
is isomorphic to
${\bf F}_{1} =
{\bf P}({\cal O}_{{\bf P}^{1}}{(-1)} \oplus
{\cal O}_{{\bf P}^{1}}$).
Denote by
$H' = H'_{X} = H_{X}^{*} - L'$
the proper preimage of the hyperplane section
$H_{X} \supset l$,
$H_{X}^{*} = {\sigma}^{*}(H_{X})$.
Let
${\pi}'_{l}:X' \rightarrow X''$
be the composition
${\pi}' = {\pi} \circ  {\sigma}_{l}$.
Clearly,
${\pi}'_{l} = {\phi}_{H'}$,
where
${\phi}_{H'}:X' \rightarrow
(\mid H_{X}' \mid )^{*} = (\mid H_{X}^{*} - L' \mid )^{*}$
is the map defined by the linear system
of $H'_{X}$ on $X'$.

{\bf (L.2.2)}
{\bf Lemma} (see [Isk6, ch.8]).

{\sl
{\bf (i).}
The image $X'' \subset {\bf P}^{5}$ is a singular
intersection of a quadric and a cubic; and the
only singularities of $X''$ are the $11$ isolated
normal singular points $x_{i}$ --
the projections of the $11$ lines $l_{i}$,
$i = 1,...,11$ which intersect $l$ (see (L.1));

{\bf (ii).}
The map ${\pi}'_{l} = {\phi}_{H'}$
is regular ( = a morphism),
i.e. ${\phi}_{H'}$ is a resolution of the projection
${\pi}_{l}$;

{\bf (iii).}
Denote the proper ${\sigma}_{l}$-preimages
of the lines $l_{1},...,l_{11}$
by $l_{1}',...,l_{11}'$. Then
the the morphism ${\phi}_{H'}$
is an isomorphism outside the
$codim.2$-cycle
$l_{1}' \cup ... \cup l_{11}'$,
and ${\phi}_{H'}$ contracts
$l_{i}'$ to the point $x_{i}$,
i = 1,...,11.
}

{\bf (L.3)}
{\it The surface $M'_{X} \subset X'$.}

Let $X = X_{10}$ and $l \subset X$,
$L' = L'_{X} = {\sigma}^{-1}(l)$, etc.
be as above.

{\bf (L.3.1)}
{\bf Proposition.}
{\sl
There exists a unique effective divisor
$M'_{X} \in \mid H'_{X} - L'_{X} \mid$ =
$\mid H^{*}_{X} - 2L'_{X} \mid$,
and $M'_{X}$ is the proper preimage of
the surface

$M_{X}$ = (the closure of)
$\cup \{ C - \mbox{ a curve on } X :
deg C = 3, p_{a}(C) = 0, \#(C \cap l) = 2 \}$;
i.e., $M_{X}$ is the union of twisted cubics
on $X$ intersecting twice the line $l$.
}

{\bf (L.3.2)}
{\bf Remarks.}

{\bf (a).}
 -- see [Isk1].
The threefold $X_{10}$ belongs to the
list of Fano threefolds of index $r = 1$,
and of the {\it 1-st series}:
$X_{2g-2} \subset {\bf P}^{g+1}$,
$6 \le g \le 12, g \neq 11$.
Here $g$ is the {\it genus} and
$d = 2g-2$ is the {\it degree} of
the threefold. In particular,
if $X = X_{10}$. i.e. $g = 6$, the existence
of $M_{X}$ -- the unique element
of the system $\mid H_{X} - 2.l \mid$
follows from the general equality
$h^{o}(H_{X_{2g-2}} - 2.l) = g -6$
(see [Isk6],[Isk7]).
In fact, the description of $M_{X}$ presented
below, proves both -- the existence and the
uniqueness of such a surface.

{\bf (b).}
Clearly, any twisted cubic $C \subset X$
which intersects twice $l$
(if such a curve on $X$ exists)
must lie on  $M_{X}$, since
$C.M_{X} = C.(H_{X} - 2l) =
deg(C) - 2.\#(C \cap l) = -1$.

Here we give an alternative proof of
the lemma, based on (b), and on the geometry
of the Grassmannian $G(2,5) \ \supset X = X_{10}$.
This will be of use for the description of the
properties of the ``line'' transformation.

{\bf Proof } of (L.3.1).
= {\it A description of $M_{X}.$}

Being a line on $X \subset G(1:{\bf P}^{4})$,
$l$ corresponds to a fixed Schubert cycle
of lines $\lambda \subset {\bf P}^{4}$:
$l = {\sigma}_{3,2}(x_{l}, {\bf P}^{2}_{l})$
=
$\{ \lambda : x_{l} \in {\lambda} \subset {\bf P}^{2}_{l} \}$.
Here
${\bf P}^{2}_{l} \subset {\bf P}^{2}$
is the ``plane'' of the Grassmann line $l$, and
$x_{l} \in {\bf P}^{2}_{l}$ is the ``center''
of $l$.

Let $Y_{o} \subset X$ be the special hyperplane
section of $X = X_{10}$ (see section A). The (general) line
$l$ intersects $Y_{o}$ in a unique point ${\lambda}_{l}$,
which does not
lie on the special conic $q_{o} \subset Y_{o}$
(see section C).
Since, outside $q_{o}$, $Y_{o}$ is a disjoint union of
the  ${\sigma}$-conics on $X$ (see section C),
there exists a unique
${\sigma}$-conic $q_{l} \subset X$ which passes through
${\lambda}_{l}$, and $q_{l}$ can be considered to be
nonsingular -- because of the generality
of the choice of $l$.
Let $y_{l}$ be the center of the
${\sigma}$-conic $q_{l}$, i.e. -- $y_{l}$ is the
vertex of the quadratic cone
$Gr(q_{l}) =
\cup \{ \lambda \in G(1:{\bf P}^{4}): \lambda \in q_{l} \}$.
Clearly,
(1). $y_{l} \neq x_{l}$ -- since $y_{l} \in q_{o}$,
and
(2). $y_{l} \in {\bf P}^{2}_{l}$ -- since
the line ${\lambda}_{l} = l \cap q_{l}$ must
pass through both centers $x_{l}$ and $y_{l}$.

Consider the line pencil
$\overline{l} =
{\sigma}_{3,2}(y_{l},{\bf P}^{2}_{l}) \subset G(1:{\bf P}^{4})$.
Clearly, $\overline{l}$ does not lie on $X$.  Indeed,
the center $y_{l}$ of $\overline{l}$ is the same as the center
of $q_{l}$, $q_{l}$ is non-singular, and the only elements
of the cycle
${\sigma}_{3,0}(y_{l}) \cap X$ (of lines in ${\bf P}^{4}$,
which pass through $Y_l$, and
which belong to $X$)  must be elements of the ${\sigma}$-conic
$q_{l}$: otherwise $X = W_{5} \cap Q$ must contain the plane
${\sigma}_{3,0}(y_{l}) \cap W_{5}$.

Let $\lambda \in \overline{l}$ be general.
Now, by the elementary Schubert
calculus on $X \cong 2.({\sigma}_{1,0})^{3}$,
we obtain that
the $1$-cycle
$\tilde{C}_{\lambda} = {\sigma}_{2,0}(\lambda) \cap X$
=
$\{ \mu \in X: \mu \cap \lambda \neq \oslash \}$
is a projectively normal connected and reduced curve,
and $\tilde{C}_{\lambda}$
has invariants $p_{a} = 2$, $deg = 6$.
However, $l$ and $q_{l}$ must be always components
of $\tilde{C}_{\lambda}$. Indeed:
(1). $\lambda \in \overline{l}$ lies
in ${\bf P}^{2}_{l}$, therefore any element of $l$
intersects $\lambda$;
(2). $y_{l} \in \lambda$, therefore
any element of $q_{l}$ intersects $\lambda$ -- in $y_{l}$.

Therefore,
$deg(C_{\lambda}) = deg(\tilde{C}_{\lambda})
- deg(q_{l}) - deg(l) = 3.$
The curve $C_{\lambda}$ must be projective normal, since $X$ --
being an intersection of quadrics -- does not contain
plane cubics; and the general $X$ does not contain planes --
since $rank \ {\bf Pic}(X) = 1$.  Now, by using similar arguments
as above, one can see that $C_{\lambda}$ intersects twice $l$.


The same result can be derived, if we use the
pencil of ${{\bf P}^3}-s$  through the plane
${\bf P}^{2}_{l}$. Indeed, let
${\bf P}^{3}({\lambda})$ be an element of the
pencil of hyperplanes in ${\bf P}^{4}$ which
pass through ${\bf P}^{2}_{l}$.
Then the elementary considerations from the
Schubert calculus on $G(2,5)$ imply that the
curve
$\tilde{C}({\lambda})$ =
${\sigma}_{1,1}({\bf P}^{3}({\lambda})).X$
=${\sigma}_{1,1}.2({\sigma}_{1,0})^{3}$ is
a quartic of arithmetical genus $1$.
In our case, $l$ is always a component of
$\tilde{C}({\lambda})$. Therefore
the curve
$C({\lambda}) := \tilde{C}({\lambda}) - l$
is a twisted cubic on $X$ which must intersect
twice the residue component $l$.

It is also not hard to see that any twisted cubic
$C \subset X$ which intersect $l$ twice must be
one of the pencil
$\{ C_{\lambda} : \lambda \in \overline{l} \}$.

Moreover, if $C = C_{\lambda}$ is such a curve,
then any element $\mu \in C$ must intersect the
plane ${\bf P}^{2}_{l}$ -- since $\mu$ intersects
$\lambda \subset {\bf P}^{2}_{l}$. Therefore
the surface $M_{X} \subset X$ lies in the
special hyperplane section
$H_{l} = {\sigma}_{1,0}({\bf P}^{2}_{l}) \cap X$.
Since $rank \ {\bf Pic(X)} = 1$, $M_{X}$ and $H_{l}$
must coincide.
Since the Schubert cycle
${\sigma}_{1,0}({\bf P}^{2}_{l}) \subset G(2,5)$
is tangent to $G(2,5)$ along the plane ${\bf P}^{2}_{l}$,
the hyperplane section $M_{X} = H_{l}$ belongs to the
linear system $\mid H_{X} - 2.l \mid$.
Clearly, $M_{X}$ is the only effective divisor in this
linear system, since any other effective divisor in it
must contain the pencil of curves
$\{ C_{\lambda} \}$.
The proper preimage $M'_{X}$ of $M_{X}$ is the only
element of the system
$\mid H'_{X} - L'_{X} \mid$
= $\mid H^{*}_{X} - 2L'_{X} \mid$.
{\bf q.e.d.}

{\bf (L.4)}
{\it The linear system $\mid 2H^{*}_{X} - 3L'_{X} \mid$.}

{\bf (L.4.1)}
It can be immediately seen that
the linear system
$\mid 2H^{*}_{X} - 3L'_{X} \mid$
= $\mid 2H'_{X} - L'_{X} \mid$
is nonempty. Indeed, let
$m'_{X}$ be the section which defines
$M'_{X}$. Then, after tensoring by
$m'_{X}$, one obtains the natural embedding
${\bf C}^{6} \cong H^{o}(H'_{X})
\subset
\mid 2H'_{X} - L'_{X} \mid$.
Clearly, any such section defines
reducible divisors, and any reducible divisor
of the last system is defined by a section
of this form.

{\bf (L.4.2)}
{\bf Proposition.}
{\sl $h^{o}(2H'_{X} - L'_{X}) \ge 8$.}

{\bf Proof.}
The linear system $\mid H'_{X} \mid$
defines the morphism
${\phi}_{H'} = {\pi}'_{l}:X' \rightarrow X''$,
where $X'' \subset {\bf P}^{5}$ is an intersection
of a quadric and a cubic (see (L.2.2)(i)).
The kernel ${\bf K}$ of the natural map
${\psi}:S^{2}H^{o}(H'_{X}) \rightarrow H^{o}(2H'_{X})$
is $1$-dimensional -- it is spanned on the only quadric
which passes through $X''$.
Therefore,
$dim.H^{o}(2H'_{X}) \ge
dim.{\psi}(S^{2}(H^{o}(H'_{X} - L'_{X}))) = 20$.

Consider the sequence

{\bf (*).}
$0 \rightarrow 2H'_{X} - L'_{X} \rightarrow
2H'_{X} \rightarrow 2H'_{X}{\mid}_{L'_{X}} \rightarrow 0$

Let $s'$ be the $(-1)$-section, and $f'$ be the fiber of
the surface $L'_{X} = {\bf F}_{1}$.  It can be seen
(see [Isk6, (8.3)])  that
$H'_{X}{\mid}_{L'_{X}} = s' + 2f'$.
Now, the inequality $h^{o}(2H'_{X}) \ge 20$,
and {\bf (*)} imply:

$h^{o}(2H'_{X} - L'_{X}) \ge
h^{o}(2H'_{X}) - h^{o}(2s' + 4f') \ge 20 - 12 = 8.$
{\bf q.e.d.}

{\bf (L.4.3)}
{\bf Proposition.}
{\sl
The linear system
$\mid 2H'_{X} - L' \mid$
has no fixed component, and the base set of
$\mid 2H'_{X} - L' \mid$
coincides with the exceptional set
$\cup \{ l'_{i} \}$
of ${\pi}'_{l} = {\phi}_{H'}$.
}

{\bf Proof.} (see [Isk6, ch.8.]), where the same
is proved for the system
$\mid H'_{X} - L'_{X} \mid =
\mid H^{*}_{X} - 2L'_{X} \mid$,
defining the double projection from
the general line
$l \subset X = X_{2g-2} \subset {\bf P}^{g+1}$,
$7 \le g \le 12, g \neq 11$.

{\bf (L.4.4)}
{\bf Remark.}
Note that -- in our case -- one has $g = 6$.  Therefore
$h^o(H^*_X - 2L'_X) = g - 5 = 1$,  and
$H^{o}(H^{*}_{X} - 2L'_{X}) = {\bf C}.m'_{X}$, where
$m'_{X}$ is the section defining $M'_{X}$;
i.e. the double projection contracts $X_{10}$
to a point. Anyway,
``the triple quadratic projection'' --
defined by the system  $\mid 2H^{*}_{X} - 3L'_{X} \mid$
-- is nontrivial.
We shall see that
the linear system $\mid 2H'_{X} - L'_{X} \mid$
defines a birational map $X' \rightarrow X_{l}$, where
$X_{l} \in {\cal X}_{10}$.

\bigskip

{\bf (L.5)}
{\it The $D'$-flop $X' \rightarrow X^{+}$.}

\centerline{----------------------}

{\bf (L.5.1).}
{\it Shortly about flops and extremal rays.}

{\bf (A).}
{\bf Definition.}

{\sc flop}:
Let ${\pi}':X' \rightarrow X''$ be a birational morphism
from the smooth $3$-fold $X'$ to the normal $3$-fold $X''$
such that:

(a'). the exceptional set of ${\phi}_{H'} = {\pi}'$  is a union
of $1$-dimensional cycles $l'_{i} \subset X'$
intersecting
trivially the anticanonical class $-K_{X'}$, i.e.
$(-K_{X'}).l'_{i} = 0, \forall i$;

(b'). the morphism ${\pi}'$ is indecomposable, i.e. ${\pi}'$
cannot be represented as a composition of morphisms
${\pi}':X' \rightarrow Z \rightarrow X''$,
$X' \neq Z \neq X''$ , $Z$ - normal.

The birational map
${\rho}:X' \rightarrow X^{+}$ is called a flop
(over $X''$) if

(1). ${\rho}$ is an isomorphism outside codimension $2$;

(2). the composition
${\pi}^{+} := {\rho}^{-1} \circ {\pi}':X^{+} \rightarrow X''$
has the properties

$(a^{+}) := (a')$, after replacing the sign $'$ by ${}^{+}$;

$(b^{+}) := (b')$, after replacing the sign $'$ by ${}^{+}$.

{\sc $D'$-flop}.
Let, moreover,
$D'$ be an effective divisor on $X'$,
such that

(c'). $D'.l'_{i} < 0, \forall i$;

and let $D^{+}$ be the proper preimage of $D'$, on $X^{+}$.

Then the flop ${\phi}'$ is called a $D'$-flop (over $X''$),
if $D^{+}$ has the property

$(c^{+})$.  $D'.l_{i}^{+} > 0 , \forall i$.

{\bf (B).}
{\bf Existence of $D'$-flops} [K].

{\sl
Let $X$ be a smooth $3$-fold, let $X''$ be a normal $3$-fold,
let ${\pi}':X \rightarrow X''$ be a birational morphism
with the properties (a') and (b'), and let
$D'$ be any effective divisor on $X'$ with the
property (c').
Then a $D'$ flop always exists, and any
sequence of $D'$-flops is finite.
}

{\bf (C).}
Denote by $\equiv$ the numerical equivalence of
cycles on the $n$-fold $X$, and -- by $[C]$
-- the numerical equivalence class of the cycle $C$.
If $C$ is effective, then -- with a probable abuse
of the notation -- we denote by  $[C]$  also
the set of effective cycles on $X$ which are
numerically equivalent to $C$.

The real space
$N(X) = \{ 1- \mbox{cycles on } X \}{/}_{\equiv}
\otimes_{\bf Z} {\bf R}$
is finite-dimensional, and the integer
${\rho}(X) = dim(N(X))$ is called the Picard number of $X$.
Let $\overline{NE}(X)$ be the closure of the convex cone
$NE(X) \ \subset N(X)$ spanned on the effective $1$-cycles
on $X$.

The half-line $R = {\bf R}_{+}.[C]$ is called an extremal
ray if $R$ is an extremal ray of the cone
$\overline{NE}(X)$,  and  $-K_{X}.C > 0$.

The rational curve $C \subset X$ is called extremal if
$0 < -K_{X}.C \le dim(X) + 1$, and
$R = {\bf R}_{+}.[C]$ is an extremal ray on $X$.

The extremal ray
$R = {\bf R}_{+}.[C]$
is called numerically effective if
$C.D \ge 0$ for any effective divisor $D$ on $X$.

{\bf (D).}
{\bf The cone theorem} [Mo].
{\sl
Let
$\overline{NE}_{K}(X) =
\{ C \in \overline{NE}(X): K_{X}.C \ge 0 \}$.

Then, for any extremal ray $R$ on $X$, there exists
an extremal curve $C$ on $X$ such that
$[C] \in N(X) - \overline{NE}_{K}(X)$
and
$R = {\bf R}_{+}.[C]$.
Moreover,
$\overline{NE}(X)$
is spanned on $\overline{NE}_{K}(X)$ and on the
set $\{ R \}$ of extremal rays on $X$; and
$\{ R \}$ is discrete outside $\overline{NE}_{K}(X)$.
}

{\bf (E).}
Let $[C] \in N(X)$.
We call the numerically effective divisor
$D$ a {\it supporting function} of $[C]$
if $D.C = 0$.

In particular the numerically effective divisor $D$
is a supporting function for any $[C]$ of the set
$D^{\perp} = \{ [C] \in N(X) : C.D = 0 \}.$

{\bf (F).}
{\bf Contraction of an extremal ray} (see [Mo])
{\sl
Any extremal ray $R$ on $X$ defines a morphism
${\phi}_{R}:X \rightarrow Y$, where $Y$ is
a normal variety, such that:

for any irreducible curve $C \subset X$,
${\phi}_{R}(C)$ is a point {\it iff}
$[C] \in R$.

Moreover, for any extremal ray $R$ there exists
a supporting function $D$ such that
$D^{\perp} \cap \overline{NE}(X) = R$.
}

\centerline{---------------------}

\newpage

{\bf (L.5.2)}
{\it The $D'$-flop ${\rho}_X : X' \rightarrow X^+$.}

Let $D'$ be any effective divisor linearly equivalent
to $2H'_{X} - L'_{X}$.

It follows from  (L.2.2)(ii)  that the morphism
${\pi}'_{l} = {\phi}_{H'}:X' \rightarrow X''$
is an isomorphism outside the $1$-dimensional
cycle
${\cup}_{i = 1,...,11} l'_{i}$.

Moreover,

$K_{X'}.l'_{i} = (- H'_{X}).l'_{i}$ =
$(- H^{*}_{X} + L'_{X}).l'_{i}$ =
$- deg.(l_{i}) + \#(l.l_{i}) = 1 - 1 = 0$;
$D'.l'_{i} = (2H'_{X} - L'_{X}).l'_{i} = -1$.

Therefore, there exists a $D'$-flop
${\rho}_{X} : X' \rightarrow X^{+}$,
over $X''$ (see above).

Let
$H^{+}_{X}, M^{+}_{X}, D^{+}_{X}$, etc.,
be the proper ${\rho}_{X}$-images of
$H'_{X}, M'_{X}, D'$, etc.

{\bf (L.5.3)}
{\bf Proposition.}
{\sl
Let $C' = C'_{\lambda} \subset X'$
(${\lambda}$-- fixed) be any ruling of
$M'_{X}$ --
the proper ${\sigma}_{l}$-preimage
of the twisted cubic $C_{\lambda}$,
$\#(C_{\lambda}.l) = 2$;
and let $C^{+} \subset X^{+}$
be the proper ${\rho}_{X}$-image
of $C'$.
Let
${\phi}' := {\phi}_{H'}$,
and let
${\phi}^{+} = {\rho}_{X}^{-1} \circ {\phi}':
X^{+} \rightarrow X''$
be as above.
Then $C^{+}$ is an extremal rational
curve on $X^{+}$.

Let $[C^+]$ be the class of effective
$1$-cycles on $X^{+}$ numerically equivalent
to $C^{+}$, and let $R := {\bf R}_{+}.[C^{+}]$ be
the extremal ray generated by $[C^{+}]$.
Then
the set of elements of $[C^+]$ coincides with
the pencil
$\{ C^{+}_{\lambda} \}$ of rulings of $M^{+}_{X}$.
In particular, the extremal morphism
${\phi}_{R}:X^{+} \rightarrow Y$ -- defined by $R$ --
contracts the rulings
${C}^{+}_{\lambda}$ of $M^{+}_{X}$
to the points of a smooth rational curve
$\overline{l} \subset Y$.
}

{\bf Proof.}

Let $C^{+}$ be any of the effective $1$-cycles
$C^{+}_{\lambda}$; let $C'$ be the proper
${\rho}_{X}$-image of $C^{+}$, and let
$C \subset X$ be the ${\sigma}_{l}$-image of $C'$.

Since ${\sigma}_{l}$ is a sigma-process along $l$,
the anticanonical class
$-K_{X'} \cong H^{*}_{X} - L'_{X} \cong H'_{X}$;
in particular $-K_{X'}$ is effective.
Then -- since
${\rho}_{X}:X' \rightarrow X^{+}$
is an isomorphism outside codimension $2$ --
the anticanonical class $-K_{X^{+}}$ is also
effective, and
$-K_{X^{+}} = H^{+}_{X}$
( = the direct ${\rho}_{X}$-image of $H'_{X}$) -- see above.

It is not hard to see that the rulings $C^{+}_{\lambda}$
are numerically equivalent on $X^{+}$. In fact,
$rk \ {\bf Pic}(X) = 2$, and we have to compute the intersections
of $C^{+}_{\lambda}$ and any two generators of ${\bf Pic}(X^{+})$.
For example,
${\bf Pic}(X^{+}) = {\bf Z}.H^{+}_{X} + {\bf Z}.L^{+}_{X}$,
and:

(1). $C^{+}_{\lambda}.L^{+}_{X} =
     C'_{\lambda}.L'_{X} =
     \#({C}_{\lambda}.l) = 2$, and

(2). $C^{+}_{\lambda}.H^{+}_{X} = C'_{\lambda}.H'_{X} =
     C'_{\lambda}.(H^{*}_{X} - L'_{X}) =
     deg(C) - \#(C.l) = 3 - 2 = 1,
     \ \forall C^{+}_{\lambda}$.

Moreover,

(3). $M^{+}_{X}.C^{+}_{\lambda} =
     (H'_{X} - L'_{X})^{+}.C^{+}_{\lambda} =
     1 - 2 = -1$;
     therefore the numerical equivalence class
     $[C]^{+}$ and the pencil
     $\{ C^{+}_{\lambda} \}$ coincide.

Let $[C^{+}]$ be the class of $C^{+}_{\lambda}$,
ant let $R = {\bf R}_{+}.[C^{+}]$.
We shall see that $R$ is an extremal ray, and
the divisor $D^{+}$ is a supporting function
for $C^{+}$.

Let $C' = C'_{\lambda}$ be irreducible.
Then $C'$ does not intersect the exceptional set
$\cup \{ l'_{i} \}$. Therefore, $C^{+}$ does not intersect
$\cup \{ l^{+}_{i} \}$.

Let $C'$ be reducible.  Then, as it follows from the
definition of the projection ${\pi}_{l}$,
$C'$ has the form
$C' = C'_{i} := q'_{i} + l'_{i} = (q_{i} + l_{i})'$
where $q_{i} \subset X$ is the unique conic
which intersects $l$ and $l_{i}$; the conic $q_{i}$
is the unique involutive of the degenerate conic
$l + l_{i}$ (see section C).

On the one hand:
$M'_{X}.C'_{i} = M'_{X}.C' = (H^{*}_{X} - 2L'_{X}).C'$ =
$deg(C) - 2.\#(l.C)$ =
$3 - 2.2 = -1$,
$D'.l_{i} = -1$; therefore
$D'.q'_{i} = 0$
(here $C'$ stands for the general $C'_{\lambda}$).

On the other hand:
$M^{+}.C^{+}_{i} = D^{+}.C^{+}$ (see (1),(2),(3)),
and
$M^{+}.C^{+} = M'_{X}.C' = -1$
(since the general $C'$
does not intersect $\cup \{ l'_{i} \}$, and
the flop ${\rho}'_{X}$ is an isomorphism outside
$\cup \{ l'_{i} \}$).
Moreover, both
$q'_{i}$ and $l'_{i}$ lie on $M'_{X}$.

Next,

$C^{+}_{i} = (q_{i} + l_{i})^{+}$ =
$q^{+}_{i} + l^{+}_{i}$, and
$M^{+}.l^{+}_{i}$ =
$(H^{+}_{X} - L^{+}_{X}).l^{+}_{i}$ =
$(D^{+} - H^{+}).l^{+}_{i}$ =
$(D^{+} + K_{X^{+}}).l^{+}_{i}$ =
$D^{+}.l^{+}_{i} + K_{X^{+}}.l^{+}_{i} \ge 0$,
by the properties $(a^{+})$ and $(c^{+})$
of the studied $D'$-flop.
Therefore $l^{+}_{i}$ does not lie on
$M^{+}_{X}$;
clearly, $q^{+}_{i} \subset M^{+}_{X}$.
Therefore,

(4). $M^{+}_{X} \subset X^{+}$ is a
(non-degenerate) ruled surface over the rational base
$\{ \lambda \}$ of the pencil $\{ C_{\lambda} \}$;
and the general $C^{+}_{\lambda}$ is the general
fiber of $M^{+}_{X}$.

We shall see that the ray
$R = {\bf R}^{+}.[C^{+}]$ is extremal.

The property (4) implies that
the ray $R = {\bf R}_{+}.[C^{+}] $
is an extremal ray the geometric cone $\overline{NE}(X)$.
Moreover,

(5). $0 < -K_{X^{+}}.C^{+} = 1 \le 4 = dim \ X^{+} + 1$.

Indeed,
$-K_{X^{+}} = H^{+}$, and
$H^{+}.C^{+} = 1$ (see above).

Properties (4) and (5) imply that
$C^{+} \subset X^{+}$ is an extremal rational curve,
i.e. $R = {\bf R}_{+}.[C^{+}]$ is an extremal ray.

We shall see that the effective divisor $D^{+}$
is a supporting function for the extremal ray $R$.
The first is to see that:

(6). The effective divisor $D^{+} \cong (2H'_{X} - L'_{X})^{+}$
     is numerically effective (nef).

{\it Remark.}
The verification of (6) is based on the same
argument used for the nef-ness for the divisor
$D_{1}^{+} \cong (H'_{X} - L'_{X})^{+}$
defining the double projection from the general line
$l \subset X_{2g-2}$, $7 \le g \le 12, g \neq 11$
(see [Isk6. ch.8]):

On the one hand, (L.4.3) implies that
$D'$ is nef outside its base set $\cup \{ l'_{i} \}$
of $X' \rightarrow X''$
Therefore, the flop-ed divisor  $D^{+}$ is
nef outside its base set $\cup \{ l^{+}_{i} \}$.

On the other hand, $D^{+}.l^{+}_{i} \ge 0$,  by the
properties of the $D'$-flop ${\rho}_{X}$.
Therefore  $D^{+}$  is nef.

The second is to see that $D^{+}$ is a supporting function
for the general $C^{+} = C^{+}_{\lambda}$.

One has:
     $\ D^{+}.C^{+}$ = $D'.C'$ = $(2H'_{X} - L'_{X}).C'$ =
     $(2H^{*}_{X} - 3L'_{X}).C'$ =
     $2.deg(C) - 3.\#(C.l)$ = $2.3 - 3.2 = 0$.
Therefore
     $R =
     {\bf R}_{+}.[C^{+}] \subset (D^{+})^{\perp}$.

It is also not hard to see that any
irreducible and reduced $1$-cycle $Z$ on
$X^{+}$ such that $D^{+}.Z^{+} = 0$
must be one of the curves $C^{+}_{\lambda}$.
Therefore

(7).  $R = {\bf R}_{+}.[C^{+}] = (D^{+})^{\perp}$,
i.e. the  nef  divisor $D^{+}$ is a supporting
function for the extremal ray $R$.

Let ${\phi}_{R}:X^{+} \rightarrow X_{l} := Y$ be the
extremal morphism defined by the extremal ray $R$.
By definition, ${\phi}_{R}$ contracts the
elements of the numerically equivalence class
$R = {\bf R}_{+}.[C^{+}] = (D^{+})^{\perp}$,
and ${\phi}_{R}$ contracts no other irreducible
curves on $X^{+}$.
Now, the property (3) implies that, in fact,
${\phi}_{R}$ contracts the ruled surface
$M^{+}_{X}$ to a rational curve
$\overline{l} \subset Y = X_{l}$.

{\bf (L.5.4)}
{\it The triple quadratic projection ${\beta}_{l}$.}

{\bf Corollary.}
{\sl
The threefold $Y$ coincides
with the image $X_{l}$ of $X'$, under the
map ${\beta}'_{l}$ defined by the linear system
$\mid 2H'_{X} - L'_{X} \mid$ =
$\mid 2H^{*}_{X} - 3L'_{X} \mid$
( = the image of $X$, under the map ${\beta}_{l}$
defined by
the system
$\mid 2H_{X} - 3l \mid$).
}

\newpage

{\bf (L.6)}
{\it The involutive property of ${\beta}_{l}$}.

{\bf Theorem.}
{\sl
Denote by
${\beta}_{l}: X \rightarrow  X_{l} = Y$
the birational map defined by
$\mid 2H_{X} - 3l \mid$.
Then
$X_{l} \in {\cal X}_{10}$,
the curve $\overline{l} \subset Y = X_{l}$
is a line on $X_{l}$, and the inverse map
$({\beta}_{l})^{-1}: Y = X_{l} \rightarrow X$
coincides with the map
${\beta}_{\overline{l}}:Y \rightarrow Y_{\overline{l}}$
defined by the system
$\mid 2H_{Y} - 3.\overline{l} \mid$
($H_{Y}$ being the hyperplane section of
$Y \in {\cal X}_{10}$).
}

{\bf Proof.}
It follows from the properties of the extremal morphism
${\phi}_{R}$ that
$rk \ {\bf Pic}(Y) = rk \ {\bf Pic}(X^{+}) - 1$ =
$rk \ {\bf Pic}(X') - 1 = rk \ {\bf Pic}(X) = 1$,
since the flop preserves $rk{\bf Pic}$ and
${\sigma}_{l} : X' \rightarrow X$ is a
${\sigma}$-process. Moreover, since
${\phi}_{R}:X^{+} \rightarrow Y$ is a blow-down
of the ruled surface $M^{+}_{X}$ to a
curve on the smooth $3$-fold $Y$ (see [Mo]),
${\phi}_{R}^{*}(-K_{Y}) = - K_{X^{+}} + M^{+}_{X}$
$\cong H^{+}_{X} + (H^{+}_{X} - L^{+}_{X})$
$\cong D^{+}_{X}$.

In particular, $-K_{Y}$ is nef.
Moreover, $-K_{Y}$ is big since its preimage
$D^{+}_{X}$ is big (since the extremal morphism
${\rho}_{R}$ is defined by a multiple of its
supporting function $D^{+}_{X}$, and
$dim(Image({\rho}_{R})) = 3$.
Therefore
$Y$ is a smooth Fano $3$-fold of the $1^{-st}$
specie (see I):
$rank \ {\bf Pic} = 1$.

Next,
${\bf Z}.H^{+}_{X} + {\bf Z}.M^{+}_{X}$ =
${\bf Pic}(X^{+})$ =
${\rho}_{R}^{*} \ {\bf Pic}(Y) \ + \ M^{+}_{X}$;
the first identity reflects the same identity
for $X'$, and the second one follows from the
properties of the ${\sigma}$-process
${\rho}_{R} = {\sigma}_{\overline{l}}$.
Moreover,
$H^{+}_{X} = -K_{X^{+}}$ and
$K_{X^{+}} = {\rho}_{R}^{*}(K_{Y}) + M^{+}_{X}$
(see above).

Therefore
${\bf Pic}(Y)$ is spanned on $-K_{Y}$ as a
${\bf Z}$-module, i.e. $Y$ is a Fano $3$-fold
of index $1$.

In particular
$deg(Y) = (-K_{Y})^{3}$.

In order to compute $(-K_{Y})^{3}$ we need to state
some identities for the intersections of divisors
on $X'$ and $X^{+}$.

(1). In ${\bf Pic}(X') = {\bf Z}. H_{X}^{*} + {\bf Z}.L'_{X}$
the following identities hold:

     $(H_{X}^{*})^{3} = 10$,
     $(H_{X}^{*})^{2}.L'_{X} = 0$,
     $H^{*}_{X}.(L'_{X})^{2} = -1$,
     $(L'_{X})^{3} = 1$.

{\it Proof of (1).}
The first identity is obvious, since
$(H_{X}^{*})^{3} = deg({\sigma}_{l}).(H_{X})^{3} =1.10 = 10$.
The next three identities are properties of the
${\sigma}$-process ${\sigma}_{l}$
(see e.g. [Yu.Manin, Uspehi Matem. Nauk.,24:5 (1969), 3-86]).
For example,
$(L'_{X})^{3} = -c_{1}(N_{l \mid X})$ =
$-deg({\cal O}(-1) \oplus {\cal O}) = 1$ (ibid.).
In particular, since
$-K_{X'} = (-K_{X})^{*} - L'_{X} = H_{X}^{*} - L'_{X}$,
and
$M'_{X} = H^{*}_{X} - 2L'_{X}$:

(2'). ${\bf Pic}(X') = {\bf Z}.(-K_{X'}) + {\bf Z}.M'_{X}$
    and

    $(-K_{X'})^{3} = (H_{X}^{*} - L'_{X})^{3} = 6$,
    $(-K_{X'})^{2}.M'_{X} = 3$,
    $-K_{X'}.(M'_{X})^{2} = -2$,
    $(M'_{X})^{3} = -10$.

{\bf (L.6.1)}
{\bf Lemma.}
{\sl
In
${\bf Pic}(X^{+}) = {\bf Z}.(-K_{X^{+}}) + {\bf Z}.M^{+}_{X}$
the following identities hold:

$(2^{+})$.
    $(-K_{X^{+}})^{3}  = 6$,
    $(-K_{X^{+}})^{2}.M^{+}_{X} = 3$,
    $-K_{X^{+}}.(M^{+}_{X})^{2} = -2$,
    $(M^{+}_{X})^{3} = 1$.
}

{\bf Proof.}
The flop ${\rho}_{X} : X' \rightarrow X^{+}$ is an
isomorphism outside the $1$-dimensional exceptional
sets
$\cup \{ l'_{i} \}$ and $\cup \{ l^{+}_{i} \}$.
This implies the first three identities.

*********************************************************

{\bf (L.6.2)}
{\bf Remark.}
{\it Description of the $D'$-flop ${\rho}_{X}$}.

In order to explain the difference between the
$4$-th identity in (2') and the $4$-th identity
in $(2^{+})$, we remember
that the set $\{ l'_{i} \}$ has $11$ elements --
the proper preimages of the $11$ lines on $X$
$l_{1},...,l_{11}$
which intersect (and not coincide with) $l$.

The $D'_{X}$-flop ${\rho}_{X} :X' \rightarrow X^{+}$
can be described as follows:

(a). Blowing-up of all the $11$ exceptional curves
     $l'_{i}, i=1,...,11$:

    ${\sigma}' :=$
    ${\Pi}_{i=1,...,11}{\sigma_{l'_{i}}}:
    \tilde{X} \rightarrow X'$.

Since $l \subset X$ is general, $l_{i}$ are general.
In particular
$N_{l_{i} \mid X} = {\cal O}(-1) \oplus {\cal O}$,
and
(since $l_{i}$ intersects the center $l$ of
the ${\sigma}$-process
${\sigma}_{l}$)  the normal sheaf
$N_{l'_{i} \mid X'} = {\cal O}(-1) \oplus {\cal O}(-1)$.

Therefore the exceptional divisor
$L_{i} = {\sigma}'(l'_{i})$ is a quadric,
and ${\sigma}'$ contracts one of the
rulings of $L_{i}$, say $f'_{i}$, onto $l'_{i}$.
Let $f^{+}_{i}$ be the $2$-nd ruling of $L_{i}$.
It is standard that $f^{+}_{i}$ defines
an extremal ray $R_{i}$ on $\tilde{X}$,
and the extremal morphism
${\rho}_{R_{i}}$ contracts
the ${\bf P}^{1}$-family $\{ f^{+}_{i} \}$.

Let
${\sigma}^{+} = {\Pi}_{i = 1,...,11}{\rho}_{R_{i}}$.
Now, the following is obvious:

(b). The composition
${\sigma}^{+} \circ ({\sigma}')^{-1}$
is a $D'_{X}$-flop (above $X''$); and we can let
${\rho}_{X} := {\sigma}^{+} \circ ({\sigma}')^{-1}$.

In particular, ${\sigma}^{+}$ contracts $L_{i}$
to the curve $l^{+}_{i}$.

It is also easy to see that

(c). $l^{+}_{i}.K_{X^{+}} = 0$,
     $l^{+}_{i}.D^{+}_{X} = +1 > 0$,
     and
     $\{ \cup l^{+}_{i} \}_{i = 1,...11}$
     is the exceptional set of the morphism
     $X^{+} \rightarrow X''$.

In particular, the $0$-cycle $(M^{+})^{3}$
is a union of the  proper image of the
$0$-cycle $(M')^{3}$ (of degree $-10$)
and the $11$ points of intersection
$z_{i} = l^{+}_{i}.M^{+}_{X}, i=1,...,11$,
i.e., $deg(M^{+}_{X})^{3} = -10 + 11 = 1$.
Next we give a direct proof of this equality.

**************************************************

For proving $4$-th identity in $(2^+)$,  we use the fact that
$M^{+}_{X}$ is a ruled surface over a rational curve.
Therefore
$(K_{M^{+}})^{2} = 8$ (always).
Now, by adjunction,

$K_{M^{+}} = (K_{X^{+}} + M^{+}){\mid}_{M^{+}}$.
Therefore

$8 = (K_{M^{+}})^{2}$ =
$(K_{X^{+}})^{2}.M^{+} + 2K_{X^{+}}.(M^{+})^{2} +
(M^{+})^{3} = 3 + 2.2 + (M^{+})^{3}$.
Therefore $(M^{+})^{3} = 1$.
{\bf q.e.d.}

{\it End of the proof of Theorem (L.6).}

Now, we can compute $deg(Y)$:

$deg(Y) = (-K_{Y})^{3} = (-K_{X^+} + M^{+})^{3}$ =
$M^{+}.(-K_{X^+} + M^{+})^{2}$ =

$(-K_{X^+})^{3} + 2M^{+}.(-K_{X^+})^{2} +
(-K_{X^+}).(M^{+})^{2}$ =
$6 + 2.3 - 2 = 10.$

This proves that the threefold
$X_{l} := Y = {\beta}_{l}(X) \in {\cal X}_{10}$.

We shall compute $deg(\overline{l})$:

$deg({\overline{l}}) = (-K_{Y}).\overline{l}$ =
$ (-K_{X^{+}} + M^{+}).(M^{+}.(-K_{X^{+}})$ =
$(-K_{X^+})^{2}.M^{+} + (-K_{X^+}).(M^+)^{2}$ =
$3 - 2 = 1$,
i.e. $\overline{l}$ is a line on $Y = X_{l}$.

The second equality in the deduction probably
needs some comments. In fact, we can replace
the $1$-cycle
$\overline{l}$ on $Y$ by the intersection
$M^{+}.(-K_{X^{+}})$ on $X^{+}$, because
the anticanonical divisor $-K_{X^{+}}$
intersects the ruled surface $M^{+}$
along a $1$-section -- the intersection
$M^{+}.(-K_{X^{+}})$; and this gives
a lifting of the intersection on $Y$
to an intersection on $X^{+}$.
In order to see that $M^{+}.(-K_{X^{+}})$
is a section of $M^{+}$ we have only to see
that $C^{+}.(-K_{X^{+}}) = 1$ for the general
fiber $C^{+} = C_{\lambda}$ of $M^{+}$.
This is easy to be computed:

$C^{+}.(-K_{X^{+}})$ =
$C^{+}.H^{+}_{X}$ =
(since $C^{+}$ is general) =
$C'.H'_{X}$ = $C'.(H^{*}_{X} - L'_{X})$ =
$deg(C) - \#(C.l) = 3 - 2 = 1$
(here $C \subset X$ is the ${\sigma}_{l}$-image
of the proper ${\rho}_{X}$-preimage $C'$ of $C^{+}$).

We have seen that the pair
$(X,l):X \in {\cal X}_{10},
l - \mbox{ a general line on } X$  defines uniquely
a pair
$(X_{l} = Y, \overline{l})$,
$\ (X_{l} \in {\cal X}_{10},
\overline{l} - \mbox{ a line on } X_{l}$,
such that:

${\bf (i,\rightarrow )}$ $X_{l} = {\beta}_{l}(X)$ =
the image of $X$ under the rational map
${\beta}_{l}$ defined by the non-complete
linear system $\mid 2H_{X} - 3l \mid$
(-- the ``triple quadratic projection'' from $l$);

${\bf (ii, \rightarrow ).}$ $\overline{l}$ is the ${\beta}_{l}$-image
of the surface $M_{X}$ swept out by the pencil of
twisted cubics $C \subset X$ intersecting twice $l$.

After obvious change of notation, we claim that

${\bf (i,\leftarrow ).}$ $X = {\beta}_{\overline{l}}(X_{l})$;

${\bf (ii,\leftarrow ).}$ $l$ =
the ${\beta}_{\overline{l}}$-image
of the surface $L_{X_{l}}$ which is swept out by
the pencil of twisted cubics $C \subset X_{l}$ which
intersect twice $\overline{l}$,
i.e.
the rational map ${\beta}_{\overline{l}}$
sends $X_{l}$ exactly to the initial threefold
$X$, and ${\beta}_{\overline{l}} = ({\beta}_{l})^{-1}$.

However, we have already proved this: the description
of the $D'$-flop from the last remark is the same
if we read it beginning from $X_{l}' := X^{+}$
and ending by $X_{l}^{+} := X'$.
Moreover, the inverse
$({\beta}_{l})^{-1}$ begins, just as ${\beta}_{l}$,
with the blow-up of a line -- the line $\overline{l} \subset X_{l}$,
and ends with the blow-down to a line -- the line $l \subset X$.
(In particular,
this proves also that $\overline{l}$
is a general line on $X_{l}$.)
{\bf q.e.d.}


\newpage


\centerline{{\bf C. Conic-transformations.}}

{\bf (C.0)}
In this section we show that the general pair
$(X,q): X \in {\cal X}_{10}, q \mbox{ a conic on } X$
defines uniquely another pair
$(X_{q}, \overline{q}): X_{q} \in {\cal X}_{10},
\overline{q}$ -- a conic on $X_{q}$.
In particular,
$X_{q} = {\alpha}_{q}(X)$
for some well-defined birational map
${\alpha}_{q}:X \rightarrow X_{q}$.
In fact, ${\alpha}_{q}$ coincides
with the ``quadruple cubic projection''
from $q \subset X$ defined by the system
$\mid 3H_{X} - 4q \mid$ ($H_{X}$ -- the
hyperplane section of $X \subset {\bf P}^7$).

Moreover, the correspondence
$(X,q) \rightarrow (X_{q},\overline{q})$
possesses a kind of involutive property --
the pair defined by the same way
by $(X_{q}, \overline{q})$ coincides with
$(X,q)$. In particular, $X$ coincides with the
birational image of $X_{q}$ under the
``quadruple cubic projection'' from
$\overline{q} \subset X_{q}$.

The study is based on the same arguments as those in
section L.  This makes it possible to leave out the
obviously identical considerations, and only to
mark out the specific differences between the line
transformations and the conic transformations.

{\bf (C.1)}
{\it The Fano surface ${\cal F}(X)$
of conics on $X \in {\cal X}_{10}$.}

{\bf (C.1.0)}
{\it Preliminaries.}

Let
$X = X_{10} \subset {\bf W}_{5} \subset {\bf P}^{7}$
be, as usual, a general smooth Fano threefold of index $1$
and of degree $10$.
The family of conics on $X$,

${\cal F}(X)$ := (the closure of)
$\{ q \subset X: q - \mbox{ an irreducible reduced 1-cycle },
deg(q) = 2 \}$

has been studied in [Pu], and -- in more detail --
in [L];
see also [Ili1] -- where ${\cal F}(X_{10})$
is studied together with the $3$-fold family of twisted cubics
on $X_{10}$, or [Ili3] -- where one can find a study of
the conics on the Gushel threefold $X'_{10}$
(see (I.0), (I.1)).

{\it Special conics on $X_{10}$.}

Let $q \subset X$ be a smooth conic on $X$.
Since the threefold
$X$ is embedded in the Grassmannian $G(2,5)$,
one can consider $q$ as a conic on $G(2,5)$.

Let
$Gr(q) := \cup \{ \lambda \}_{\lambda \in q}$
be the union of lines $\lambda \subset {\bf P}^{4}$
which are points of $q \subset G(2,5)$.
{}From this point of view, the conics on $G(2,5)$
can be separated into $3$ types:

${\tau}$-conics: $Gr(q) = Q(q)$ is a smooth quadratic surface,
and the elements $\lambda$ of $q$ are the elements of
one of the two rulings of $Q(q)$;

${\sigma}$-conics: $Gr(q) = C(q)$ is a quadratic cone,
and the elements $\lambda$ of $q$ are the elements
of the ruling of $C(q)$;

${\rho}$-conics: $Gr(q) = {\bf P}^{2}(q)$ is a plane,
and the elements $\lambda$ of $q$ are the tangent lines
of a fixed conic $q(0) \subset {\bf P}^{2}(q)$.

This notation can be applied, via degeneration,
also to reducible conics -- sum of two lines,
and to non-reduced conics -- double lines.
It is not hard to see that some of these degenerate
conics can belong to more than one of the three
types: e.g., if $q = 2l$  then $q$ is a
$\sigma$-conic, as well -- a $\rho$-conic.

{\bf (C.1.1)}
{\bf Lemma.}
{\sl
Let ${\cal F} = {\cal F}(X)$ be the family of conics
-- the Fano surface -- of the general $X \in {\cal X}_{10}$.
Then:

{\bf (i).}
${\cal F}$ is a smooth algebraic surface of general type
(and of irregularity
$10 = dim {\bf Alb}(S)$ = $dim \ J(X)$).
The base family $\{q_t \}$ of ${\sigma}$-conics on $X$
is a smooth rational rational curve on ${\cal F}(X)$;
and there exist exactly one $\rho$-conic $q_o$ on $X$.
Moreover,
the $\sigma$-conics sweep out a hyperplane
section $S_o \subset X$;
and (if $X$ is general then) the $\rho$-conic
$q_o$ is smooth.

{\bf (ii).}
The general conic $q \subset X$ is a smooth
${\tau}$-conic, with a normal bundle
$N_{q \mid X} \cong {\cal O} \oplus {\cal O}$.

{\bf (iii).}
Let $q \subset X$ be general (esp., $q$ is a $\tau$-conic).
Then there exist:

-- exactly one $\tau$-conic $\tilde{q} \subset X$
which intersects $q$ in $2$ points.

{\it Definition.}
We call $\tilde{q}$ the {\it involutive conic} of $q$;

-- exactly two $\sigma$-conics on $X$,
say $q_{1}$ and $q_{2}$
which intersect $q$; moreover $q_{1}$ and $q_{2}$
intersect the involutive conic $\tilde{q}$;

-- exactly $20$ lines $l_{1},...,l_{20}$ on $X$ which
intersect $q$.
}

{\bf (C.1.2)}
{\bf Remarks.}

{\bf (1).}
For a proof of (C.1.1) -- see [Pu], [L]; or [Ili3] --
where similar result has been proved for the Fano surface
of the hyperelliptic Fano threefold $X'_{10}$.
In particular, the Albanese variety of ${\cal F}(X_{10})$ is
$10$-dimensional, since the Abel-Jacobi map
$Alb \ {\cal F}(X) \rightarrow J(X)$ is an isomorphism
(see [L] and [Ili3]), and $dim \ J(X) = h^{1,2}(X) = 10$
(see section I).

{\bf (2).}
Let $X = W \cap Q$ be a representation of $X$
as a quadratic section of the fourfold $W = W_5$.
Concerning the $\sigma$-conics $q_t \subset X$, and the
$\rho$-conic $q_o \subset X$, these conics are obtained
as intersections of the $\sigma$-planes ${\bf P}^2_t$,
and the $\rho$-plane on $W$, with the quadric $Q$
(see (A.1), and thereafter).

{\bf (C.2)}
{\it The projection from a general conic $q \subset X$.}

{\bf (C.2.1)}
{\it Notation.}

Let $q \in {\cal F}(X)$ be a general conic on $X$.
In particular, $q$ is a smooth $\tau$-conic,
and
$N_{q \mid X} = {\cal O} \oplus {\cal O}$.
Let
${\pi}_{q}:X \rightarrow X''$
be the projection of $X$ from $q$.
Let
${\sigma}_{q}:X' \rightarrow X$ be the
blow-up of $q$, and
let $Q'_{X} = ({\pi}_{q})^{-1}(q) \subset X'$
be the exceptional divisor
$Q'_{X} \cong {\bf P}(N_{q \mid X})
\cong {\bf P}^{1} \times {\bf P}^{1}$
of ${\sigma}_{q}$.
Let $H_{X}$ be the hyperplane section of $X$,
let $H^{*}_{X} = {\sigma}^{*}(H_{X})$,
and let
$H'_{X} \cong H^{*}_{X} - Q'_{X}$ be the
proper preimage of
$H \in \mid H_{X} - q \mid$.
Since the rational projection
${\pi}_{q}$ is defined by the
system $\mid H_{X} - q \mid$,
the map
${\pi}'_{q} = {\pi}_{q} \circ {\sigma}_{q}:
X' \rightarrow X''$
is defined by the proper preimage
$\mid H'_{X} \mid = \mid H^{*}_{X} - Q'_{X} \mid$
of
$\mid H_{X} - q \mid$.

{\bf (C.2.2)}
{\bf Lemma.}
{\sl
Let
${\sigma}_{q}$,
${\pi}'_{q} = {\pi}_{q} \circ {\sigma}_{q}$,
$H_{X}$, $H'_{X}$, $l_{i}$, $q_{1},q_{2}$,
$\tilde{q}$,
etc., be as above.
Then:

{\bf (i).}
the projection $X'' \subset {\bf P}^{4}$ is
a quartic threeold with a finite number
of isolated normal singular points:
the images $x_{i}$ of the $20$ lines
$l_{i}$ which intersect $q$,
the image $z$ of the involutive conic
$\tilde{q}$, and the images $y_{1}$ and $y_{2}$
of the two ${\sigma}$-conics $q_{1}$ and $q_{2}$
which intersect $q$ (and $\tilde{q}$);

{\bf (ii).}
the map
${\pi}'_{q} = {\pi}_{q} \circ {\sigma}_{q} = {\phi}_{H'}:
X' \rightarrow X''$,
defined by
$\mid H'_{X} \mid$,
is a resolution of the projection ${\pi}_{q}$,
i.e. ${\pi}'_{q}$ is regular (= a morphism);

{\bf (iii).}
the regular map
${\pi}'_{q}: X' \rightarrow X''$
is an isomorphism, outside the $1$-dimensional
set
${\cup}_{i=1,...,20} l'_{i}
\cup q'_{1} \cup q'_{2} \cup \tilde{q}'$,
where $l'_{1},..., \tilde{q}'$ are the proper
${\sigma}_{q}$-preimages of
$l_{1},...,\tilde{q}$.
}

{\bf (C.3)}
{\it The surface $M'_{X}$.}

{\bf (C.3.1)}
{\bf Proposition.}
{\sl
There exists a unique effective divisor
$M'_{X}$, in the linear system
$\mid 2H'_{X} - Q'_{X} \mid$ =
$\mid 2H^{*}_{X} - 3Q'_{X} \mid$;  and
$M'_{X} \subset X'$ is the
proper ${\sigma}_{q}$-preimage of the surface
$M_{X}$ = (the closure of)
$\cup \{ C \subset X: C$
is an effective irreducible and reduced smooth
$1$-cycle on $X : deg(C) = 4, \#(C.q) = 3 \}$;
that is, $M_{X}$ is the union of twisted
quartics on $X$ which intersect $q$ in $3$ points.
Moreover, the set of all these quartics form
a rational pencil
$\{ C_{t}:t \in {\bf P}^{1} \}$
with only degenerations:
the $20$ curves $l_{i} + C_{i}, i = 1,...,20$,
and the $2$ curves $q_{i} + \tilde{q}, i = 1,2$;
here $C_{i} \subset X$ is the twisted cubic
involutive to the (degenerate) twisted cubic
$q + l_{i} \subset X$ -- see (C.3.3), or [Ili1].
}

{\bf Proof of (C.3.1).}

{\bf (C.3.2)}
{\it Existence $\Rightarrow$ uniqueness.}

Assume the existence of (at least one) effective
divisor
$M'_{X} \in \mid 2H'_{X} - Q'_{X} \mid$ =
$\mid H^{*}_{X} - 3Q'_{X} \mid$,
and let
$M_{X} ={\sigma}_{q}(M'_{X}) \in$
$\mid 2H_{X} - 3q \mid$
be the image of $M'_{X}$.
Let $C \subset X$ be (if exists) a twisted
quartic intersecting $q$ in $3$ points.
Then
$C.M_{X} = 2.deg(C) - 3.\#(C.q) = -1$.
Therefore $C$ must lie on $M_{X}$.
This implies the uniqueness of $M_{X}$
provided the system
$\mid 2H'_{X} - Q'_{X} \mid$ is non-empty,
and the union of twisted quartics
intersecting $q$ in $3$ points is at least
a surface on $X$.

{\it Description of $M_{X}$ = Proof of $(C.3.1)$.}

Being general in ${\cal F}(X)$, $q$ is a ${\tau}$-conic.
In particular
$Q(q) = Gr(q) =
\cup \{ \lambda \subset {\bf P}^{4}: \lambda \in q \}$
is a smooth quadratic surface in ${\bf P}^{4}$, and
the family
$\{ \lambda \subset {\bf P}^{4}: \lambda \in q \}$
coincides with one of the two rulings of $Q(q)$,
say $\Lambda$.

Let
$\overline{\Lambda}$ =
$\{\overline{\lambda} \}$
be the second ruling of $Q(q)$, and let
$\overline{\lambda} \in \overline{\Lambda}$.

Then the Schubert calculus, on
$X = 2.{\sigma}_{2,0}(\overline{\lambda}) \subset G(2,5)$,
show that the curve
$\tilde{C}(\overline{\lambda})$ =
${\sigma}_{2,0}(\overline{\lambda}).X$,
must be projectively normal, of degree $6$, and
of arithmetical genus $2$. The conic $q$ is always
a component of
$\tilde{C}(\overline{\lambda})$ .
Therefore, the residue curve
${C}(\overline{\lambda})$ =
$\tilde{C}(\overline{\lambda}) - q$
must be a twisted quartic intersecting
the conic $q$ in three points.

Now, the description of the pencil
$\{ C(\overline({\lambda}) \}$ shows
that the surface
$M_{X} := \cup \{ C(\overline{\lambda}) \} \subset X$
belongs to the complete linear system
$2H_{X}$.
Indeed, $M_{X}$
is a section of $X$ with the special quadratic
section
${\cal Q}_{q}$ := $\{ \mu \in G(2,5):$
$\mbox{ the line } \mu \mbox{ intersects the quadratic
surface } Q(q) \subset {\bf P}^{4}$
of $G(2,5) \}$.

Moreover, $M_{X}$ passes through $q$ with a multiplicity
(at least) $3$.  In fact, let
$m := mult_{q}(M_{X})$, and let $M'_{X}$ be the
proper ${\sigma}_{q}$-preimage of $M_{X}$.
Then $M'_{X} \cong 2H^{*}_{X} - m.Q'_{X}$,
and $M'_{X}$ is swept out by the pencil of
proper preimages $C'(\overline{\lambda})$
of the curves $C(\overline{\lambda})$.
Therefore
$C'(\overline{\lambda}).M'_{X}$ =
$C'(\overline{\lambda}).(2H^{*}_{X} - m.Q'_{X})$ =
$2.deg(C(\overline{\lambda}) -
m.{\#(C(\overline{\lambda}).q)}$ =
$8 - 3m$.  Moreover, the number $8 -3m$ must be
non-positive, since the curve $C'(\overline{\lambda})$
lies on $M'_{X}$.  Therefore $m \ge 3$. In particular
$M_{X} \in \mid 2H_{X} - m.q \mid
\subset \mid 2H_{X} - 3.q \mid .$

In order to see that $m =3$, we can use either the
considerations in [Isk6, ch.8] -- involving the
multiplicity of $M'_{X}(l) \in \mid H'_{X} - L'_{X} \mid$
by considering the adjoint sequences defined by
adding of $L'_{X}$,  or to use a degeneration
argument as follows:

Assume that $q = l_{1} + l_{2}$ is a sum of two lines.
Since  $m \ge 3$, the uniqueness of the already
existing effective divisor
$M_{X}(q) \in \mid 2H_{X} - m.q \mid$
follows from the same argument as in (L.1.2).
Moreover, we know that the line $l_{i}$ defines
the unique divisor
$M_{X}(l_{1}) \in \mid H_{X} - 2.l_{i} \mid$,
$i = 1,2$. Clearly, $M_{X}(l_{1})$ passes through
$l_{2}$  (with multiplicity $1$:  Otherwise
the uniqueness of
$M_{X}(l_{1})$ and $M_{X}(l_{2})$ implies that
this multiplicity will be $2$,  and
$M_{X}(l_{1}) = M_{X}(l_{2})$. In particular,
the rulings of $M_{X}(l_{1})$ and
$M_{X}(l_{2})$ must coincide, i.e.
any twisted cubic $C$ on $X$ which intersects
twice $l_{1}$ will intersect twice $l_{2}$.
The last is impossible by the construction
of the pencils of these cubics:
the residue components of the singular elliptic quartics
${\sigma}_{1,1}({\bf P}^{3}).X \supset l_i$,  for
${\bf P}^{3} \supset {\bf P}^{2}(l_i)$,
$i = 1,2$. Therefore, any ${\bf P}^{3}$
which passes through ${\bf P}^{2}(l_{1})$
must pass through ${\bf P}^{2}(l_{2})$,
i.e. the planes
${\bf P}^2(l_1)$ and ${\bf P}^2(l_2)$,
of $l_{1}$ and $l_{2}$,
must coincide.
Therefore $q = l_{1} + l_{2}$ is
a ${\rho}$-conic on $X$.  However, there exists
unique ${\rho}$-conic $q_{o} \subset X$, and
$q_{o}$ is smooth (see (C.1.1)) --
{\it contradiction.}  Therefore

$M_{X}({l_{1}}) + M_{X}({l_{2}})$
$\in \mid H_{X} - 2.l_{1} - l_{2} \mid
  +  \mid H_{X} - 2.l_{2} - l_{1} \mid$
$\subset \mid 2H_{X} - 3.(l_{1} + l_{2}) \mid$,
and the multiplicity $3$ is exact.
Now the uniqueness of $M_{X}(l_{1} + l_{2})$
implies that
$M_{X}(l_{1} + l_{2})$  =
$M_{X}({l_{1}}) + M_{X}({l_{2}})$.
Therefore $m = 3$, since the integer
multiplicity $3$ of the special conic $l_{1} + l_{2}$,
in $M_{X}(l_{1} + l_{2})$, cannot be less than
the multiplicity $m$ of the general conic $q$
in $M_{X}(q)$.

{\bf (C.3.3)}
{\it The involutions between conics ([Pu],[L]), and between
twisted cubics on $X_{10}$} ([Ili1]).

Let $q$ be a conic on $X = X_{10}$, $q \neq q_o$.
Then $<Gr(q)>$ is a $3$-space ${\bf P}^3_q \subset {\bf P}^4$.
The Schubert cycle ${\sigma}_{1,1}({\bf P}^3_q)$
intersects $X$ along a singular projectively
normal elliptic quartic $E = E^1_4$, and $q$ is a component of
$E$.  The involutive conic $\tilde{q}$ is the residue
$\tilde{q} = E - q$.
It is not hard to see that: (1). if $q$ is a $\tau$-conic,
then $\tilde{q}$ is also $\tau$-conic; (2). if $q$ is any
$\sigma$-conic, then $\tilde{q}$ is the unique $\rho$-conic
$q_o$.

Let $C$ be a general twisted cubic on $X$.  Then
-- see [Ili1] -- there exists a line
$\lambda \subset {\bf P}^4$, s.t.
$C \subset {\sigma}_{2,0}{\lambda} \cap X$
(i.e. all the elements of $C$, being lines in ${\bf P}^4$,
intersect the line $\lambda$.
The cycle ${\sigma}_{2,0}({\lambda}) \cap X$ is a singular
projectively normal curve of degree $6$, and of
arithmetical genus $2$. The involutive twisted cubic
$\tilde{C}$ is the residue
$\tilde{C} = {\sigma}_{2,0}({\lambda}) \cap X - C$.

{\bf (C.4)}
{\it The linear system
$\mid 3H'_{X} -  Q'_{X} \mid$ =
$\mid 3H'_{X} - 4Q'_{X} \mid .$}

{\bf (C.4.1)}
The following statements are the analogs of
(L.4.2), (L.4.3).

{\bf (C.4.2)}
{\bf Proposition.}
{\sl
$h^{o}(3H'_{X} - Q'_{X}) \ge 8$.
}

{\bf (C.4.3)}
{\bf Proposition.}
{\sl
The linear system
$\mid 3H'_{X} - Q'_{X} \mid$
has no fixed components, and no base
points -- outside the $1$-dimensional set
$l'_{1} \cup ... \cup l'_{20} \cup q'_{1} \cup q'_{2}
\cup \tilde{q}'$
}

{\bf (C.5)}
{\it The $D'_{X}$-flop $X' \rightarrow X^{+}$.}

{\bf (C.5.1)}
The following statements $(C.5.{}_*)$ for conic transformations,
are the analogs of the statements $(L.5.{}_*)$ for line
transformations.

{\bf (C.5.2)}
{\sl
Let $D'_{X}$ be an effective divisor, linearly equivalent
to $3H'_{X} - Q'_{X}$.  Then there exists a $D'_{X}$-flop
${\rho}_{X}: X' \rightarrow  X^{+}$, over $X''$.
}

{\bf (C.5.3)}
{\bf Proposition.}
{\sl
Let
$C' = C'(\overline{\lambda}) \subset M_{X}$
be the proper preimage
of the general rational quartic
$C(\overline{\lambda})$
which intersects $q$ in $3$ points, and let
$C^{+} \subset M^{+}_{X}$
be the proper ${\rho}_{X}$-image of $C'$.
Then the rational curve $C^{+} \subset X^{+}$
is extremal, and the (effective) elements of the class
$[C^{+}]$ are the elements of the pencil
$\{ C^{+}(\overline{\lambda}) \}$; i.e.
the extremal morphism
${\phi}_{R}:X^{+} \rightarrow Y =: X_{q}$,
defined by the extremal ray
$R = {\bf R}_{+}.[C^{+}]$,
contracts the rulings of $M^{+}_{X}$ and
no other irreducible curves on $X^{+}$.
}

The proof of (C.5.3) is the same as the proof of
(L.5.3).
In particular, the effective divisor
$D^{+}_{X} \cong (3H'_{X} - Q'_{X})^{+}$
is a supporting function for the elements of
the extremal ray $R$, and
$R = (D^{+}_{X})^{\perp}$
-- see the equality (7) in the proof of (L.5.3).

{\bf (C.5.4)}
{\it The quadruple cubic projection ${\alpha}_{q}$.}

The following is the analog of (L.5.4):

{\bf Corollary.}
{\sl
The threefold
$Y = {\rho}_{R}(X^{+})$
coincides with the image $X_{q}$ of $X'$
under the map ${\alpha}'_{q}$ defined
by the linear system
$\mid 3H'_{X} - Q'_{X} \mid$ =
$\mid 3H^{*}_{X} - 4Q'_{X} \mid$
 ( = the image of $X$ under the map
${\alpha}_{q}$ defined by the system
$\mid 3H_{X} - 4.q \mid$).
}

{\bf (C.6)}
{\it The involutive property of ${\alpha}_{q}$.}

{\bf Theorem.}
{\sl
Denote by
${\alpha}_{q}:X \rightarrow X_{q} = Y$
the birational map defined by
the system
$\mid 3H_{X} - 4.q \mid$.
Then $X_{q} \in {\cal X}_{10}$,
the curve
$\overline{q} = {\rho}_{R}(M^{+}_{X})$
is a conic on $X_{q}$,
and the inverse map
$({\alpha}_{q})^{-1}:Y = X_{q} \rightarrow X$
coincides with the map
${\alpha}_{\overline{q}}:X_{q} =Y \rightarrow Y_{\overline{q}}$
defined by the system
$\mid 3H_{Y} - 4.{\overline{q}} \mid$
( $H_{Y}$ being the hyperplane section of
$Y \in {\cal X}_{10}$).
}

{\bf Proof.}
The proof of (C.6) is similar to the
proof of (L.6) --  up to the change of notation, up to the
integers in the analogs of the
equalities $(L.6)(1)-(2')-(2^{+})$,
and up to the differences of the (unnecessary) description
of the flop attached to the conic transformation.

Here we write down these differences:

In
${\bf Pic}(X')$ =
${\bf Z}.H^{*} + {\bf Z}.Q'$ =
${\bf Z}.(-K_{X'}) + {\bf Z}.M'$,
the following analogs of (L.6)(1)-(L.6)(2') hold:

(1).    $(H^{*})^{3} = 10$,
        $(H^{*})^{2}.Q' = 0$,
        $H^{*}.(Q')^{2} = -2$,
        $(Q')^{3} = 0 (= c_{1}(N_{q \mid X})$.

(2').   $-K_{X'} = H'$, and

        $(-K_{X}')^{3} = 4$,
        $(-K_{X'})^{2}.M' = 4$,
        $(-K_{X'}).(M')^{2} = -2$,
        $(M')^{3} = -28$.

{\bf (C.6.1)}
{\bf Lemma.}
{\sl
In ${\bf Pic}(X^{+}) = {\bf Z}.(-K_{X^+}) + {\bf Z}.M^{+}$
the following analog of (L.6)$(2^{+})$ holds:

$(2^{+})$. $(-K_{X^+})^{3} = 4$,
           $(-K_{X^+})^{2}.M^{+} = 4$,
           $(-K_{X^+}).(M^{+})^{2}) = -2$,
           $(M^{+})^{3} = 0$.
}

On the base of (1),(2'),$(2^{+})$ we prove that
$deg(Y) = 10$, $deg(\overline{q}) = 2$:

$deg(Y) = 2g(Y) - 2 = (-K_{Y})^{3}$ =
$(-K_{X^+} + M^{+})^{3}$ =
$(-K_{X^+})^{3} + 3.(-K_{X^+})^{2}.M^{+} +
3.(-K_{X^+}).(M^{+})^{2} + (M^{+})^{3}$ =
$4 + 3.4 + 3.(-2) + 0 = 10$;

$deg(\overline{q})$ =
$(-K_{Y}.(\overline{q})$ =
$(-K_{X^+} + M^{+}).(M^{+}.(-K_{X^+})$ =
$(-K_{X^+})^{2}.M^{+} + (-K_{X^+}.(M^{+})^{2}$ =
$4 + (-2) = 2$;

that is
$X_{q} := Y \in {\cal X}_{10}$,
and
$\overline{q} = {\rho}_{R}(M^{+})$
is a conic on $Y = X_{q}$.


****************************************************

{\bf (C.6.2)}
{\it Description of the $D'$-flop ${\rho}_{X}$.}

Just as in (L.5.2), one can describe a $D'$ flop
(over the quartic $X''$)
${\rho}_{X}:X' \ rightarrow X^{+}$ for the
effective divisor
$D' \in \mid 2H'_{X} - 3Q'_{X} \mid$:

(a). Blow-up the $20$ exceptional curves
$l'_{1},...,l'_{20}$ -- the (disjoint) proper preimages
of the $20$ lines $l_{1},...,l_{20}$ which
intersect the conic $q$:
${\sigma}^{a} = {\Pi}_{i=1,...,20}{\sigma}_{l'_{i}}:
X^{a} \rightarrow X'$;

(b). Blow-up the proper preimage
$\tilde{q}^{a} \subset X^{a}$
of the involutive conic $\tilde{q}$:
${\sigma}^{b} = {\sigma}_{\tilde{q}}:
X^{b} \rightarrow X{a}$;

(c). Blow-up of the (disjoint) proper preimages
$q_{1}^{b} \subset X^{b}$ and
$q_{2}^{b} \subset X^{b}$ of the two
${\sigma}$-conics $q_{1}$ and $q_{2}$
which intersect $q$ (and $\tilde{q}$):
${\sigma}^{c} =
{\sigma}_{q_{1}^{b}} \circ {\sigma}_{q_{2}^{b}}:$
$X^{c} \rightarrow X^{b}$.

Since $q$ is general, the lines $l_{i}$,
the involutive conic $\tilde{q}$, and
the two ${\sigma}$-conics $q_{1}$ and $q_{2}$
can be considered as general lines, conic,
and ${\sigma}$-conics on $X$.

In particular,
$N_{l_{i} \mid X} = {\cal O}(-1) \oplus {\cal O}$,
$i = 1,...,20$,
$N_{\tilde{q} \mid X} = {\cal O} \oplus {\cal O}$,
and
$N_{{q_{j}} \mid X} = {\cal O}(-1) \oplus {\cal O}(1)$,
$j = 1,2$.

The base set $q$ of ${\sigma}_{q}$ intersects simply
each of $l_{i}$ and $q_{j}$, and $q$ intersects twice
$\overline{q}$. Therefore

$N_{l'_{i} \mid X'} = {\cal O}(-1) \oplus {\cal O}(-1)$,
$i = 1,...,20$,
$N_{\tilde{q}' \mid X'} = {\cal O}(-1) \oplus {\cal O}(-1)$,
and
$N_{q'_{j} \mid X'} = {\cal O}(-1) \oplus {\cal O}$,
$j = 1,2$.

In particular, the components
$L_{i}^{a} = ({\sigma}^{a})^{-1}(l'_{i})$ of the exceptional
divisor of ${\sigma}^{a}$ are isomorphic to smooth quadrics.

The base set $l'_{1} \cup ... \cup l'_{20}$ of
${\sigma}^{a}$ is disjoint from
$\tilde{q}'$, $q'_{1}$,
and $q'_{2}$. Therefore

$N_{\tilde{q}^{a} \mid X^{a}} = {\cal O}(-1) \oplus {\cal O}(-1)$,
and
$N_{q^{a}_{j} \mid X^{a}} = {\cal O}(-1) \oplus {\cal O}$,
$j = 1,2$.

In particular, the exceptional divisor
$\tilde{Q}^{b} = ({\sigma}^{b})^{-1}({\tilde{q}}^{a})$
is isomorphic to a smooth quadric.

The base set $\tilde{q}^{a}$ of ${\sigma}^{b}$ intersects simply
$q_{1}^{a}$ and $q_{2}^{a}$. Therefore

$N_{q^{b}_{j} \mid X^{b}} = {\cal O}(-1) \oplus {\cal O}(-1)$,
$j = 1,2$.

In particular, the components
$Q_{j}^{c} = ({\sigma}^{c})^{-1}(q_{j}^{c})$
of the exceptional divisor of ${\sigma}^{c}$
are isomorphic to smooth quadrics.

(d). Blow-down the two disjoint quadrics
$Q_{1}^{c}$ and $Q_{2}^{c}$
along their rulings which are transversal
to the rulings defining the $\sigma$-processes
${\sigma}_{q_{1}^{b}}$ and ${\sigma}_{q_{2}^{b}}$:
${\sigma}^{c+} :X^{c} \rightarrow X^{b+}$.

Let,
$q_{j}^{b+} = {\sigma}^{c+}(Q_{j}^{c}) (j = 1,2)$,
and let (as usual) $(*)^{c+}$ be the proper
${\sigma}^{c+}$-image of
the cycle $(*)^{c}$ on $X^{c}$.

It is not hard to see that the proper image
$\tilde{Q}^{c+}$ of the quadric $\tilde{Q}^{b}$
is a quadric.

(e). Blow-down the quadric $\tilde{Q}^{c+}$
along the fibering induced from the
ruling which is transversal to
the ruling defined by ${\sigma}^{b}$:
${\sigma}^{b+}: X^{b+} \rightarrow X^{a+}$.
In the obvious notation (see (d)),
the surfaces $L_{i}^{a+}$ became quadrics.

(f). Blow-down all the $L_{i}^{a+}$
(just as in (d) and (e)):

${\sigma}^{a+}: X^{a+} \rightarrow X^{+}$.

Now, it is easy to check that all the intersections
of $D^{+}$ and the curves $l_{i}^{+}$,
$q_{j}^{+}$, and $\tilde{q}^{+}$ are positive.
Moreover, $(-K_{X^{+}})$ intersects trivially
any of these curves, i.e. the composition
$\rho$ :=
${\sigma}^{a+}\circ{\sigma}^{b+}\circ{\sigma}^{c+}
\circ{\sigma}^{c}\circ{\sigma}^{b}\circ{\sigma}^{a}$
is a $D'$ flop over $X''$.

*****************************************************


\newpage


\centerline{{\bf N.  Nodal $X_{10}$.}}

{\bf (N.0)}
In this section we find the general fiber of the Griffiths map
$j$, which contain a nodal $X_{10}$.
This will be used in the next section, to describe the general
fiber of the Griffiths intermediate jacobian map
for general Fano threefolds of degree $10$.

{\bf (N.1)}
{\it The projection from the node of $X_{10}$.}

{\bf (N.1.0)}
Let
$X = X_{10} = W_{5} \cap Q \subset {\bf P}^{7}$
be a general nodal Fano $3$-fold of degree $10$,
and let $o$ be the node of $X$.

Let
${\sigma}_{o}:X' \rightarrow X$
be the blow-up of $o$, and let
$Q' = ({\sigma}_{o})^{-1}(o)$
be the exceptional quadric of
${\sigma}_{o}$.  Obviously, the composition
${\pi}'_{o} = {\pi}_{o} \circ {\sigma}_{o}:
X' \rightarrow X''$
is defined by the linear system
$\mid H' \mid \cong H^{*} - Q'$,
where
$H' \cong H^{*} - Q'$
is the proper
preimage of the hyperplane section
$H$ of $X$ through the node $o$.

Denote by $H''$ and $Q''$
the images of $H'$ and $Q'$ on $X''$.
Clearly,  $H''$ is a hyperplane section
of $X''$, and any hyperplane section of $X''$
can be represented in this way.

It is standard that the restriction
of ${\pi}'$ on $Q'$ is an isomorphism;
i.e., the surface $Q''$
is a smooth quadratic surface on $X''$.

The following lemma is an analog of (L.2.2) and (C.2.2):

{\bf (N.1.1)}
{\bf Lemma.}
{\sl

{\bf (i).}
There are exactly $6$ lines
$l_{1},...,l_{6}$ on $X \subset {\bf P}^{7}$
which pass through the node $o$.

{\bf (ii).}
The map ${\pi}':X' \rightarrow X''$ is regular,
and ${\pi}'$ is an isomorphism outside
the $1$ dimensional set
$\cup \{ l'_{i} : i = 1,...,6 \}$
where $l'_{i}$ is the proper
${\sigma}_{o}$-preimage of the line
$l_{i}, i = 1,...,6$.

{\bf (iii).}
The threefold $X''$ is a singular Fano
threefold of degree $8$ which is a complete
intersection of three quadrics in ${\bf P}^6$.

{\bf (iv).}
The only singularities of $Y$ are the $6$
isolated points $x''_{i} = {\pi}'(l'_{i})$ =
the projections of the $6$ lines
$l_{i} \subset X$
which pass through $o$, and all these points
lie on the quadric $Q'' \subset X''$.
}

{\bf Proof.}

(i).
Any Fano threefold of the $1^{-st}$ specie,
which has degree $10$ and index $1$,
can be represented as an intersection of
the  cone
${\bf C}(W) \subset {\bf P}^{8}$ (over $W$)
with a hyperplane ${\bf H}$ and a quadric ${\bf Q}$
(see [Gus]).
In this setting, the Gushel threefold $X_{10}^g$ is a
smooth projective degeneration of the
non-hyperelliptic $X_{10}$, inside ${\bf C}(W)$:

$X_{10}^{g} = {\bf C}(W) \cap {\bf H} \cap {\bf Q}$
-- s.t. the hyperplane ${\bf H}$ passes through the
vertex $v$ of the cone ${\bf C}(W)$ (ibid.).

The projection, from the vertex $v$ of the cone
${\bf C}(W)$, defines a natural double covering
$f: X_{10}^{g} \rightarrow Y$,
where
$Y = Y_{5} \subset {\bf P}^{6}$
is the Fano threefold of index $2$ and of degree $5$ --
see [Isk1].
The branch locus $B$ of $f$ is an intersection of
$Y$ and a quadric in ${\bf P}^{6}$.
Clearly, $X_{10}^{g}$ has a node
{\it iff}
the branch locus
$B \in \mid {\cal O}_{Y}(2) \mid$
has a node.
Moreover, the number of lines on the general nodal
$X_{10}^{g}$
which pass through the node of $X_{10}^g$
is the same as the number
of lines on the general nodal $X_{10}$
which pass through the node of $X_{10}$ --
since $X_{10}^g$ is a projective degeneration
of $X_{10}$, inside ${\bf C}(W)$.
Now, the lines on $X_{10}^{g}$ are the components
of the splitting $f$-preimages of lines on $Y = Y_{5}$.
In particular, the lines through the node $o$ of
the nodal $X_{10}^{g}$ are the components of the
$f$-preimages of the lines on $Y$ which pass through
the node $p(o)$ of $B$. Now, the number of lines on $Y_5$
through the general point $p(o) \in Y_{5}$ is {\it three}
(see e.g. [Isk1], [FN], [Ili2]), and
the $f$-preimage of any of these $3$ lines
is a sum of $2$ lines through the node
$o$ of $X_{10}^{g}$.
{\it q.e.d.};

(ii). standard;

(iii).
Let
${\sf I} = {\oplus}_{d \ge 0} {\sf I}_{d}$
be the graded ideal of polynomial equations
of the Grassmannian of lines in ${\bf P}^4$:
$G = G(2,5) \subset {\bf P}^9$.

Since $G$ is an intersection of quadrics,
${\sf I}$ is generated by the component
${\sf I}_{2}$.

Let ${\bf P}({\sf I_{2}})$
be the projectivization
of the space ${\sf I_{2}}$.
There exists a natural
linear isomorphism
${\bf P}^4 \rightarrow {\bf P}({\sf I_{2}})$:
$x \mapsto Pl({\bf P}^{4}/v)$, where
$Pl(v) := Pfaff({\bf P}^{4}/v)$ is
the unique Pl\"ucker quadric
-- the Pfaffian --
of the projective $3$-space
${\bf P}^{4}/v$.

The same is true for the graded ideal of
$W = G \cap {\bf P}^7$, and for the restrictions
of the Pl\"ucker quadrics on ${\bf P}^7$.
In particular, if
${\bf P}^{1}_{o} \subset {\bf P}^{4}$
be the line representing the node $o$,
then
$Pl := Pl({\bf P}^{1}_{o})$
is a projective line in
${\bf P}{\sf I_{2}} \cong {\bf P}^{4}$.

Considered as a quadric in ${\bf P}^{9}$,
the general element
$Pl(v), v \in {\bf P}^{1}_{o}$
is a quadric of rank $6$ with a vertex
-- the ${\sigma}$-space
${\bf P}^{3}(v) = {\sigma}_{3,0}(v)$
of lines in ${\bf P}^4$ through the point
$v \in {\bf P}^{1}_{o}$.
In particular,
$o \in \mbox{ the vertex of } Pl(v)$.
The same is true for the restriction
of $Pl(v)$ to ${\bf P}^{7}$.
Therefore the linear pencil
of quadrics
$Pl = \{ Pl(v) : v \in {\bf P}^{1}_{o} \}$
can be regarded as a pencil of quadrics in
${\bf P}^{6}$ -- the projection of ${\bf P}^{7}$
through the node $o$.

Let $W'' \subset {\bf P}^{6}$ be the projection
of $W = W_{5}$.  It is easy to see
that $W''$ is a complete intersection
of any two quadrics of the pencil $Pl$,
and $X'' \subset W''$.
Therefore
$Pl = {\bf P}(\sf I_{2}(W''))
\subset {\bf P}(\sf I_{2}(X''))$,
and $Pl$ is a line in the projective space
${\bf P}(\sf I_{2})(X'')$.

Moreover, as it follows from (i)-(ii),
$-K_{X'} \cong (-K_{X})^{*} - Q') \cong H'$ --
since $o$ is a double point of $X$,
and
$-K_{X''} = {\pi}'_{*}(-K_{X'}) = H''$
$\ \Leftarrow$ the regular map ${\pi}'$
is an isomorphism outside codimension $2$.
By construction, the divisor $H''$
is a hyperplane section  of $X''$.
It follows that $X''$ is a (singular) Fano threefold --
of of degree $8 = deg(X) - mult_{o}(X)$,
and $X'' \subset W''$, where the fourfold $W''$
is an intersection of two quadrics.
This implies that $X''$ is an intersection of three
quadrics in ${\bf P}^6$, i.e.
${\bf P}({\sf I_{2}}(X'')) = {\bf P}^{2}$.
This is equivalent to the existence of a quadratic
cone $Q$, such that
$sing(Q) = o$,
and
$X = W \cap Q$.  In particular, such a cone can
be also regarded as a quadric in ${\bf P}^6$,
and the plane
${\bf P}({\sf I_{2}}(X''))$ is spanned on
the line $Pl$ and on the quadric
$Q \subset {\bf P}^6$;

(iv). standard.

{\bf (N.1.2)}
{\bf Remark.}

In addition, we shall see that
{\sl
$x''_{i}$ are ordinary double points (nodes) of $X''$.
}

{\it Check-up.}

Let
$\tau: \tilde{X} \rightarrow X'$
be the blow-up of the $6$ exceptional
curves $l'_{1},...,l'_{6}$, and let
$L_{i} \subset \tilde{X}$ be the exceptional
rational ruled surfaces. The composite map
$\tilde{\pi} = {\pi}'_{o} \circ {\tau}:
\tilde{X} \rightarrow X''$ contracts
$L_{i}$ to $x''_{i}$. Since $\tilde{X}$
is non-singular, $\tilde{\pi}$ is a resolution
of $X''$, and the irreducibility of $L_{i}$
shows that this resolution is minimal.

Let $mult_{x''_{i}}(X'') = m_{i}$
for some integers
$m_{i} \ge 2, i = 1,...,6$.
Therefore
$(-K_{\tilde{X}}) = {\tilde{\pi}}^{*}(-K_{X''}) -
{\Sigma}_{i=1,...,6}(3 - m_{i}).L_{i}$.

The map
$\tau : \tilde{X} \rightarrow X'$.
blows-up the smooth irreducible curves
$l'_{i}, i =1,...,6$ on the smooth $3$-fold $X'$.
Therefore
$(-K_{\tilde{X}}) = {\tau}^{*}(-K_{X'}) -
{\Sigma}_{i=1,...,6}L_{i}$.
Moreover,
$-K_{X'} = H'$,
and
$-K_{X''} = H''$ = the proper image of $H'$.
Therefore the first summand in the two expressions
for $-K_{\tilde{X}}$ coincide, hence
${\Sigma}_{i=1,...,6}(3 - m_{i}).L_{i}$
=
${\Sigma}_{i=1,...,6}L_{i}$,
i.e.
$m_{i} = 2, i = 1,...,6$.
The double points $x''_{i}$ are ordinary nodes since
the exceptional divisors are smooth and irreducible
ruled surfaces (of degree $2 = m_{i}$) = smooth
quadrics.
{\it q.e.d.}

{\bf (N.2)}
{\it The double projection from the node $o \in X_{10}$.}

{\bf (N.2.1)}
{\it The linear system $\mid H^* - 2.Q' \mid .$}

It follows from (N.1) that any hyperplane section
of $H''$ of $X''$ is a image of a unique effective
divisor
$H' \in \mid H^* - Q' \mid$,
and the image of any such $H'$ is a hyperplane
section of $X''$.
Let, in particular,
$D'' \in \mid H'' \mid$
passes through the quadric $Q''$,
i.e. $D'' \in \mid H'' - Q'' \mid$.
Then, by the same argument,
$D''$ is a proper image of a unique effective
divisor $D' \in \mid H' - 2.Q' \mid$, and the
${\pi}'$-image
of any such $D'$ belongs to the system
$\mid H'' - Q'' \mid$.
The common image on $X$ of these systems
is the noncomplete linear system
$\mid H - 2.o \mid$
of hyperplane sections of $X$ which pass ''twice''
through the node $o$.  In particular,
$\mid H - 2.o \mid \cong
\mid H'' - Q' \mid \cong {\bf P}^2$,
i.e. it determines a rational map
$p : X \rightarrow {\bf P}^2$.
We shall see that $p$ defines a birational conic bundle
structure on $X$.

{\bf (N.2.2)}
Let
$D'' \in \mid H'' - Q'' \mid$
be general, and let
$D' \subset X'$
be the proper preimage of $D''$.
Any component $l'_{i}$ lies on $D'$,
since
$l'_{i}.D' = l'_{i}.(H^* - 2Q')$ =
$deg \ l_{i} - 2.mult_{o}(l_{i}) = -1$.
Therefore
${\pi}' : D' \rightarrow D''$
is a contraction of the six
$(-1)$-curves $l'_{i} \subset D'$
to the points $x''_{i}$, $i = 1,...,6$,
and these six points lie on the
$1$-cycle $C'' = Q'' \cap D''$.
Now, the formal adjunction on the
complete intersection divisor
$Q'' + D''$ on the Fano threefold
$X''$, and arguments involving general
position, imply:

{\bf (i).} The general $D''$ is a smooth
del Pezzo surface of degree $6$
embedded anticanonically in ${\bf P}^6$.

{\bf (ii).} The $1$-cycle $C'' = Q'' \cap D''$
is an anticanonical curve on $Q''$ and $D''$.
In particular, $C''$ is a $(2,2)$-curve
on $Q''$ which pass through
$x''_{1},...,x''_{6}$. Moreover:

{\bf (iii).} The restriction map
$D'' \rightarrow C'' = D''_{\mid Q''}$
defines an isomorphism between the
${\bf P}^2$-systems
$\mid {\cal O}_{Q''}((2,2) -
{\Sigma}_{i=1,...,6}x''_{i}) \mid$
and
$\mid H'' - Q'' \mid \cong \mid H^* - 2.Q' \mid$.

(i),(ii), and (iii) imply:

{\bf (N.2.3)}
{\bf Proposition.}
{\sl
The linear system
$\mid H^* - 2Q' \mid$
defines a rational map
$p':X' \rightarrow {\bf P}^{2}$
without fixed components, and
without fixed points outside
the $1$-dimensional exceptional set
of ${\pi}'$:
$\cup \{ l'_{i} , i = 1,...,6 \}$.
}

{\bf (N.2.4)}
{\it The $D'$-flop ${\rho}:X' \rightarrow X^+$.}

It follows from the (N.2.1)-(N.2.3) that
the general element
$D' \in \mid H^* - 2Q' \mid$
is a smooth irreducible surface
passing through
$\cup \{ l'_{i} , i = 1,...,6 \}$.
Moreover, the map ${\pi}'$ is
an isomorphism outside
$\cup \{ l'_{i} , i = 1,...,6 \}$,
and
$(-K_{X'}),l'_{i} = H'.l'_{i}
= H^*.l'_i - Q'.l'_i =
deg(l'_i) - mult_o(l_i) = 0$,
$D'.l'_i = H^*.l'_i - 2Q'.l'_i = -1 < 0$.

Therefore -- see (L.5) -- there exists a $D'$-flop
${\rho}:X' \rightarrow X^+$ over $X''$.

In particular,  let
$(*)^+$ be the ${\rho}$-image of the
cycle $(*)'$ on $X'$,
and let $\cup \{ l^+_i \}$ be the
$1$-dimensional exceptional set of
${\pi}^+:X^+ \rightarrow X''$.
Then
$-K_{X^+}.l^+_i = 0$, and
$D^+.l^+_i > 0$ --
by the properties of $\rho$.

We can describe such a $D'$-flop:

As it follows from Remark (N.1.2),
the components $L_i$ of the exceptional
set of
$\tau : \tilde{X} \rightarrow X'$
are smooth quadrics, and one can blow-down
each of these $6$ quadrics along the
ruling transversal to this defined by
$\tau$.  One obtains the map
${\tau}^+ : \tilde{X} \rightarrow X^+$.
Now, it can be seen directly that the composition
$\rho = {\tau}^+ \circ \tau : X' \rightarrow X^+$
is a $D'$-flop over $X''$.
In particular,
the components $l^+_i$ of the exceptional
set of $\rho$ are the bases of the
$6$ $\sigma$-processes which define ${\tau}^+$,
and
$D^+.l^+_i = 1 > 0, i = 1,...,6$.

{\bf (N.2.5)}
{\it The double projection from $o$ --
a rational conic bundle structure on $X$.}

Let $l_{o} \subset {\bf P}^4$ be the line representing
the node $o \in X \subset G(2,5)$.
Consider the ${\bf P}^2$-family of spaces:
${\cal P}$ :=
$\{ {\bf P}^3 \subset {\bf P}^4 : l_o \subset {\bf P}^3 \}$.
By the elementary Schubert calculus on $X = 2({\sigma}_{1,o})^3$,
the map

${\cal P} \rightarrow \{ \mbox{ 1-cycles on } X \}$:
$\ {\bf P}^3 \mapsto {\sigma}_{1,1}({\bf P}^3) \cap X$,

states an isomorphism between ${\cal P}$ and the
family ${\cal C}_4[2]$ of rational quartic curves $C$ on $X$
of arithmetical genus $1$ which have double point in the
node $o$.

Let
$C \in {\cal C}_4[2]$
and
$D' \in \mid H^* - 2Q' \mid$
be general, and let $C'$ be the proper
${\sigma}_o$-preimage of $C$.  Then
$D'.C' = (H^* - 2Q').C' =
deg(C) - 2.mult_o(C) = 4 - 2.2 = 0.$
This implies that the double projection
from $o$ contracts the curve $C$.
Now, the following is the analog of
(L.5.3) and (C.5.3):

{\bf Proposition.}
{\sl
The proper image $C^+$ of the general
$C \in {\cal C}_4[2]$
is an extremal rational curve
on $X^+$ defining the extremal ray
$R = {\bf R}_+.[C^+]$.
The general element of the numerical equivalence class
$[C^+]$ is the proper $\rho$-image of
the general $C'$, $C \in {\cal C}_4[2]$,
and the effective divisor
$D^+ \cong (H^* - 2Q')^+$
is a supporting function for
the corresponding extremal
morphism ${\phi}_R$.

Moreover, ${\phi}_R$ defines a
standard (see e.g. [Isk4])
conic bundle structure
$p^+ := {\phi}_R : X^+ \rightarrow {\bf P}^2$
on $X^+$.
In particular, the double projection
from the node $o \in X$ -- corresponding to
the system $\mid H - 2.o \mid$ --
defines a rational conic bundle structure
$p:X \rightarrow {\bf P}^2$.
}

{\bf Proof.}
\ = the same as the proofs of (L.5.3) and (C.5.3).

It rests to be seen that ${\rho}_R$ is a conic bundle
projection.
According to [Mo], this is equivalent to see
that $R$ is an extremal ray ''of type (C1)'' :=
($dim.Image({\rho}_{R}) = 2$, and
$(-K_{X^+}).C^+ = 2$).
The $2$-dimensionality of $Image({\phi}_{R})$
has been already established.  As regards
the equality,
$(-K_{X^+}).C^+ = H^+.C^+ = H'.C' =
(H^* - Q').C' = deg(C) - mult_{o}(C) = 2$.
{\bf q.e.d.}

{\bf (N.2.6)}
{\it The discriminant of $p^+$.}

Let
$\Delta = \{ x \in {\bf P}^2 : \mbox{ the conic }
q_x = (p^+)^{-1}(x) \mbox{ is singular} \}$
be the discriminant curve of $p^+$, let
$\tilde{\Delta}$
be the double discriminant curve --
the curve of components of singular fibers of
$p^+$, and let
${\pi}:\tilde{\Delta} \rightarrow {\Delta}$
be the induced double covering. Let
$(J,\Theta)$
be the principally polarized (p.p.)
intermediate jacobian of $X^+$ (= the abelian part of
the jacobian $J(X)$), and let
$(P,{\Xi})$ be the p.p. Prym variety of the pair
$(\tilde{\Delta},{\Delta})$.

The varieties $J$ and $P$ are isomorphic as
p.p. abelian varieties -- see e.g. [B1].
In particular,
$dim \ P = dim \ J = 9$ -- since $X$ has one node,
and $h^{2,1}(X_{10}) = 10$ for the general $X_{10}$.
Therefore, if one proves that $deg \ {\Delta} = 6$,
then $\Delta$ will be smooth. Indeed, if
$deg \ {\Delta} = 6$, then
$dim \ P \le p_{a}({\Delta}) -1 = 9$,
and the equality takes place
{\it iff} \ $\Delta$ is smooth, and the double
covering$\pi :\tilde{\Delta} \rightarrow \Delta$ is
unbranched.

{\bf Proposition.}
{\sl The discriminant $\Delta$ of the conic bundle
$p^+:X^+ \rightarrow {\bf P}^2$
is a smooth plane sextic, and the induced double
covering
${\pi}:\tilde{\Delta} \rightarrow \Delta$
is unbranched.
}

{\bf Proof.} It is enough to see that
$deg \ {\Delta} = 6$ -- see above.

Just as in the proof of (N.1.1)(i), it is enough to
prove this equality for the general nodal Gushel
threefold $X'$, which is known by [Ili3]
\ ${\bf q.e.d.}$

{\bf (N.3)}
{\it Nodal $X_{10}$ and plane sextics.}

{\bf (N.3.0)}
We shall see that the general plane sextic
is a discriminant curve of
the conic bundle structure $p$
of a (non-unique) nodal
$X \in {\cal X}_{10}$.
To that purpose, we prove that
the set ${\cal X}''$, of projections $X''$ of
these $X$ from their nodes,
fills in the space ${\cal Z}$
of complete intersections of three quadrics
containing a smooth quadratic surface,
and that ${\cal Z}$ admits
a natural surjective map
to the space ${\cal R}_6$ of Beauville
double coverings of plane sextics.

This will imply the two-dimensionality of the
general fiber of the Griffiths
intermediate jacobian map
$j:{\cal X}_{10} \rightarrow {\cal A}_{10}$.
In addition, the discovered by A.Verra splitting
of the fiber of the Prym map for ${\cal R}_6$ into
two Beauville pairs, reflects in the splitting of
the fiber of $j$ into two irreducible components.

{\bf (N.3.1)}
Let, as usual,
$X'' \subset {\bf P}^6$
be the projection of the nodal $X = X_{10}$ through
the node $o$.  According to (N.1.1),
$X''$ is an intersection of three
quadrics containing the smooth quadratic surface
$Q''$, and the $6$ nodes $x''_i , i = 1,...,6$
of $X''$ which lie on $Q''$.
On the contrary, let $Z$ be an intersection of
three quadrics in ${\bf P}^6$ containing a quadric
$Q$ and $6$ nodes $z_i , i = 1,..,6$ on $Q$.
Then -- as one can easily see -- after blowing-up of
$z_i , i = 1,...,6$, and then blowing-down of the
the obtained $6$ exceptional quadrics along
the appropriate rulings,
one obtains a smooth threefold $X'$,
and a birational morphism
${\pi}': X' \rightarrow Z$,
such that the restriction
${\pi}': Q' \rightarrow Q$ is an isomorphism
(here $Q'$ is the proper preimage of $Q$).
Moreover, as one can easily check, the normal
sheaf
$N_{Q' \mid X'} \cong
{\cal O}_{{\bf P}^1 \times {\bf P}^1} (-1,-1)$,
i.e. $Q'$ is contractible.  The blow-down of
$Q'$ defines a threefold $X$ with a unique
node.
We shall see that:
{\sl
$X$ is a nodal Fano threefold
of degree $10$.
}

{\it
Sketch of proof.
}
Since:

(1). \ $X' \rightarrow Z$ is an isomorphism outside
codimension $2$, and

(2). $-K_{Z} = H_Z$ is a hyperplane section of $Z$,

then
$-K_{X'} = H'$ = the proper preimage of $H_Z$.
Moreover,
${\bf Pic}(X') = {\bf Z}.H' + {\bf Z}.Q'$.
Then, after the blow-down of $Q'$:

(3). $-K_{X} = H$ = the proper image of $H'$, and

(4). ${\bf Pic}(X) = {\bf Z}.(-K_{X})$,

i.e. $X$ is a Fano $3$-fold of index $1$.

Therefore,
$deg(X) = (-K_{X})^3 = (-K_{X'})^3 + 2 =
(-K_{Z})^3 + 2 = deg(Z) + 2 = 10$
(The $2^{nd}$ equality follows from the property
of blow-up of a node. The $3^{rd}$ equality is true
since $X' \rightarrow Z$ is an isomorphism outside
codimension $2$.)

In fact, requiring -- in addition -- an existence of $6$
nodes $z_{1},...,z_6$ is superfluous.
We shall verify this.

{\bf (N.3.2)}
{\bf Proposition.}
{\sl
Let
$Z \subset {\bf P}^6$
be a complete intersection of three quadrics
containing a fixed smooth quadratic surface $Q$.
Then:

{\bf (i).}
There exists a linear pencil ${\bf P}$
of quadrics which contains $Z$,
such that
$rank({\cal Q}_t) = 6, \forall {\cal Q}_t \in {\bf P}$,
and one of the following two alternatives
takes place:

{\bf (a).}
the vertices of the general two different quadrics of ${\bf P}$
do not coincide;

{\bf (b).}
two different (hence -- the general two) quadrics
of ${\bf P}$ have a common vertex.

Moreover, if $Z$ is general, then $Z$ fulfills (a).

{\bf (ii).}
The rational projection $X''$ from the node of
the general nodal $X \in {\cal X}_{10}$ is
a complete intersection of three quadrics, containing
a quadratic surface,
for which the alternative (a) takes place.

{\bf (iii).}
The general $Z$ is canonically birational to a nodal
$X \in {\cal X}_{10}$,
and the birational isomorphism
$Z \rightarrow X$
is the same described in (N.3.1);
in particular, $Z$ has $6$ nodes lying on $Q$.
}

{\bf Proof.}

{\bf (i).}
Let ${\Pi} = {\bf P}({\sf I_{2}}(Z))$
be the plane of quadrics containing
$Z \subset {\bf P}^6$, and let
${\bf P}^3_o = Span({Q)}$.

The
${\bf P}^2$-family of quadrics
$\{ {\it Q} : {\it Q} \in \Pi \}$
intersects on ${\bf P}^3_o$
the single quadratic surface
$Q$.
Therefore, there exists a pencil
${\bf P} \subset \Pi$ of quadrics
which contain the space
${\bf P}^3_o$.  The pencil ${\bf P}$ is
linear:
 \ If
${\cal Q}_o , {\cal Q}_{\infty} \in {\bf P}$,
then
${\cal Q}_t = {\cal Q}_o + t.{\cal Q}_{\infty} \in {\bf P},
\forall t \in {\bf C} \cup {\infty}$
-- since ${\cal Q}_t \supset {\bf P}^3_o$.
That is, ${\bf P}$ is a line in
$\Pi \cong {\bf P}^2$.

The general choice of $Z$ implies that
the general ${\cal Q}_t$ is a quadric of rank $6$.
Let $v_t = sing({\cal Q}_{t})$
be the vertex of the general ${\cal Q}_{t}.$
 Clearly, $v_t \in {\bf P}^3_o$.

There are two possibilities:

$(a').$ \  $v_o \neq v_{\infty}$, and
\ $(b').$ \  $v_o = v_{\infty}$.

The possibility $(b')$ corresponds to the alternative $(b).$
In this setting, it is clear that case $(b)$ is a
degeneration of the case $(a)$ -- if two (hence -- the general
two) quadrics, containing ${\bf P}^3_o$, have  coincident
vertices.


We shall see that for the general  $Z$ all the quadrics
${\cal Q}_t \in {\bf P}$ are of rank $6$.

Let $Q$ be a (fixed) smooth quadratic surface in
${\bf P}^6$, and let
${\cal Z}$ be the moduli space

${\cal Z}$ =
$\{ Z \subset {\bf P}^6 : Z \supset Q$,
and is a complete intersection of thee quadrics
$\} / \cong .$

Let
${\cal Z}^a \subset {\cal Z}$
be the open subset of these $Z \in {\cal Z}$
for which the alternative (a) takes place.


Let $Z \in {\cal Z}$ be general,
in particular, $Z \in {\cal Z}^a$;
and let $D_k$ be the determinantal set
of quadrics of rank $\le k$, in ${\bf P}^6$,
$k = 1,...,7$.
It is well-known
(or, one can see this directly)
that

$D_k \subset D_{k+1}$, $k = 1,...,6$, and

$codim.(D_{k} \subset D_{7})$ =
$(7 - k)(7 - k + 1)/2$, $k = 1,...,6$.

Now, the general position of the line
${\bf P} \subset D_6$
gives that ${\bf P}$
does not intersect $D_5$ -- since
$codim.(D_{5} \subset D_{6}) = 2
> dim({\bf P}) = 1$,
i.e. all the quadrics
$P_t \in {\bf P} \subset D_6$
have to be of rank $6$.

It follows from the preceding that the vertex map
$v : {\bf P}^1 \cong {\bf P} \rightarrow {\bf P}^3_o$
is regular at any point $P_{t} \in {\bf P}$.

Let $C_v = v({\bf P}) \subset {\bf P}^3_o$
be the image of $v$.
Then  $C_v$ is an irreducible curve in ${\bf P}^3_o$.
The map $v$ is an isomorphism,
since otherwise the vertices of two elements of ${\bf P}$
will coincide ( = the alternative (b)).  In particular,
$C_v$ is reduced -- as a scheme-image of ${\bf P}$.
Therefore $C_v \subset {\bf P}^3_o$ is a smooth irreducible
rational curve.
This proves (i).

{\bf (iii).}
Let $d = deg(C_v)$, and let
${\cal Q} \notin {\bf P}$ be a general
quadric which contains $Z$.
It follows, from the definition of the
vertex curve $C_v$, that the fourfold
$W$ -- defined by the pencil of quadrics ${\bf P}$ --
has a double singularity along $C_v$.
Moreover, the $2d$ points
$z_1,...,z_{2d}$ of intersection
of $C_v$ and ${\cal Q}$ lie on
$Q = {\cal Q} \cap {\bf P}^3$,
and these $2d$ points are nodes of
$Z = W \cap {\cal Q}$;
and the general choice of $Z$ implies that
the only singularities of $Z$ are the
$2d$ nodes $z_{1},...,z_{2d}$.

Just as in (N.3.1), one can blow up
$z_{1},...,z_{2d}$, etc.  The preimage
$Q' \subset X'$ of $Q$ is a contractible
quadric.  After blowing down  $Q'$, one obtains
a $3$-fold $X$ which has a node $o$ on the place
of $Q$. Now, the same argument as in (N.3.1)
(see the proof of (3)-(4), in (N.3.1)) implies
that $X$ is a Fano $3$-fold of index $1$
and of degree $10$.

In particular,
$d = deg(C_v) = 3$,
i.e. the vertex curve $C_v$ is a twisted cubic,
for the general element of ${\cal Z}$.
This proves (iii).

{\bf (ii).}
Let $X'' \subset {\bf P}^6$ be
the projection from the node of the general nodal
$X \in {\cal X}_{10}$. Clearly,
$X''$ represents an element
$[X''] \in {\cal Z}$,
and we have to see that
$[X''] \in {\cal Z}^a$.

In this case, the line ${\bf P} \subset \Pi$
coincides with the line $Pf$ of Pfaffians defining
the projection $W''$ of the fourfold $W = W_5$.
We shall find this line,  in canonical coordinates
on $W$  -- see Lemma (A.2.2):

Let $x_{ij}$ be canonical coordinates on $W$.
In these coordinates
$W \subset G(2,5)$
is defined by the hyperplanes:

$H_0: \ x_{03} = x_{14}$,
$H_1: \ x_{04} = x_{23}$.

We can fix the general point $o \in W$, which is a node
of $X = X_{10}$, to be $o = e_{34}$.  Indeed, the
prehomogeneous fourfold $W$ is a closure of the
$Aut(W)$-orbit $Orb(e_{34})$ -- and all the
points of $Orb(e_{3,4}) = W - Y_o$ are
$Aut(W)$-equivalent to each other (see Proposition (A.5.1)).
Then, one can choose
$(x_{01},x_{02},x_{12},x_{03},x_{04},x_{13},x_{24})$
to be coordinates on ${\bf P}^6$.

The line $Pf$ is spanned on the quadrics:

$P_o = Pfaff(e_3)$ =
$x_{01}x_{24} - x_{02}x_{14} + x_{04}x_{12}$ =
$x_{01}x_{24} - x_{02}x_{03} + x_{04}x_{12}$
and

$P_{\infty} = Pfaff(e_4)$ =
$x_{01}x_{23} - x_{02}x_{13} + x_{03}x_{12}$ =
$x_{01}x_{04} - x_{02}x_{13} + x_{03}x_{12}$

Clearly,
$rank(P_o) = rank(P_{\infty}) = 6$,
$v(P_o) = sing(P_o) = e_{03}$,
$v(P_{\infty}) = sing(P_{\infty}) = e_{24}$,
i.e. $(a')$ takes place.

Let
$t \in {\bf C} \cup \infty$.
Then
$P_t = P_o + t.P_{\infty}$ =

$x_{01}(x_{24} + t.x_{04})
-
 x_{02}(x_{03} + t.x_{13})
+
 x_{12}(x_{04} + t.x_{03})$.

Therefore
$rank(P_t) = 6, \forall t \in {\bf C} \cup \infty$.

The common $3$-space of the quadrics
$P_t$ is
$\ {\bf P}^3_o =
Span(e_{03},e_{04},e_{13},e_{24})$.

In the projective base
$(e_{03},e_{04},e_{13},e_{24})$
of ${\bf P}^3_o$,
the vertex
$v_t = sing(P_t) =
e_{13} - t.e_{03} + t^2.e_{04} - t^3.e_{24}$.
Therefore, the vertex curve $C_{v}$
is a rational twisted cubic in
${\bf P}^3_o$, i.e. $d = 3$.
In particular, the general $X''$
is subject to the alternative (a).
This proves (ii).

{\bf (N.3.3)}
{\it Description of the space
${\cal Z} = {\cal X}''$.}

Let ${\cal X}''$ be the moduli space of
projections of the nodal $X_{10}$ from their
node, and let
${\cal Z}^a \subset {\cal Z}$ be as in (N.3.2).

As it follows from (N.3.2),
the general elements of
${\cal Z}$ and of ${\cal X}''$
coincide; in particular
$dim \ {\cal Z} =dim \ {\cal Z}^a =
dim \ {\cal X}''$ =
$dim \ {\cal X}_{10}^{nodal} =
dim \ {\cal X} - 1 = 21$
(see also (A.6.2)).

There is also a direct way to describe
the moduli space ${\cal Z}$.

{\it Description of ${\cal Z}$.}

By definition, the elements of ${\cal Z}$
are the isomorphism classes $[Z]$ of
complete intersections
$Z \subset {\bf P}^6$
of $3$  quadrics, containing a smooth quadratic
surface $Q$.  We can let $Q$ be fixed,
and let
${\bf P}^3 = span(Q) \subset {\bf P}^6$.
It can be seen directly, e.g. -- by regarding
reducible $Z \in {\cal Z}$, that the general
$Z \in {\cal Z}$ has no biregular automorphisms
(see, for example, the proof of Lemma (A.6.1));
i.e., one can identify the general $Z$ and
the class $[Z]$ of $Z$.
It follows that ${\cal Z}$ is generically isomorphic
to the quotient space

$\tilde{\cal Z}$ =
${\cal P} / St_{{\bf P}^7}(Q)$,
where
${\cal P}$ is the variety of
complete intersections of $3$ quadrics through $Q$,
and
$St_{{\bf P}^7}(Q)$ is the stabilizer of $Q$ in
$Aut({\bf P}^6) = {\bf PGL}(7)$.

On the other hand,
the affinization
${\bf C}^*.St_{{\bf P}^7}(Q)$
of
$St_{{\bf P}^7}(Q)$
is an extension of
${\bf C}^*.St_{{\bf P}^7}({\bf P}^3$,
via
${\bf C}^*.{\bf P}(SO(4)) = {\bf C}^*.St_{{\bf P}^3}(Q)$.
\ Therefore
$\ dim \ St_{{\bf P}^7}(Q)$ =
$dim \ St_{{\bf P}^7}({\bf P}^3 +
dim \ {\bf P}(SO(4)) +1$ =
$21 + 7 - 1 = 27$.

On the other hand,
the variety ${\cal P}$ of complete intersections of
three quadrics through the fixed quadratic surface
$Q$ is naturally isomorphic to the
$48$-dimensional grassmannian $G(3,19)$
of $3$-dimensional subspaces in
the $19$-dimensional space
$H^o({\cal O}_{{\bf P}^6}(2) - Q)$.

Therefore
$dim \ {\cal Z} = dim \ \tilde{\cal Z}$ =
$dim \ {\cal P} - dim \ St_{{\bf P}^7}(Q)$ =
$48 -27 = 21.$
{\bf q.e.d.}

{\bf (N.3.4)}
{\it The determinantal sextic of $Z \in {\cal Z}$.}

Let
$Z = X'' \in {\cal Z}^a = {\cal X}''$
be general, and let
${\bf P} = {\bf P}({\sf I}_{2}(W''))$
be the line of Pfaffians defining
the projection $W''$ of the fourfold
$W = W_5$ (see the proof of (N.3.2)(ii)).

Let
$\Pi = {\Pi}(Z) =
\{ {\cal Q}(s) : (s) \in {\bf P}^2 \}$
be the plane of quadrics through $Z$,
and let
${\Delta}_7(Z) =
\{ {\cal Q} \in {\Pi}(Z) : {\cal Q} \in D_6 \}$ =
$\{ (s) \in {\bf P}^2 : det({\cal Q}(s)) = 0 \}$
$\subset {\bf P}^2$
be the determinantal curve
(of degree $7$)
of the plane of quadrics $\Pi$.

For the general
$Z \in {\cal Z}$,
any
${\cal Q} \in {\Delta}_7$
is a quadric of rank $6$
(see the proof of (N.3.2)(ii));
and the two rulings of
${\cal Q}, {\cal Q} \in {\Delta}_{7}$,
define a unbranched double covering
${\pi}_7 : \tilde{\Delta}_{7} \rightarrow {\Delta}_{7}$.
(The covering is unbranched also in the double points of
${\Delta}_7$ -- the preimage of any double
point of ${\Delta}_7$ is a pair of double points of
$\tilde{\Delta}_7$.)
\ Obviously, if
$\Pi$ were a general plane in the space
$\mid {\cal O}_{{\bf P}^6}(2) \mid$,
then ${\Delta}_7$ will be smooth,
and ${\pi}_7$ will be a unbranched
double covering.
It is well-known that the intermediate jacobian
$(J(X_8),\Theta)$ of a threefold
$X_8$, defined by such $\Pi$, is isomorphic,
as a p.p. abelian variety, to the
Prym variety
$(P,\Xi)$ of the pair
$(\tilde{\Delta}_7 , {\Delta}_7)$
(see [B1]).

However,
if $Z \in {\cal Z} = {\cal X}''$
is general,
and $\Pi = {\Pi}(Z)$ is the
plane of quadrics containing $Z$,
then the determinantal curve
${\Delta}_7(Z)$ is reducible --
it contains the line of Pfaffians
${\bf P}$.
Moreover, any quadric
${\cal Q}(s) \in {\bf P}$
has rank $6$
(see (N.3.2)(ii)).
It is also not hard to see that,
for such $Z$,
the residue curve
${\Delta}(Z) = {\Delta}_7(Z) - {\bf P}$
is a smooth plane sextic
(see e.g. the proof of Proposition (N.2.6)),
which intersects ${\bf P}$
in $6$ disjoint points:
$(s)_1,...,(s)_6$.

{\it Definition.}
We call the residue curve
${\Delta} = {\Delta}(Z)$
the determinantal sextic of
$Z \in {\cal Z}^a$.

The curve
$\tilde{\Delta}_7$
splits into a union of three smooth
irreducible curves
${\Delta}_7 = \tilde{\Delta} + {\bf P}' + {\bf P}''$,
such that :

{\bf (i).}
   the double covering ${\pi}_7$ induces
   a unbranched double covering
   $\pi : \tilde{\Delta} \rightarrow \Delta$;

{\bf (ii).}
   ${\bf P}'$ and ${\bf P}''$ are disjoint
   smooth rational curves, and ${\pi}_7$
   induces on ${\bf P}'$ and ${\bf P}''$
   the isomorphisms
   ${\pi}':{\bf P}' \rightarrow {\bf P} \cong {\bf P}^1$,
   and
   ${\pi}'':{\bf P}'' \rightarrow {\bf P} \cong {\bf P}^1$.

{\bf (iii).}
   $({\pi}_{7})^{-1}(\{ (s)_1 ,...,(s)_6 \})$ =
   $\{ (s)'_1, .., (s)'_6 \}
   \cup
   \{ (s)''_1, .., (s)''_6 \}$,
   where
   $\{ (s)'_1, .., (s)'_6 \} =
   \tilde{\Delta} \cap {\bf P}'$,
   $\{ (s)''_1, .., (s)''_6 \} =
   \tilde{\Delta} \cap {\bf P}'$,
   and
   ${\pi}_7(s)'_i) = {\pi}_7(s)''_i = (s)_i,
   i = 1,...,6$.

{\bf (N.3.5)}
{\bf Proposition.}
{\sl
Let $Z \in {\cal Z}$ be general,
let $\Delta$ be the determinantal sextic
of $Z$, and let
$\pi : \tilde{\Delta} \rightarrow \Delta$
be the induced by ${\pi}_7$ unbranched
double covering, and let
$(J(Z),\Theta)$ be (the abelian part of)
the intermediate jacobian of $Z$.
Then $(J(Z), \Theta )$ is isomorphic,
as a p.p. abelian variety, to the
Prym variety
$(P,\Xi)$ of the pair
$(\tilde{\Delta},{\Delta})$.
}

The proof is a straightforward
specialization of the same statement
for the general complete intersection of
three quadrics
(see e.g. [Tju, ch.3, sect.2]).

{\bf (N.3.6)}
{\it The map
$\tilde{det}:{\cal Z} \rightarrow {\cal R}_{6+1}^{even}$.}

It follows from (N.3.4) that the general
$Z \in {\cal Z}$
defines uniquely the determinantal pair
$(\tilde{\Delta}_7,{\Delta}_7)$ =
$(\tilde{\Delta} \cup {\bf P}' \cup {\bf P}'',
{\Delta} \cup {\bf P})$
satisfying the properties
(N.3.4)(i),(ii),(iii).
Moreover, since
$Z$ is a complete intersection of three quadrics,
the determinantal pair of $Z$ must be {\it even}
(see e.g. [Tju, ch.3,sect.2]).

We shall see that the opposite is also true.

Let
${\cal R}_d (d = 6,7)$
be the moduli space of
double coverings of plane curves of degree $d$,
and let
${\cal R}_{6+1} \subset {\cal R}_7$
be the moduli space of double coverings
satisfying (N.3.4)(i),(ii),(iii).
Let
${\cal R}_{7}^{even} \subset {\cal R}_{7}$,
(resp.
${\cal R}_{6+1}^{even} \subset {\cal R}_{6+1}$)
be their components of even
(i.e. admissible -- see [FS]) coverings
of corresponding plane septics.

Let
${\cal X}_{2.2.2}$
be the moduli space of complete intersections
of three quadrics in ${\bf P}^6$.
The determinantal map
$\tilde{det}:{\cal X}_{2.2.2} \rightarrow {\cal R}_7^{even}$,
$X_{2.2.2} \rightarrow (\tilde{\Delta},{\Delta})(X_{2.2.2})$
is well-defined in the generic point $X_{2.2.2}$
of ${\cal X}_{2.2.2}$, as well -- in the general
point $Z$ of the subspace
${\cal Z} \subset {\cal X}_{2.2.2}$.
Moreover,
$\tilde{det}({\cal Z}) \subset {\cal R}_{6+1}^{even}$.

{\bf Proposition.}
{\sl
(The restriction of) the determinantal map
$\tilde{det}: {\cal Z} \rightarrow {\cal R}_{6+1}^{even}$
is surjective and of degree $1$ in
the general point.
}

{\bf Lemma.}
{\sl
Let ${\bf P}$ be a general linear pencil of quadrics of
rank $6$ in ${\bf P}^6$.
Then the quadrics of ${\bf P}$ have a common ${\bf P}^3$.
}

{\bf Proof} of the lemma.
Assume that
${\bf P} = \{ {\cal Q}(t) : (t) \in {\bf P}^1 \}$
is a general linear pencil of quadrics of rank $6$
in ${\bf P}^6$, and let
$W_{2.2}$ be the base set of ${\bf P}$.
Any ${\cal Q}(t)$ is a cone with a $0$-dimensional
vertex -- a point $v(t) \in {\bf P}^6$.
Let $C_v = \{ v(t) : t \in {\bf P}^1 \}$
be the set of  vertices of ${\cal Q}(t)$.
By the theorem of Bertini,
$C_v \subset W_{2.2}$.
The general choice of ${\bf P}$ implies that
$C_v$ is a curve. Moreover,
$C_v \subset Sing(W_{2.2})$.

Let $Span(C_v)$ be the projective linear span
of $C_v \subset {\bf P}^6$.  Since the base set
$W_{2.2}$ is an intersection of quadrics, and
$C_v \subset Sing(W_{2.2})$, then any bisecant
line of $C_v \subset Sing(W_{2.2})$ lies in $W_{2.2}$.
Therefore
$Span(C_v) \subset W_{2.2}$, i.e. the projective
space
${\bf P}^d = Span(C_v)$ lies on any quadric
${\cal Q}(t) \in {\bf P}$. In particular,
$d = dim \ {\bf P}^d =: d({\bf P}) \le 3$, since ${\bf P}^3$
is the projective subspace of maximal dimension,
in a quadric of rank $6$ in ${\bf P}^6$.
Clearly, the integer-valued function
$d : {\bf P} \rightarrow d({\bf P})$ takes the
maximum value $d_o$ in the general pencil ${\bf P}$.
If we find a pencil ${\bf P}_{o}$ such that
$d({\bf P}_o) = 3$, then $d_o = 3$.
This is the case, when ${\bf P}$ comes from the
set of Pfaffians defining the projection
$W_{2.2}$ of the fourfold $W_5$, through a general
point $o \in W_5$ (see the proof of (N.3.2)(ii)).
Therefore $d_o = 3$.
{\it q.e.d.}

{\bf Proof}
of the proposition.
On the one hand, it is known by [FS] that
$\tilde{det}$
defines
a generically one-to-one correspondence
$\tilde{det} \subset {\cal X}_{2.2.2} \times {\cal R}_7^{even}$
(= a birational isomorphism
$\tilde{det} : {\cal X}_{2.2.2} \rightarrow {\cal R}_7^{even}$);
see also [D1].
In particular, any fiber of $\tilde{det}$
has to be connected -- since the general fiber
(= a point) is connected.


On the other hand,
the lemma implies that
if
$(\tilde{\Delta},{\Delta}) \in {\cal R}_{6+1}^{even}$,
then the elements of the line ${\bf P}$, occurring as
a component of ${\Delta}_7$, have a common
${\bf P}^3$.
Therefore, any $X \in {\bf X}_{2.2.2}$
defined by a plane of quadrics
$\Pi \supset {\bf P}$, must contain
a quadratic surface, i.e.
$X \in {\cal Z}$.

In particular,

(1). \ the fiber $F$ of $\tilde{det}$,
    over the general element  ${\cal R}_{6+1}^{even}$
    must lie entirely in ${\cal Z}$,
    i.e. the set of these $Z \in {\cal Z}$
    which have the same determinantal pair, must
    be connected (see above).

The surjectivity of
$\tilde{det}:{\cal X}_{2.2.2} \rightarrow {\cal R}_7^{even}$,
together with (1), imply:

(2). $\tilde{det}: {\cal Z} \rightarrow {\cal R}_{6+1}^{even}$
     is surjective;

(3). the $\tilde{det}$-preimage of ${\cal R}_{6+1}^{even}$
     lies entirely in ${\cal Z}$.

It rests to see that
$dim \ {\cal R}_{6+1}^{even} = 21 (= dim \ {\cal Z}$.
Then the general fiber of
$\tilde{det}{\cal Z}\rightarrow {\cal R}_{6+1}^{even}$
will be a point.

According to (N.3.4)(i),(ii),(iii), the variety
${\cal R}_{6+1}$ admits a natural surjective map

${\psi}: {\cal R}_{6+1}^{even} \rightarrow {\cal R}_6$
,
$(\tilde{\Delta}_6 \cup {\bf P}' \cup {\bf P}'')
\mapsto
(\tilde{\Delta}_6 \cup {\bf P}).$

As it follows from the preceding, the general fiber
$({\psi})^{-1}(\tilde{\Delta}_6 \cup {\bf P})$
is isomorphic to the $2$-dimensional set
${\cal S}$
of
$6$-tuples:

${\cal S} = {\cal S}({\Delta})$
$\{ (s)'_1,...,(s)'_6 \}
\in Symm^6(\tilde{\Delta}_6):
{\pi}((s)'_1 + ... + (s)'_6) \in
\mid {\cal O}_{{\Delta}_6}(1) \mid \}$.

It is known -- see [B3] -- that
${\cal S}$ is a union of two irreducible
$2$-dimensional components
${\cal S} = {\cal S}^+ \cup {\cal S}^-$,
dependent on the parity of the
elements of
$({\psi})^{-1}(\tilde{\Delta}_6 \cup {\bf P})$.
Therefore,
one of these two components,
say ${\cal S}^+$
must belong to ${\cal R}_{6+1}^{even}$,
the restriction
${\psi}^{even} := {\psi}_{\mid {\cal R}_{6+1}^{even}}:
{\cal R}_{6+1}^{even} \rightarrow {\cal R}_6$
is surjective,
and
the general fiber of $\psi$ is isomorphic to
${\cal S}^+$.
In particular,
$dim \ {\cal R}_{6+1} = dim \ {\cal R}_6 + 2 = 21$
{\bf q.e.d.}

{\bf (N.3.7)}
{\bf Corollary.}
{\sl
The natural map
${\psi}\circ \tilde{det}:{\cal Z} \rightarrow {\cal R}_6$
is surjective,
and the general fiber
of ${\psi}\circ \tilde{det}$,
over the covering
${\pi}:\tilde{\Delta} \rightarrow {\Delta}$,
is isomorphic to
one of the two irreducible components of
the surface
${\cal S}({\Delta}) =
{\pi}^*(\mid {\cal O}_{\Delta}(1) \mid )$.
}


{\bf (N.3.8)}
{\it The fiber of the intermediate jacobian map
$j:{\cal X}_{10}^{nodal} \rightarrow {\cal A}_{9}$.}

Let $j:{\cal X}_{10}^{nodal} = {\cal X}'' = {\cal Z}
\rightarrow {\cal A}_9$
be the intermediate jacobian map:
$j:Z \mapsto j(Z) = (J(Z),{\Theta})$,
where $J(Z)$ is the abelian part of the generalized
intermediate jacobian $J(Z)^{gen}$ of the singular
threefold $Z$. Clearly, $J(Z) = J(X)$, where
$X$ is the nodal Fano threefold of degree $10$
canonically birational to $Z$ -- see (N.3.2).

We shall identify ${\cal Z} = {\cal X}''$ and
${\cal X}_{10}^{nodal}$, as well their
intermediate jacobian maps $j$.

{\bf Proposition.}
{\sl
Let
$j:{\cal X}_{10}^{nodal} \rightarrow {\cal A}_9$
be the intermediate jacobian map, defined in the
general $X \in {\cal X}_{10}^{nodal}$ by
$\ j:X \mapsto (J(X),{\Theta})$
-- where $J(X)$ is the abelian part of the generalized
intermediate jacobian of $X$, and ${\Theta}$
is the theta-divisor of principal polarization
on $J(X)$.
Let ${\cal J}_{10}^{nodal}$ be the image of $j$,
let $j(X)$ be a general point of
${\cal J}_{10}^{nodal}$, and let
$j^{-1}(j(X)) \subset {\cal X}_{10}^{nodal}$
be the fiber of $j$ through $X$. Then
$j^{-1}(j(X)) = F \cup F'$,
where $F$ and $F'$ are irreducible
surfaces. Moreover, let $Z = X''$ be the
projection of $X$ from the node $o \in X$,
let $\Delta$ be the determinantal sextic of $Z$,
and let
${\pi}:\tilde{\Delta} \rightarrow {\Delta}$
be the induced double covering.
Then one of the components of the fiber, say $F$,
is isomorphic to one of the two irreducible
components, say ${\cal S}^+$, of the variety
${\cal S} = \{ D \in Symm^6(\tilde{\Delta}):
{\pi}_*(D) \in \mid {\cal O}_{\Delta}(1) \mid \}$.
Moreover, $F'$ is isomorphic to a similar
component ${\cal S}'^+$, obtained from
the pair
$(\tilde{\Delta}',{\Delta}')$
which is the unique involutive of
$(\tilde{\Delta},{\Delta})$
under the Dixon correspondence
(see [Ve])
in ${\cal R}_6$.
}

{\bf Proof.}
Let
$p:{\cal R}_6 \rightarrow {\cal A}_9$ be the Prym map
$p:(\tilde{\Delta},{\Delta}) \mapsto$
the p.p. Prym variety
$Prym(\tilde{\Delta},{\Delta}) = (P,\Xi )$.
As it follows from the recent result of Verra,
the fiber
$p^{-1}(p(\tilde{\Delta},{\Delta}))$
has exactly two elements -- the pair
$(\tilde{\Delta},{\Delta}) $, and the pair
$({\tilde{\Delta}}',{\Delta}')$, obtained from
$(\tilde{\Delta},{\Delta}) $ by the Dixon
correspondence for double coverings of plane sextics
-- see [Ve].

Let ${\cal S}'$ be the analog of the surface
${\cal S}$, for the pair
$(\tilde{\Delta}',{\Delta}')$.
As it follows from (N.3.6), one of the two
irreducible components of ${\cal S}'$, say
${\cal S}'^+$, is in $1$-by-$1$ correspondence
with the family of these
$Z \in {\cal Z}$ the determinantal pair of
which coincides with $(\tilde{\Delta}',{\Delta}')$.

The elements
$Z \in {\cal Z} \cong {\cal X}_{10}^{nodal}$,
which are in $1$-by-$1$ correspondence with the
elements of the two described $2$-dimensional
families ${\cal S}^+$ and ${\cal S}'^+$,
sweep out one entire fiber of $j$.
Indeed,
$j(Z) = p(\tilde{\Delta},{\Delta})$ =
$p(\tilde{\Delta}',{\Delta}')$,
for any such $Z$.  Moreover, the $1$-by-$1$
correspondence between ${\cal Z}$
and ${\cal R}_{6+1}^{even}$, the irreducibility
of the fiber ${\cal S}^+$
of the natural projection
${\cal R}_{6+1}^{even} \rightarrow {\cal R}_6$,
and the Verra's result that $deg(p) = 2$,
imply that the elements $F \cup F'$ sweep out
one complete fiber of $j$.
{\bf q.e.d.}


{\bf (N.3.9)}
{\bf Remark.}

Let $Z \in {\cal Z}$, and let $X$ be the
nodal Fano $3$ fold of degree $10$
for which $Z$ = $X''$ = the projection of $X$ from its
node.  Define $F(Z)$ to be this component of the
family of conics on $Z = X''$, the general
element of which is a projection of a conic on $X$.
The numerical definition of such a conic $q \subset Z$ is:
$deg(q) = 2$ and $(q.Q) = 0$, where $Q \subset Z$ is the
quadratic surface which lie on $Z$.
That is, $F(Z)$ is the projection of the Fano surface
${\cal F}(X)$ of conics on $X$.

Let $(\tilde{\Delta},{\Delta})$ be the determinantal
pair of $Z$, and let
${\cal S}^+ \subset Symm^6(\tilde{\Delta})$
be the ''even'' component of ${\cal S}$
parameterizing the elements of the component $F$
of $j^{-1}(j(Z))$ which pass through $Z$  -- see (N.3.6).

It can be seen that there exists a natural
birational isomorphism
${\phi}:{\cal F}(X) \rightarrow {\cal S}^+$
-- see [L], or (F.5.4).
Assuming the birational isomorphism ${\phi}$,
it follows that the Fano surfaces $F(Z)$,
of the elements $Z$
of the component $F$ (see (N.3.8), are birational
to each other (any $F(Z)$ is birational to  $F$).
Moreover, let $X$ be the unique nodal
Fano $3$-fold of degree $10$, for which
$Z = X''$ = the projection of $X$ from its node.

Since, by definition, the general conic of $F(Z)$
is a projection of a conic from the Fano surface
${\cal F}(X)$, this implies
that ${\cal F}(X)$ is birational to the component
$F$ of $j^{-1}(j(X))$.
(Here the families ${\cal X}_{10}^{nodal}$ and
${\cal Z}$, as well their intermediate jacobian
maps to ${\cal A}_9$, are identified.)

The same is true for the elements $Z'$ of the
residue component
$F'$ of the fiber $j^{-1}(j(Z))$, and for the
Dixon involutive pair
$(\tilde{\Delta}',{\Delta}')$.

\centerline{*  \  *  \  *}

In the next section,  we shall see that similar
results are true also for the general Fano threefold
of degree $10$ :
The fiber of the intermediate jacobian map
$j:{\cal X}_{10}$
is a union of two surfaces $F \cup \overline{F}$, s.t.
$F$ (resp. $\overline{F}$) is birational to the Fano
surface $F(X)$ of any threefold $X$ representing
a point of $F$ (resp. $\overline{F}$).


\newpage


\centerline{
{\bf F.
The fiber of the Griffiths intermediate jacobian map
$j:{\cal X}_{10} \rightarrow {\cal A}_{10}$.
}}

{\bf (F.0)}
In this section we describe the general fiber of the
Griffiths map on ${\cal X}_{10}$ -- see Theorem (F.6).

{\bf (F.1)}
{\bf Proposition.}
{\sl
Let
$X \in {\cal X}_{10}$
be general, let $q$ be a general conic on $X$, and let
$
{\alpha}_q : X \rightarrow X_q \in {\cal X}_{10}$
be the conic transformation of $X$ defined by $q$.
Then $X_q$ is not biregular to $X$,
i.e. ${\alpha}_q$ is not a birational automorphism of $X$.
}

{\bf Proof.}
Let ${\cal F} = {\cal F}(X)$ be the Fano surface
of conics on $X$
-- see (C.1),
and let $q \in {\cal F}$ be general.
According to
(C.6),
${\alpha}_q$
transforms the pair
$(X,q), q \in {\cal F}(X)$
to a pair
$(X_q,\overline{q}), \overline{q} \in {\cal F}(X_q)$.

Assume that $q$ is a general conic on $X$, and
$X_q = X$.

Let
$M = M_q \in \mid 2.H - 3.q \mid$
and
$\overline{M} =
M_{\overline{q}} \in \mid 2.H - 3.\overline{q} \mid$
be the exceptional divisors of
${\alpha}_q$ and ${\alpha}_{\overline{q}}$
-- see (C.3.1)-(C.3.2),
and let
${\sigma} = {\sigma}_q:X' \rightarrow X$
and
$\overline{\sigma} = {\sigma}_{\overline{q}}:
X^+ \rightarrow X_q = X$
be the blow-ups of $q$ and $\overline{q}$,
defined in the construction of the map
${\alpha}_q$
(see (C.2.1)).
Let
$Q' = {\sigma}^{-1}(q) \subset X'$
and
$\overline{Q}^+ = {\overline{\sigma}}^{-1}(\overline{q})
\subset X^+$
be the exceptional divisors of
${\sigma}$ and $\overline{\sigma}$.
Let $M' \subset X'$ be the proper preimage
of $M \subset X$, and let
$\overline{M}^+ \subset X^+$ be the proper
preimage of $M^+ \subset X_q = X$.

Let ${\rho}:X' \rightarrow X^+$
be the flop, defined by ${\alpha}_q$.
Then, by the involutive property (C.6) of ${\alpha}_q$,
the flop $\overline{\rho}:X^+ \rightarrow X'$,
defined by ${\alpha}_{\overline{q}}$,
coincides with ${\rho}^{-1}$.
Let
$M^+ \subset X^+$
be the proper ${\rho}$-image of $M'$ on $X^+$, and let
$\overline{M}' \subset X'$
be the proper ${\rho}$-preimage of $\overline{M}$ on $X'$.
Then, according to
(C.5.3) (, or (C.6.2)),
$M^+ = \overline{Q}^+$,
and
$\overline{M}' = Q'$.
There are $4$ possibilities:

(a). $\overline{q} \neq q$,
     and $\overline{q} \cap q = \oslash$;

(b). $\overline{q} \neq q$, but
     $\overline{q}$ intersects $q$ in a single point;

(c). $\overline{q} = \tilde{q}$ =
     the unique involutive conic of $q$
     -- see (C.3.3);

(d). $\overline{q} = q$.

In case (a), the coincidences
$M^+ = \overline{Q}^+$ (= the exceptional divisor of
the blow-up of $\overline{q}$),
and
$\overline{M}' = Q'$  (= the exceptional divisor of
the blow-up of $q$),
imply that  the divisors
$Q' \subset X'$ and
$M'$ (= the proper ${\rho}$-preimage, in $X'$,
of $M^+ = \overline{Q}^+$)
are disjoint, which contradicts the definition
of the divisor $M \in \mid 2.H - 3.q \mid$.

In cases (b) and (c), the same coincidences as in
(a), imply that $Q'$ and $M'$ intersect
each other in a union of fibers of the natural
projections
${\sigma}:Q' \rightarrow q$
and
$\overline{\sigma}\circ {\rho}:
M' \rightarrow M^+ = \overline{Q} \rightarrow \overline{q}$.
This is also impossible
by the definitions of $M$ and $\overline{M}$:
\ the surface $M \subset X$ coincides with the union of
rational quartic curves $C \subset X$ which intersect
$q$ in $3$ points
-- see (C.3.1)-(C.3.2).
As it follows from the construction of
${\alpha}_q$
-- see (C.6.2),
the last implies that the intersection cycle
$M' \cap Q'$ is a $3$-section of
${\sigma}: Q' \rightarrow q$ -- contradiction.

In case (d), the conics $q$ and $\overline{q}$
are coincident.
Therefore
$X'$ and $X^+$
can be identified, i.e. $X' = X^+$.
In particular
${\sigma} = \overline{\sigma}$, and
$\overline{M}' = Q' = \overline{Q}^+ = M^+$
(the left and the right coincidence -- by $X_q = X$,
and $Q' = \overline{Q}^+$ -- by $q = \overline{q}$
and $X' = X^+$).
However, $M' \neq M^+$, since the flop
${\rho}:M' \rightarrow M^+$
is not a biregular isomorphism --
for example, ${\rho}$ blows-down
the exceptional curves $l'_i, q'_j , \tilde{q}' \subset M'$
-- see (C.1.1), (C.2.2), (C.6.2).
Therefore, case (d) is also impossible.
{\bf q.e.d.}

{\bf (F.2)}
{\bf Proposition.}
{\sl
Let $q_1$ and $q_2$ be two general conics on $X$,
s.t. $q_1 \neq q_2$, and let
${\alpha}_{q_1}:X \rightarrow X_{q_1}$
and
${\alpha}_{q_2}:X \rightarrow X_{q_2}$
be the conic transformations of $X$,
defined by $q_1$ and $q_2$.
Then $X_{q_1} \neq X_{q_2}$.
}

{\bf Proof:}
(see the proof of (3.1).
Just as in (F.1), the proof of (F.2) can be based
on the intersection properties of exceptional divisors
''of type $M$'' and ''of type $Q$'':

Let
${\alpha}_{q_1}:X \rightarrow X_{q_1}$,
and
${\alpha}_{q_2}:X \rightarrow X_{q_2}$
be the conic transformations of $X$,
defined by $q_1$ and $q_2$.
By assumption $q_1 \neq q_2$, but
$X_{q_1} = X_{q_2} =:\overline{X}$.
Let
$\overline{q_1} \subset \overline{X}$
and
$\overline{q_2} \subset \overline{X}$
be the conics on
$Y = X_{q_1} = X_{q_2}$,
defined by the pairs
$(X,q_1)$ and $(X,q_2)$ -- see (C.6).
Then, by the involutive property (C.6)
of ${\alpha}$,
$\ {\alpha}_{\overline{q}}(\overline{X}) = X$,
${\alpha}_{\overline{q}}(\overline{X}) = X$,
and the pairs $(X,q_1)$ and $(X,q_2)$
are defined respectively as conic transformations
of the pairs
$(\overline{X}, \overline{q_1})$
and
$(\overline{X},\overline{q_2})$ (see ........).

Let $H$ and $\overline{H}$ be the hyperplane sections
of $X$ and $\overline{X}$.  Let
$M_1 \in \mid 2.H - 3.q_1 \mid$,
$M_2 \in \mid 2.H - 3.q_2 \mid$,
$\overline{M}_1 \in
\mid 2.\overline{H} - 3.\overline{q_1} \mid$,
and
$\overline{M} \in
\mid 2.\overline{H} - 3.\overline{q_2} \mid$
be the exceptional divisors of the corresponding
conic transformations,
defined by $q_1$, $q_2$, $\overline{q_1}$,
and $\overline{q_2}$.
Let
${\sigma}_1:X'_1 \rightarrow X$,
${\sigma}_2:X'_2 \rightarrow X$,
$\overline{\sigma}_1:X^+_1 \rightarrow \overline{X}$,
and
$\overline{\sigma}_2:X^+_2 \rightarrow \overline{X}$
be the blow-ups of
$q_1 , q_2$,
$\overline{q_1}$,
and
$\overline{q_2}$,
and let
$Q'_1, Q'_2$,
$\overline{Q}^+_1$,
and
$\overline{Q}^+_2$
be the exceptional divisors
of the corresponding blow-ups of
$q_1,q_2$,
$\overline{q_1}$,
and $\overline{q_2}$.
Just as in the proof of (F.1),
the coincidence
$X_{q_1} = X_{q_2} (= \overline{X})$
imply that the proper preimages
of $M_i$ and $\overline{Q}_i$ coincide
on $X'_1$, $X'_2$, $X^+_1$, and $X^+_2$,
$i =1,2$.
Now, corresponding to the cases

(a). $q_1$ and $q_2$ are disjoint,

(b). $q_1 \cap q_2$ = (a point),

(c). $q_2 = \tilde{q_1}$ = the unique the involutive of
     the conic $q_1$,

these coincidences
imply that the preimages of
the surfaces
$M_1$ and $M_2$ must intersect each other
along a (possible empty) union
of fibers -- rational quartic curves
intersecting $q_1$ and $q_2$
in $3$ points. As in (3.1), this is impossible.
Indeed, let $C \subset M_1$ be a general
fiber of $M_1$ -- a rational quartic
curve on $X$, intersecting $q_1$ in $3$ points.
Then $C$ intersects the quadratic section
$M_2 \in \mid 2.H - 3.q_2 \mid \subset \mid 2.H \mid$,
in $8 = 2.deg(C)$ points.  Therefore
the intersection
$M_1 \cap M_2$ is a $8$-section
(not a union of fibers), on any of the
mentioned proper preimages of $M_1$ and $M_2$ --
contradiction.  Therefore $X_p \neq X_q$.
{\bf q.e.d.}

{\bf (F.2.1)}
{\bf Corollary.}
{\sl
Let $X \in {\cal X}_{10}$ be general,
and let ${\cal F} = {\cal F}(X)$ be the
Fano surface of conics on $X$.
Let
${\alpha}:{\cal F} \rightarrow {\cal X}_{10}$
be the rational map, defined by the closure
(in ${\cal F} \times {\cal X}_{10}$)
of the graph-correspondence:

${\Sigma} =
\{ (q,Z) \in {\cal F} \times {\cal X}_{10}:
q \mbox{ is a general conic on } X,
\mbox{ and } Z = X_q \}$.

Then ${\alpha}$ defines a birational
isomorphism
${\alpha}:{\cal F} \rightarrow {\alpha}({\cal F})$
onto its image
${\alpha}({\cal F}) \subset {\cal X}_{10}$.
}

{\bf (F.3)}
{\bf Proposition.}
{\sl
Let $X \in {\cal X}_{10}$ be general,
and let
${\alpha}:{\cal F} \rightarrow {\alpha}({\cal F})
\subset {\cal X}_{10}$
be the map defined in (F.1.1).

Let
$j:{\cal X}_{10} \rightarrow {\cal A}_{10}$
be the intermediate jacobian map on ${\cal X}_{10}$,
defined (in the general point $X \in {\cal X}_{10}$) by

$j:X \mapsto \mbox{ the p.p. intermediate jacobian }
(J(X),{\Theta}) \mbox{ of } X .$
Let ${\cal J}_{10} \subset {\cal A}_{10}$
be the image of $j$,
and let
$j^{-1}(j(X)) \subset {\cal X}_{10}$
be the fiber of $j$ over the (general) point
$j(X) \in {\cal J}_{10}$.

Then the surface ${\alpha}({\cal F}(X))$
is an irreducible component of maximal dimension
of the fiber $j^{-1}(j(X))$.
}

{\bf Proof.}
For the general $X in {\cal X}_{10}$,
${\cal F}(X)$ is a smooth irreducible surface
-- see  [L], or [Ili1].
According to (F.2.1),
${\alpha}({\cal F}(X)$ is birational to ${\cal F}$.
Therefore ${\alpha}({\cal F}(X)$ is an irreducible
surface in ${\cal X}_{10}$.

The inclusion
${\alpha}({\cal F}) \subset j^{-1}(j(X))$
will follow from the coincidence of the
p.p. intermediate jacobians
$J(X)$ and $J(X_q)$, for the general
conic $q \subset X$.
This is obvious:
By construction -- see (C.6.2) --
the conic transformation
${\alpha}_q:X \rightarrow X_q$
is a product of ${\sigma}$-processes
centered on smooth rational curves;
and any such ${\sigma}$-process does
not change the intermediate jacobian
of the threefold;
see also (I.1).

In rests to see that
$dim \ j^{-1}(j(X)) = 2$.
\ Then ${\alpha}({\cal F}(X))$ will be a component
of maximal dimension of $j^{-1}(j(X))$.

Consider the restriction of $j$, on the
subvariety ${\cal X}_{10}^{nodal}$ of
nodal Fano $3$-folds of degree $10$.
Clearly, the general Lefschetz pencil
$\{ X(t) \subset W_5 : t \in {\bf P}^1 \}
\subset {\cal X}_{10}$
contains a finite number of nodal threefolds,
i.e.
$dim({\cal X}_{10}^{nodal}) =
dim({\cal X}_{10}) - 1 = 22 - 1 = 21$
(see also (N.3.8)).

Let
${\cal J}_{10}^{nodal} = j({\cal X}_{10}^{nodal})
\subset {\cal J}_{10} \subset {\cal A}_{10}$
be the subvariety of generalized intermediate
jacobians of nodal $X_{10}$.
If $X \in {\cal X}_{10}^{nodal}$ is general,
then the node $o$ of $X$ defines $j(X)$
uniquely as an extension of ${\bf C}^*$
by the abelian part
$(J(X),{\Theta}) \in {\cal A}_{9}$.

On the one hand, according to (N.3.8),
the general fiber of the map
$j^{n}:{\cal X}_{10}^{nodal} \rightarrow {\cal A}_9$,
is a union of two irreducible surfaces;
in particular,
$dim \ j^{n}({\cal X}_{10}^{nodal}) = 19$.
Therefore,
$dim \ {\cal J}_{10} \ge
dim \ j^{n}({\cal X}_{10}^{nodal}) + 1 = 19 + 1 = 20$.
On the other hand,
$dim \ {\cal J}_{10} \le
dim \ {\cal X}_{10} - 2 = 20$,
since the general fiber
$j^{-1}(J(X))$, $X \in {\cal X}_{10}$,
contains the surface
${\alpha}({\cal F}(X))$.
Therefore
$dim \ {\cal J}_{10} = 20$,
and the components of maximal dimension,
of the general fiber of $j$, are of dimension $2$
{\bf q.e.d.}

The next if the analog of (F.1)-(F.2)-(F.2.1)-(F.3),
for line transformations. The proofs are
identical:

{\bf (F.4)}
{\bf Proposition.}
{\sl
Let $X \i {\cal X}_{10}$ be general,
and let $l_1$ and $l_2$ be a pair of two
general lines on $X$.

Let
${\beta}_{l_1} : X \rightarrow X_{l_1}$
and
${\beta}_{l_2} : X \rightarrow X_{l_2}$
be the transformations defined by $l_1$ and
$l_2$.
Then $X_{l_1} \neq X$, $X_{l_2} \neq X$;
and if $l_1 \neq l_2$ then
$X_{l_1} \neq X_{l_2}$.

Let ${\Gamma} = {\Gamma}(X)$ be the
curve of lines on $X$
-- see (L.1).
\ Let
${\beta}:{\Gamma} \rightarrow {\cal X}_{10}$
be the map defined by:
${\beta}:l \mapsto X_l$, for general $l \in {\Gamma}$
(see the definition of ${\alpha}$ in (F.2.1)).
Then ${\beta}$ defines an embedding
of ${\Gamma}$, in the fiber
$j^{-1}(j(X))$ of the intermediate jacobian map $j$,
over the point $j(X) \in {\cal J}_{10}$.
}

{\bf (F.5)}
{\bf Proposition.}
{\sl
Let $X \in {\cal X}_{10}$ be general.
Then the fiber $j^{-1}(j(X))$ has exactly
two irreducible components of maximal
dimension = $2$.

Let $F$ and $\overline{F}$ be
the closures, in ${\cal X}_{10}$, of
the two $2$-dimensional irreducible components
of $j^{-1}(j(X))$.
Then
$j^{-1}(j(X)) = F \cup \overline{F}$.

Denote by $F$ the component of
$j^{-1}(j(X))$ defined by the
condition $X \in F$.
Let ${\cal F}(X)$ and ${\Gamma}(X)$
be the Fano surface and the curve of lines
on $X$.
Then
${\alpha}({\cal F}(X)) = F$,
and
${\beta}({\Gamma}(X)) \subset \overline{F}$,
i.e. the conic transformations
do not interchange the components of $j$,
while the line transformations interchange
the components of $j$.
}

{\bf Proof.}
To begin with, we shall see that
${\alpha}({\cal F}(X)$
and ${\beta}({\Gamma}(X)$
lie in different components f
$j^{-1}(j(X))$.
Assume, on the contrary, that
${\alpha}({\cal F}(X)$
and ${\beta}({\Gamma}(X)$
lie on the same component of
$j^{-1}(j(X))$,
say $F$.
Then the closed set ${\alpha}({\cal F}(X))$
coincides with $F$
-- see (5.3).
In particular,
${\beta}({\Gamma}(X)) \subset {\alpha}({\cal F}(X))$.

Turning back to the definition of the conic
transformation in section C,
one can see that ${\alpha}_q$ can be defined,
and $X_q \in {\cal X}_{10}$,
for any ${\tau}$-conic $q \subset X$
(see the description of $M$ in (C.3.1)-(C.3.2)) --
including the reducible ${\tau}$-conics
$q = l + m$, where the intersecting lines
have non-coincident centers --
see Remark (F.5.1).

In particular, we obtain that
$X_q = {\alpha}_q(X)$ is a well-defined Fano threefold
of degree $10$, for any
${\tau}$-conic $q \subset X$.
Indeed, the only exceptions of ${\tau}$-conics,
where ${\alpha}$ may be not well-defined,
arise from the finite number of involutive pairs
of reducible $\tau$-conics on the general $X$,
which form a subset of codimension $2$
in the open set of ${\tau}$-conics on $X$.
By the Hartogs's theorem, the map ${\alpha}$
can be completed to a regular map also
in these ${\tau}$-conics (all they are, in general,
disjoint from the curve of ${\rho}$-conics,
where ${\alpha}$ has not been described explicitly.

Therefore, the only elements of the Zariski closure
of ${\alpha}(X)$,
which may lie on the boundary of the moduli space
${\cal X}_{10}$, must correspond to
${\sigma}$-conics, or --
to the unique ${\rho}$-conic on $X$.

Since the set of ${\sigma}$-conics on $X$ is a smooth
rational curve on ${\cal F}(X)$
(see (C.1.1)),
the isomorphic image
${\beta}({\Gamma}(X))$ of the non-rational
curve ${\Gamma}(X)$
(see (F.4))
does not lie entirely outside
the ${\alpha}$-image of the open
set of ${\tau}$ conics on $X$.
Indeed,
by the decomposition of birational isomorphisms
in $dim = 2$ into ${\sigma}$-processes
and their inverses,
the surface ${\alpha}({\cal F}(X))$,
which is birational to the smooth
surface ${\cal F}$,
must be biregular to ${\cal F}$,
except at most in a union of rational curves.

Therefore, the general line transformation
$X_l$ of $X$ must coincide with some
conic transformation $X_q$ of $X$,
where $q$ is a ${\tau}$-conic on $X$
-- which is impossible by the same arguments
as these in the proof of (F.2).
Therefore
the curve
${\beta}({\Gamma}(X))$ does not lie in the
component ${\alpha}({\cal F}(X))$
of $j^{-1}(j(X))$.

The next is to see that ${\alpha}({\cal F}(X)$
must coincide with the component $F$ of
$j^{-1}(j(X)) = F \cup \overline{F}$,
defined by $X \in F$.
Assume, on the contrary, that
${\alpha}({\cal F}(X) = \overline{F}$.
According to the preceding, this
implies that
${\beta}({\Gamma}(X)) \subset F$.

Let $q = l + m$ be a general
reducible conic on $X$,
i.e. the components $l$ and $m$
of $q$ are two different lines,
intersecting each other in a single point $x \in X$.
Since $X$ is general, we can assume
that
$N_{l \mid X} \cong N_{m \mid X}
\cong {\cal O}(-1) \oplus {\cal O}$,
and the ``centers'' (see the proof of (L.3.1))
of $l$ and $m$
do not coincide.

Therefore $q = l + m$ is a conic of
${\tau}$-type
(see (C.1.0)).
Moreover, by the general choice of $X$,
one can suppose that the unique involutive
conic $p$, of the ${\tau}$-conic
$q = l + m$, is smooth.

By Remark (F.5.1),
the transformation
${\alpha}$ is defined on $l+m$,
and the threefold
$X_{l+m} = {\alpha}_{l+m}(X)$
is an element of ${\cal X}_{10}$.

Let ${\beta}_l:X \rightarrow X_l$
be the line transformation of $X$
defined by $l$,
and let $\overline{l} \subset X_l$
be the line defined by the pair $(X,l)$
-- see (L.6).
As it follows from the description of
${\beta}_l$
-- see e.g. (L.6.2),
$\ {\beta}_l$ transforms $m$
to a line $m_1 \subset X_l$ which
intersects $\overline{l}$.
The pair $(X_l,m_1)$
defines the line transformation
${\beta}_{m_1}:X_l \rightarrow X_{l,m_1}$,
and the line
$\overline{m_1} \subset X_{l,m_1}$.
Similarly, the proper ${\beta}_{m_1}$-image of
$\overline{l}$ is a line
$\overline{l}_1 \subset X_{l,m_1}$,
which intersects $\overline{m_1}$.

On the one hand, the birational map
${\beta}_{m_1} \circ {\beta_l}:
X \rightarrow X_{l,m_1}$
is a composition of two line transformations,
which -- by assumption -- leaves the image
$X_{l,m_1}$ in the component $F$ which
contains $X$.
On the other hand, the composition
${\beta}_{m_1} \circ {\beta_l}$,
and the map
${\alpha}_{l+m}$
coincide.
In fact, by the same argument as in
[Isk3,(4.6.1),(4.6.2)]
(where conic and line transformations are defined for
the intersection $V_{2.3}$ of a quadric and a cubic in
${\bf P}^5$),
the conic transformation ${\alpha}_{l+m}$,
and the composition
of the line transformations
${\beta}_{m_1} \circ {\beta}_l$,
differs from each other by a biregular
isomorphism of $X$.  Since the general $X$
has no biregular automorphisms
-- see Lemma (A.6.1),
this proves the coincidence.
In particular,
$X_{l,m_1} = X_{l+m} .$

Now, according to the assumption,
the line transformations
${\beta}_l$ and ${\beta}_{m_1}$
leave the images $X_{l}$
and $X_{l,m_1}$ in the component $F$,
while $X_{l+m} \in \overline{F}$.
However $X_{l+m} = X_{l,m_1}$ --
contradiction.

It rest to see that the fiber of $j$, through
the general $X \in {\cal X}_{10}$, cannot have other
components, besides the two already known components $F$
and $\overline{F}$ of the birational orbit of $X$.

As it follows from section N ,
the fiber of $j$, over the generalized intermediate
jacobian of the general nodal $X_{10}$, is a union
of two irreducible $2$-dimensional families
of nodal $X_{10}$-s  (see e.g. Proposition (N.3.8) --
formulated in terms of the $9$-dimensional abelian part
of the $10$-dimensional generalized jacobian of the nodal
$X_{10}$).
Therefore, the $2$-dimensional general fiber
$j^{-1}(j(X))$
must be a union of two irreducible surfaces.
However, we already
know these surfaces -- they are, in fact,
the already known components $F$ and $\overline{F}$
of the birational orbit  $Orb_{bir}(X)$.
In particular, all the points of $Orb_{bir}(X)$
must lie on $F \cup \overline{F}$.
{\bf q.e.d.}


{\bf (F.5.1)}
{\bf Remark.}
{\it Description of the conic transformation
in a general reducible conic.}

It is not hard to find the
degeneration of the
conic transformation for such
degenerate conics.
We shall describe it
for a reducible conic $l + m$
for which the unique involutive conic
$p = \tilde{l + m}$ is smooth.
By the general choice of $X$,
this is the general case of
a reducible conic on $X .$

To the first,
one have to replace $X$ with
$\tilde{X}$ = (the blow-up of $X$ in $x$).
Then, the proper preimages
$\tilde{l} \subset \tilde{X}$
and
$\tilde{m} \subset \tilde{X}$,
of $l$ and $m$,
become disjoint, and one can blow-up
$\tilde{l} \subset \tilde{X}$
and
$\tilde{m} \subset \tilde{X}$
simultaneously to a threefold $\tilde{X}'$.

Then one has to perform  all the
${\sigma}$-processes over the
proper preimages of the contractible lines
and conics for the projection from $l + m$,
to obtain the threefold $\tilde{X}^+$
(which is the analog of the
flop-ed threefold $X^+$
for a smooth $\tau$-conic
-- see (C.6.2) and (L.6.2).
Then the proper images
$\tilde{M_l}^+$ and $\tilde{M_m}^+$,
on $\tilde{X}^+$
-- of the unique divisors
$M_l \in \mid H - 2.l \mid$
and
$M_m \in \mid H - 2.m \mid$ --
become disjoint and
contractible along their rulings.
This way, one obtains the threefold
$\tilde{X}_{l + m}$, which
(by symmetry -- just as $\tilde{X}$)
contains a plane $P^+$ -- the proper
image of the exceptional divisor
which lie over the unique involutive conic
$p$ of $l + m$.
The base curves
$\tilde{l}^+$ and $\tilde{m}^+$
of the obtained birational morphism
$\tilde{X}^+ \rightarrow  \tilde{X}_{l + m}$
are disjoint, and each of these two curves
intersects $P^+$ in a single point.
The plane $P^+$
defines a contraction
$\tilde{X}_{l + m} \rightarrow X_{l + m}$,
and the images $\overline{l}$ and
$\overline{m}$, of the two base curves,
joint each other in the point-image
$\overline{x}$ of $P^+$.
The reducible conic
$\overline{l} + \overline{m} \subset X_{l + m}$
is the counterpart of the conic
$l + m \subset X$, under the degenerate
conic transformation ${\alpha}_{l + m}$.


{\bf (F.5.2)}
{\bf Corollary.}
{\sl
In the notation of (F.5),
let $j^{-1}(J(X)) = F + \overline{F}$
be the decomposition of the general
fiber of $j$ into irreducible components.
In particular
$X \in F = {\alpha}({\cal F}(X))$,
where ${\cal F}(X)$ is the Fano surface of conics
on $X$.

Let $Z \in {\cal X}_{10}$ be another threefold
which belong to the component $F$.
Assume that $Z$ is otherwise general, and let
${\cal F}(Z)$ be the Fano surface of conics
on $Z$.  Then
${\cal F}(Z)$ is birational to ${\cal F}(X)$.
}

{\bf Proof.}
Since $Z \in F$ is general, and
$F = {\alpha}({\cal F}(X))$,
there exists a unique
conic $q \in {\cal F}(X)$
such that $Z = X_q = {\alpha}_q(X)$
(see (F.5) and (F.2)).
Since $X$ and $Z = X_q$ belong to the
same component $F$ of
$j^{-1}(j(X)) = j^{-1}(j(X_q))$,
both
${\cal F}(X)$ and ${\cal F}(X_q))$
are birational to $F$.
This proves (F.5.2).

{\bf (F.5.3)}
{\it The birational isomorphism
${\phi}_q:{\cal F}(X) \rightarrow {\cal F}(X_q)$.}

One can find directly a natural birational
isomorphism
${\phi}_q:{\cal F}(X) \rightarrow {\cal F}(X_q)$.

We shall define ${\phi}_q$ in the general conic
$p \in {\cal F}(X)$.
Since $p$ is general, one can assume that
$q$ and $p$ are disjoint,
and the planes $Span(q)$ and $Span(p)$
do not intersect each other in
${\bf P}^7 = Span(X)$.
Let
${\bf P}^5(p,q) := Span(p \cup q)$.
Since $X$ is a Fano threefold of
index $1$, and of degree $10$,
the intersection
${\bf P}^5(p,q) \cap X$ is a canonical
curve of degree $10$
(and of arithmetical genus $6$), which
contains the conics $p$ and $q$ as components.
The general choice of the pair $(p,q)$,
and the Riemann-Roch formulae for
singular curves, imply that the residue
curve
$C(p,q) = {\bf P}^5 \ \cap X - p - q$
\ is a smooth rational sextic which intersects
each of $p$ and $q$ in $4$ points.
Since $q$ is fixed, the choice of the general
conic $p$ defines uniquely $C(p,q)$.
Now, observe that the conic transformation
${\beta}_q:X \rightarrow X_q$
sends $C(p,q)$ to a conic $r(p,q) \subset X_q$.
In fact, is enough to compute $deg \ r(p,q)$:

Since ${\beta}_q$ is defined by the system
$\mid 3.H - 4.q \mid$,
then
$deg \ r(p,q) = deg \ {\beta}_q(C(p,q))$ =
$3.deg \ C(p,q) - 4.\#(C(p,q) \cap q)$ =
$3.6 - 4.4 = 2$, i.e. $r(p,q)$ is a conic.
Let, moreover,  $\overline{q} \subset X_q$
be the conic defined by the pair $(X,q)$.
It can be easily seen that if $p$ and $q$
are disjoint, then the conics $r(p,q)$
and $\overline{q}$ are disjoint.

For the general conic $p$,
we define the birational map
${\phi}_q: {\cal F}(X) \rightarrow {\cal F}(X_q)$
by
${\phi}:p \mapsto r(p,q)$.

{\bf (F.5.4)}
{\bf Remark.}

Let $\{ X(t) \subset W_5 : t \in {\bf P}^1 \}$
be a general Lefschetz pencil
in ${\cal X}_{10}$, and let $t = 0$
be one of the finite number of values
of $t$ for which $X(t)$
acquires a node.

Let
$(\tilde{\Delta},{\Delta})$,
and
${\tilde{\Delta}}', {\Delta}')$
be the determinantal pair,
and its Dixon involutive pair,
defined by $X(0)$
\ (see [Ve], or (N.3.8)-(N.3.9)).

According to Proposition (N.3.8),
the fiber of the intermediate
jacobian map for nodal $X_{10}$,
is isomorphic to the union of the
two special surfaces
${\cal S}^+ \subset Symm^6(\tilde{\Delta})$,
and
${\cal S}'^+ \subset Symm^6({\tilde{\Delta}}')$
-- see also the proof of (N.3.6).
Moreover, since the Lefschetz pencil
defines a projective deformation of
$X(t) \subset W_5 \subset {\bf P}^7$,
the Fano surface of conics on ${\cal F}(X(t))$
degenerates to the Fano surface
${\cal F}(X(0))$.
As it follows from (F.5),
${\cal F}(X(0))$ has to be birational to
the component $F(0)$, of the fiber
$j^{-1}j(X(0)) = F(0) \cup \overline{F}(0)$.
Here the intermediate jacobian map $j$
is regarded equally -- as taking
the abelian part of the generalized
Jacobian of $X(0)$, or -- as taking
the generalized jacobian of $X(0)$.
Since
$j^{-1}(j(X(0)) \cong
{\cal S}^+ \cup {\cal S}'^+$,
the non-ordered pairs
$(F(0), \overline{F}(0))$
and
$({\cal S}^+,{{\cal S}'^+})$
must coincide.
In particular, the Fano surface
${\cal F}(X(0))$ must be birational
to at least one of
${\cal S}^+$ and ${\cal S}'^+$.
This is consistent with Remark (N.3.9),
declaring the birational isomorphism
${\cal F}(X(0)) \rightarrow {\cal S}^+$
\ (see also (N.3.9)).

\newpage

{\bf (F.6)}
{\it The fiber of the Griffiths map for non-hyperelliptic
Fano threefolds of genus $6$.}

The following summarizes the main results in the paper:

{\sc Theorem.}
\ {\sl
Let ${\cal X}_{10}$ be the
moduli space of non-hyperelliptic Fano threefolds
of genus $g = 6$
(and -- of degree $10 = 2g-2$).
Let $j:{\cal X}_{10} \rightarrow {\cal A}_{10}$
be the Griffiths period map on ${\cal X}_{10}$,
sending the general $X \in {\cal X}_{10}$
to its $10$-dimensional p.p.intermediate
jacobian $j(X) = (J(X),{\Theta})$,
and let ${\cal J}_{10} \subset {\cal A}_{10}$
be the image of $j$.

Then
$dim \ {\cal J}_{10} = 20$ (section N),
and ${\cal X}_{10}$ is birational to
the $22$-dimensional orbifold
$\mid {\cal O}_W(2) \mid / Aut(W)$ --
where $W \subset {\bf P}^7$ is the
prehomogeneous Fano fourfold of
degree $5$, and $Aut(W)$ is its
$8$-dimensional automorphism
group (section A).
The fiber
$J^{-1}(J,{\Theta}) \subset {\cal X}_{10}$,
over the general
$(J,{\Theta}) \in {\cal J}_{10}$,
is a union of two irreducible
families
of threefolds, and:

{\bf (1).}
If $X$ is any threefold of $j^{-1}(J,{\Theta})$,
then the fiber
$j^{-1}(J,{\Theta})$
is isomorphic to the birational
orbit
$Orb_{bir}(X)$
(= the set of all birational images of $X$,
in the class ${\cal X}_{10}$) -- see (I.1)(*).
In particular, all the threefolds
of the fiber
are birational to each other.

{\bf (2).}
Let $F$ be one of the components
of $j^{-1}(J,{\Theta})$,
and let $\overline{F}$ be the
residue component.

Then, for the general pair
$(X_1,X_2) \in F \times F$,
there exists a unique conic
$q_1 \in X_1$,
and a unique conic
$q_2 \in X_2$,
such that the pairs
$(X_1,q_1)$ and $(X_2,q_2)$
are obtained from each other by the
involutive conic transformation
${\alpha}$ on ${\cal X}_{10}$.
In particular,
$X_j$ is a
birational image of $X_i$,
under the
quadruple cubic projection from $q_i$,
defined by the noncomplete linear system
$\mid {\cal O}_{X_i}(3) - 4.q_i \mid$,
where $(i,j) = (1,2), (2,1)$
-- see (C.6), (F.1)-(F.3), (F.5).

Moreover, if
$X \in F$,
then the general line
$l \subset X$ defines a birational
isomorphism
${\beta}_l:X \rightarrow X_l$,
where $X_l \in \overline{F}$
(see (C.6) and (F.5)) --
i.e. -- the line transformations
${\beta}_l$, $l - \mbox{ a line on } X$,
interchange the components $F$
and $\overline{F}$ of the fiber of the
Griffiths map.
}


\newpage


\centerline{\sc References}

\bigskip

[B1]  A.Beauville, {\it Vari\'et\'es de Prym et jacobienne
interm\'ediaires.}
Ann. de l'AENS, 4 ser.,10 (1977), 149-196.

\smallskip

[B2]  A.Beauville, {\it Les singularit\'es du diviseur theta
de la jacobienne interm\'ediaire de l'hypersurface cubique
dans ${\bf P}^{4}$.}
Lect. Notes in Math., Vol. 947 (1982), 190-208.

\smallskip

[B3]  A.Beauville, {\it Sous-vari\'et'es speciales des
vari\'et\'es de Prym.}
Comp. Math., 45 (1982), 357-383.

\smallskip

[C1]  H.Clemens, {\it Double solids.}
Adv. Math., 47:2 (1983), 107-230.

\smallskip

[C2]  H.Clemens, {\it The Quartic Double Solid Revisited.}
Proc. Symp.in Pure Math., Vol. 53 (1991), 89-101.

\smallskip

[C3]  H.Clemens, {\it The fiber of the Prym map and the
period map for double solids, as given by Ron Donagi.}
Univ. of Utah - Preprint.

\smallskip

[CG]  H.Clemens, P.Griffiths, {\it The intermediate jacobian
of the cubic threefold.}
Ann. Math., 95 (1972), 281-356.

\smallskip

[D1]  O.Debarre, {\it Le th\'eor\`eme de Torelli pour les
intersections de trois quadriques.}
Inv. Math., 95 (1989), 507-528.

\smallskip

[D2]  O.Debarre, {\it Sur le theoreme de Torelli pour les
solides doubles quartiques.}
Comp. Math., 73 (1990), 161-187.

\smallskip

[F]  M.Furushima, {\it Mukai-Umemura's example of
the Fano threefold with genus $12$ as a compactification
of ${\bf C}^3$.}
Nagoya Math. J., 127 (1992), 145-165.

\smallskip

[FS]  R.Friedman, R.Smith, {\it Degenerations of Prym
varieties and intersections of three quadrics.}
Inv. Math., 85 (1986), 615-635.

\smallskip

[GH]  P.Griffiths, J.Harris, {\it Principles of algebraic
geometry.}
J.Willey \& , New York (1978).

\smallskip

[Gua]  D.Guan, {\it Toward a classification of almost
homogeneous spaces.}
Princeton Univ., Sept. 1994 -- preprint.

\smallskip

[Gus]  N.Gushel, {\it On Fano varieties of genus $6$.}
Math. USSR - Izv., 21 (1983), 445-459.

\smallskip

[H1]  S.Hashin, {\it Birational automorphisms of the
double Veronese cone of dimension three.}
Vestnik MGU, Ser.I: Mat. - Mech., No.1 (1984), 13-16
(in Russian); translated in: Moscow Univ. Math. Bull.

\smallskip

[H2]  S.Hashin, {\it Birational automorphisms of Fano
manifolds of index $1$ and of degree $10$}
(announcement).
Proc. Conference ``Algebra, Logics, and Number Theory'',
Feb. - March 1985,
Moscow Univ. Publ. House (1986), ed. Yu.Manin \&
L.Skornjakov, 84-86 (in Russian).

\smallskip

[Isk1]  V.Iskovskikh, {\it Fano threefolds I.}
Math. USSR - Izv., 11:3 (1977), 485-527.

\smallskip

[Isk2]  V.Iskovskikh, {\it Fano threefolds II.}
Math. USSR - Izv., 12:3 (1978), 469-506.

\smallskip

[Isk3]  V.Iskovskikh, {\it Birational automorphisms
of three-dimensional algebraic varieties.}
J. of Soviet Math., Vol.13:6 (1980), 815-868.

\smallskip

[Isk4]  V.Iskovskikh, {\it On the rationality problem for
conic bundles.}
Duke Math. J., 54:2 (1987), 271-294.

\smallskip

[Isk5]  V.Iskovskikh, {\it Algebraic threefolds with special
regard to the problem of rationality.}
Proc.Int.Congress of Math., Aug.16-24 (1983), Warszawa,
733-746.

\smallskip

[Isk6]  V.Iskovskikh, {\it Lectures on algebraic manifods:
Fano manifolds.}
Moscow Univ. Publ. House (1988).

\smallskip

[Isk7]  V.Iskovskikh, {\it Double projection from a line
on three-dimensional Fano varieties of the first specie.}
Matem. Sbornik, Vol.180:2 (1989), 260-278.
(in Russian); translated in: Math. USSR - Sbornik.

\smallskip

[IM]  V.Iskovskikh, Yu.Manin {\it Three-dimensional quartics
and counterexamples to the L\"uroth problem.}
Math. USSR - Sbornik, 15 (1971), 141-166.

\smallskip

[Ili1]  A.Iliev, {\it Geometry of the Fano threefold of
degree $10$ of the first type.}
Contemp. Math., Vol.136 (1992), 209-254.

\smallskip

[Ili2]  A.Iliev, {\it Lines on the Gushel threefold.}
Indag. Math., N.S., 5(3) (1994), 307-320.

\smallskip

[Ili3]  A.Iliev, {\it The Fano surface of the Gushel threefold.}
Comp. Math., 94 (1994), 81-107.

\smallskip

[Ili4]  A.Iliev, {\it The theta divisor of the bidegree (2,2)
threefold in ${\bf P}^{2} \times {\bf P}^{2}$.}
preprint (1994).

\smallskip

[K]  J.Kollar, {\it Flops.}
Nagoya Math. J., 113 (1989), 14-36.

\smallskip

[L]  D.Logachev, {\it Abel-Jacobi isogeny for the Fano threefold
of genus $6$.}
Yaroslavl' Gos. Ped. Inst./ K.D.Ushinskij (1982) -- Thesis
(in Russian); \ or:
Coll. -- ``Sb.Constructive Algebraic Geometry'',
vol. 200 -- Yaroslavl' (1982), 67-76 (in Russian).

\smallskip

[Mo]  S.Mori, {\it Threefolds whose canonical bundles are
not numerically effective.}
Ann. Math., 116 (1982), 133-176.

\smallskip

[Mu]  S.Mukai, {\it Curves, $K3$ surfaces, and Fano $3$-folds
of genus $\le 10$.}
in: ``Alg. Geom. and Comm. Alg. in Honor of M.Nagata'',
(1988), 357-377 -- Kinokuniya, Tokyo.

\smallskip

[MU]  S.Mukai, H.Umemura, {\it Minimal rational threefolds.}
L. N. in Math., Vol.1016 (1983), 490-518.

\smallskip

[Pu]  P.J.Puts, {\it On some Fano threefolds that are sections
of Grassmannians.}
Indag. Math., 44 (1982), 77-90.

\smallskip

[PS1]  T.Peternell, M.Schneider, {\it Compactifications
of ${\bf C}^3$ I.}
Math. Ann., 280 (1988), 129-146.

\smallskip

[PS2]  T.Peternell, M.Schneider, {\it Compactifications
of ${\bf C}^3$ II.}
Math. Ann., 283 (1989), 121-137.

\smallskip

[Pr1]  Yu.Prokhorov, {\it Geometric properties
of Fano manifolds.}
Moscow Univ., Dept. Math. \& Mech.
(1990) -- Thesis (in Russian).

\smallskip

[Pr2]  Yu.Prokhorov, {\it Compactifications of ${\bf C}^4$
of index $3$}.
Proc. $8^{th}$ Alg. Geom. Conf. -- Yaroslavl' (1992),
Ser. Aspects of Math., Publ.: Steklov Inst. Math.,
eds. A.Tikhomirov \& A.Tyurin.

\smallskip

[Ti1]  A.Tikhomirov, {\it Geometry of the Fano surface
of the double space ${\bf P}^3$ ramified in a quartic.}
Izv. AN SSSR - Ser. Mat., 44:2 (1980), 415-442
(in Russian); translated in: Math. USSR - Izv.

\smallskip

[Ti2]  A.Tikhomirov, {\it The middle jacobian of the double
space ${\bf P}^3$ ramified in a quartic.}
Izv.AN SSSR - Ser.Mat., 44:6 (1980), 1329-1377
(in Russian); translated in: Math. USSR - Izv.

\smallskip

[Ti3]  A.Tikhomirov, {\it The Abel-Jacobi map of sextics of genus
$3$ on double spaces of ${\bf P}^{3}$ of index two}.
Soviet Math. Dokl., Vol. 33:1 (1986), 204-206.

\smallskip

[Tju]  A.Tjurin (A.Tyurin),
{\it The Middle Jacobian of Three-Dimensional Varieties.}
J. of Soviet Math., Vol.13:6 (1980), 707-744.

\smallskip

[Tr]  S.Tregub, {\it Construction of a birational isomorphism
of three-dimensional cibic and a vaiety of Fano of the first
specie with $g =8$, connected with a rational normal curve
of degree $4$.}
Vestnik MGU, Ser.I: Math. - Mech., No.6 (1985), 99-101.
(in Russian); translated in: Moscow Univ. Math. Bull.

\smallskip

[Ve]  A.Verra, {\it The Prym map has degree two on plane sextics}.
preprint (1991).

\smallskip

[Vo]  C.Voisin, {\it Sur la jacobienne interm\'ediaire du double
solide d'indice deux}.
Duke Math. J., 57:2 (1988), 629-646.

\smallskip

[W]  G.E.Welters, {\it Abel-Jacobi isogenies for certain types
of Fano threefolds}.
Math. Centre Tracts 141, Math. Centrum  Amsterdam (1981).

\vspace{2 cm}

Atanas Iliev

\ Institute of Mathematics

\ Bulgarian Academy of Sciences

\ Acad.G.Bonchev Str., bl.8

\ 1113  Sofia,  Bulgaria

E-mail address:  algebra@bgearn.bitnet

\end{document}